\documentclass[11pt]{article}
\usepackage{amsfonts,color}
\usepackage{amsmath,amssymb,tipa,enumerate,txfonts,mathrsfs,algorithm,float,indentfirst}
\usepackage{bm}
\usepackage{a4wide}
\usepackage{epsfig,multirow,epstopdf}
\usepackage{graphicx}
\usepackage{ragged2e}

\makeatletter
  \newcommand\figcaption{\def\@captype{figure}\caption}
  \newcommand\tabcaption{\def\@captype{table}\caption}
\makeatletter

\newtheorem{lemma}{Lemma}[section]
\newtheorem{theorem}{Theorem}[section]
\newtheorem{remark}{Remark}[section]
\newtheorem{corollary}{Corollary}[section]

\def\be{\begin{equation}}
\def\en{\end{equation}}
\def\beq{\begin{eqnarray}}
\def\eq{\end{eqnarray}}
\def\beqx{\begin{eqnarray*}}
\def\eqx{\end{eqnarray*}}

\def\12{{1\over 2}}

\def\0{{\bf 0}}

\def\prod{{\bf \sqcap}}

\def\Pi{\pi}

\begin{document}
\title{An unconditionally convergent discontinuous Galerkin neural network method activated by plane waves for Helmholtz equation and Maxwell's equations }
\author{
Long Yuan\footnote{College of Mathematics and Systems Science, Shandong University of Science and Technology,
Qingdao 266590, China. (email:\it{sdbjbjsd@163.com}).
}
~~~ Menghui Wu\footnote{College of Mathematics and Systems Science, Shandong University of Science and Technology,
Qingdao 266590, China. (email:\it{1293242185@qq.com}).
}
~~~and~~~
Xiaoqiang Yue\footnote{Corresponding author. National Center for Applied Mathematics in Hunan,
Key Laboratory of Intelligent Computing $\&$ Information Processing of Ministry of Education,
Hunan Key Laboratory for Computation and Simulation in Science and Engineering, Xiangtan University, Xiangtan 411105, China. (email:\it{ yuexq@xtu.edu.cn}).
}
}

\date{\today }
\maketitle

\begin{abstract} In this paper we propose a discontinuous Galerkin plane wave neural network (DGPWNN) method for approximately solving Helmholtz equation and Maxwell's equations. In this method, we define an elliptic-type variational problem as in the plane wave least square method with $h-$refinement and introduce the adaptive construction of recursively augmented discontinuous Galerkin subspaces whose basis functions are realizations of element-wise neural network functions with hp-refinement, where the activation function is chosen as a complex-valued exponential function like the plane wave function. A sequence of basis functions approaching the unit residuals are recursively generated by iteratively solving quasi-maximization problems associated with the underlying residual functionals and the intersection of the closed unit ball and discontinuous plane wave neural network spaces. The convergence results of the DGPWNN method are unconditional only if nonlinear initial parameters for the quasi-maximizers are selected by suitably distributed propagation directions, and established without the assumption on the boundedness of the neural network parameters. Numerical experiments confirm the effectiveness of the proposed method.
\end{abstract}


\vskip 0.1in


\section{Introduction}

The neural network method has emerged as a powerful tool for approximating solutions to partial differential equations (PDEs). To obtain an approximate solution of a PDE using neural networks, the key step is indeed to minimize the PDE residual effectively. Several approaches have been proposed to achieve this \cite{berg, raissi, hexu, Kharazmi, eyu, zang}. While such physics-informed neural networks and variational approaches have enjoyed success in particular cases, the stagnation of relative $L^2$ error of neural network approximation solutions around $0.1-1.0\%$ still exists, no matter how many neurons or layers are used to define the underlying network architecture (see \cite[the second paragraph of Introduction]{AD}). Pushing beyond this limit requires addressing optimization, architecture and sampling challenges. Particularly, {\it continuous} neural networks combined with the Galerkin framework, based on the adaptive construction of a sequence of finite-dimensional subspaces whose basis functions are realizations of a sequence of deep neural networks, have been explored in \cite{AD}. This approach relies on the existing assumption that the bilinear form of the associated variational problem satisfies {\it symmetric, bounded} and {\it coercive} properties. The key aspects of the sequential algorithm are to iteratively improve the accuracy of an approximation and provide a useful indicator for the error in the energy norm that can be used as a criterion for terminating the sequential updates. The convergence results established in \cite{AD} rely on a key assumption that the neural network parameters are bounded. Besides this work, the very recent work \cite{yuanhu} proposed a {\it discontinuous} plane wave neural network (DPWNN), where the desired approximate solution is recursively generated by iteratively solving the quasi-minimization problem associated with the quadratic functional and the sets spanned by element-wise neural network functions with a single hidden layer.

The study of wave propagation is fundamental in various scientific and engineering disciplines, including radar, sonar, acoustics, medical imaging, and seismology. Mathematically, these problems are often described using partial differential equations (PDEs), with the Helmholtz equation and Maxwell's equations being two key models(see \cite{ALZOU, BZZOU, zhao}). The oscillatory nature of solutions in propagating problems (such as wave propagation, quantum mechanics, or time-dependent PDEs) indeed poses significant challenges in numerical simulations. Recent work \cite{MP, chung, yuan2,yuanyue} has led to the development of algorithms that incorporate the highly oscillatory nature of the solutions into the basis functions. In particular, plane wave methods, which are non-polynomial finite element methods for time-harmonic wave problems, have been designed in the last years in order to reduce the computational cost, with respect to more classical polynomial-based methods. Indeed, the Plane Wave Discontinuous Galerkin (PWDG) method \cite{Git, ref21,pwdg,yuan2, yuanyue} and the Plane Wave Least Squares (PWLS) method \cite{ref12,hy2,hy3,peng,peng2} are two prominent numerical techniques that leverage plane wave bases for solving wave propagation problems, particularly in the context of Helmholtz equations (time-harmonic wave problems). All these methods differ from the traditional finite-element method (FEM) and the boundary-element method (BEM) in the sense that the plane wave functions are chosen as basis functions, satisfying the governing differential equation, to generate the finite dimensional subspaces for a standard discontinuous Galerkin approximation of the variational equation.

We would like to point out that, a key advantage of the plane wave methods over Lagrange finite elements when discretizing the Helmholtz equation and time-harmonic Maxwell equations \cite{ref21,pwdg,hy2,HMPsur,peng,peng2} is that: plane waves are tailored to the PDE's structure, allowing them to represent the solution more compactly than generic basis functions, and thereby relatively smaller number of degrees of freedom are enough for the same accuracy, provided boundary conditions and potential ill-conditioning are addressed carefully. In particular, an advantage of the PWLS method over the other plane wave methods is that the sesquilinear form associated with the PWLS method is Hermitian positive definite, so it is natural to combine the PWLS method with the Galerkin neural networks for wave equations. In this paper we define a sequence of sets of {\it discontinuous} plane wave neural networks, where the activation function on each element is chosen as $e^{i\omega x}$ ($\omega$ denotes the wavenumber) and a single hidden layer is considered, and apply it to the discretization of the minimization problem for the time-harmonic Helmholtz equation and Maxwell equations. However, differently from the standard PWLS method, for the new method the plane wave directions are not fixed and can be corrected in the iterative process of searching an approximate maximizer of the proposed residuals.

We introduce the adaptive construction of recursively augmented discontinuous Galerkin subspaces whose basis functions are realizations of element-wise neural network functions with $hp$-refinement. A sequence of basis functions approaching the unit residuals are recursively generated by iteratively solving {\it quasi-maximization} problems associated with the underlying residual functionals defined by the former {\it quasi-maximizer} and the intersection of the closed unit ball and discontinuous plane wave neural network spaces. For convenience, we call the method as DGPWNN method. The convergence results of the DGPWNN method are unconditionally established without the assumption on the boundedness of the neural network parameters. The performance of the algorithms will confirm the effectiveness of the proposed method. While the proposed method described above falls into the class of Galerkin neural networks, it differs from the Galerkin neural network method introduced in \cite{AD} both in algorithm design and theoretical analysis. Particularly, new activation function-the plane wave function is defined element-wisely, and plane wave direction angles and activation coefficients are alternatively computed by iterative algorithms. Such choice of activation functions can enhance the accuracies of the resulting approximate solutions for the underlying models and accelerate the convergence of the iterative algorithm.

It is worth noting that, if we augment the discretization subspace without training directions procedures for the approximate maximizers in Algorithm 2, the proposed DGPWNN method degenerates into the same greedy residual-enrichment algorithm but with fixed (equispaced) directions. Since such a fixed-direction space is contained in the nonlinear dictionary, the greedy discontinuous Galerkin method has notable significance when sufficiently accurate quasi-maximizers are found. The quasi-maximizers found by the optimization algorithm always ensures the unconditional convergence of the training algorithm and high accuracies, only if nonlinear initial parameters are selected according to suitably distributed propagation directions. {\it To our knowledge, it is the first unconditional convergence result of the Galerkin neural networks.}

Compared with the very recent work \cite{yuanhu}, the proposed method has a completely different way of augmenting the discontinuous Galerkin subspaces. The latter is to seek the approximate maximizer $\varphi_i^{NN}$ of the new defined residual $ \langle r(u_{i-1}), v\rangle$ in $V_{n_i}^{\sigma}\cap B$, and use it to further augment the discontinuous Galerkin subspace
 $\Phi_i^{NN}$. While the former is to seek the approximate quasi-minimizer $\xi_r^{\theta}$ of the original loss functional $J(u_{i-1} + \cdot)$ in the subset $V_{n_i}({\mathcal T}_h)$, and use it to directly get a better approximation $u_i=u_{i-1}+\xi_r^{\theta}$ to $u$. Besides, for the DGPWNN method, the equivalent computable quantity $\eta(u_{i-1}, \varphi_i^{NN})$ to the true error $|||u-u_{i-1}|||$ in the energy norm is additionally obtained, which can be used as a stopping criterion for the proposed recursive solution procedure. Moreover, for the DGPWNN method, the condition number of the associated discrete system induced by basis functions $\{ \varphi_j^{NN}\}_{j=1}^i$ in the augmented discontinuous Galerkin subspace $\Phi_i^{NN}$ can be bounded by a constant slightly greater than 1.
Table \ref{COMmethod} lists the main algorithmic differences among the proposed DGPWNN method, DPWNN method \cite{yuanhu}, GNN method \cite{AD} and plane wave methods \cite{ref21,pwdg,hy2,HMPsur,peng,peng2}.

\vskip 0.1in
\begin{center}
       \tabcaption{Comparison among different methods.}
\label{COMmethod}
       \vskip 0.1in
\begin{tabular}{|c|c|c|c|} \hline
   Method  & CG / DG & Training Parameters &  Augmenting Galerkin space by  \\ \hline
   DGPWNN &   DG &  Yes & Residual \\ \hline
    DPWNN &   DG &  Yes & Loss functional \\ \hline
    GNN &   CG &  Yes & Residual  \\ \hline
    PWM &   DG &  No & Pre-defined  \\ \hline
   \end{tabular}
     \end{center}

The paper is organized as follows: In section 2, we define a minimization problem for the underlying Helmholtz equation. In section 3, we describe discontinuous plane wave neural networks. In section 4, we describe a discretization for the minimization problem, and design iterative algorithms to generate approximate solutions of the discrete  quasi-maximization problem, and give error estimates of the approximate solutions. In section 5, we apply the proposed discontinuous Galerkin plane wave neural network to discretize the time-harmonic Maxwell equations. Finally, in section 6, we report some numerical results to confirm the
effectiveness of the proposed method.

\section{A minimization method for underlying Helmholtz model} \setcounter{equation}{0}

Let $\Omega \subset \mathbb{R}^d~(d=2, 3)$ be a bounded polygonal/polyhedral domain with boundary $\gamma=\partial\Omega$. Consider Helmholtz equation which is formalized by
\begin{eqnarray}
\left\{\begin{array}{ll} -\Delta u-\omega^2u=0& \text{in}\quad
\Omega\\
(\partial_\text{\bf n} + \text{i}\omega)u=g& \text{on}\quad\gamma=\partial\Omega
\end{array}\right.
\label{helm1}
\end{eqnarray}
where $g \in L^2(\partial\Omega)$. The outer normal derivative is referred to by $\partial_{\bf n}$ and the angular frequency by $\omega$.


Let \( {\cal T}_h\) be a finite element partition of \(\Omega\), with possible hanging nodes. \(h\) is the meshwidth of the triangulation. We will define $\Gamma_{lj}=\partial\Omega_l\bigcap\partial\Omega_j$ for $l\not=j$ and $\Omega_l, \Omega_j\in {\cal T}_h$, and set $\gamma_k=\overline{\Omega}_k\bigcap\partial\Omega ~(k=1,\ldots,N).$ Then $\gamma=\bigcup^N_{k=1}\gamma_k$.

The space $V(\Omega_k)$ typically denotes a function space consisting of solutions to the homogeneous Helmholtz equation on each element $\Omega_k$:
\be
 V(\Omega_k) = \bigg\{v_k \in H^1(\Omega_k); -\Delta v_k - \omega^2 v_k = 0 \bigg\}.
 \en
Define
\be
V({\cal T}_h) = \{ v\in L^2(\Omega): ~v|_{\Omega_k}\in V(\Omega_k) ~\forall \Omega_k\in {\cal T}_h  \}.
\en
Set $v|_{\Omega_k}=v_k$ $(k=1,\cdots,N$). Define the functional
\begin{eqnarray} \label{funct}
 J(v) &=& \sum_{k=1}^N\int _{\gamma_k} |(\partial_{\bf
n}+i\omega )v_k-g|^2 ds  \cr
& + & \sum_{j\not=k} \bigg( \alpha \int_{\Gamma_{kj}}
|v_k-v_j|^2 ds + \beta \int_{\Gamma_{kj}} | (\partial_{\text{\bf n}_k}
v_k+\partial_{\text{\bf n}_j} v_j ) |^2ds \bigg), ~~\forall v\in V({\mathcal T}_h),
\end{eqnarray}
where $\alpha$ and $\beta$ are two given relaxed factors. 

Consider the minimization problem: find $u\in {V}({\cal T}_h)$ such that
\begin{equation} \label{minmum}
J({u})=\min\limits_{{v}\in {V}({\cal T}_h)}~J(v).
\end{equation}

The variational problem related to the minimization problem (\ref{minmum}) typically involves finding a function $u \in V({\mathcal T}_h)$ such that (refer to \cite{hy}):

\begin{eqnarray}
 &&\sum_{k=1}^N\int _{\gamma_k}\big((\partial_{\bf
n}+i\omega )u_k-g\big)\cdot\overline{(\partial_{\bf n}+i\omega) v_k}ds+
 \sum_{j\not=k} \bigg(\alpha\int_{\Gamma_{kj}}
(u_k-u_j)\cdot\overline{(v_k-v_j)}ds \cr
&&\quad\quad+\beta\int_{\Gamma_{kj}}(\partial_{\text{\bf n}_k}
u_k+\partial_{\text{\bf n}_j} u_j
)\cdot\overline{(\partial_{\text{\bf n}_k} v_k+\partial_{\text{\bf
n}_j} v_j )} ds\bigg)=0, ~~\forall v\in V({\mathcal T}_h).
\label{va1}
\end{eqnarray}
Define the sesquilinear form $a(\cdot, \cdot)$ by
\beq
& a(u,v) = \sum_{k=1}^N\int _{\gamma_k} (\partial_{\bf
n}+i\omega )u_k \cdot\overline{(\partial_{\bf n}+i\omega) v_k}ds
+ \sum_{j\not=k} \bigg(\alpha\int_{\Gamma_{kj}}
(u_k-u_j)\cdot\overline{(v_k-v_j)} ds  \cr
& + \beta\int_{\Gamma_{kj}}(\partial_{\text{\bf n}_k}
u_k +\partial_{\text{\bf n}_j} u_j ) \cdot \overline{(\partial_{\text{\bf n}_k} v_k+\partial_{\text{\bf n}_j} v_j )} ds \bigg), ~\forall v\in V({\mathcal T}_h),
\eq
and semilinear form $L(v)$ by
\be
L(v)=\sum_{k=1}^N\int _{\gamma_k} g ~\overline{(\partial_{\bf
n}+i\omega )v_k} ds, ~\forall v\in V({\mathcal T}_h).
\en
 Then (\ref{va1}) can be written as the following variational formulation:
 \begin{eqnarray}\label{pwlsvar}
\left\{\begin{array}{ll} \text{Find} \ u\in V({\mathcal T}_h) ~~s.t. \\
a(u,v)=L(v) \quad \forall v\in V({\mathcal T}_h).
\end{array}\right.
\end{eqnarray}

We can see that $a(v,v)$ is a norm on $V({\mathcal T}_h)$ (see \cite{hy}), and denote it by $||| v |||^2$. Using the Cauchy-Schwarz inequality, it can be proved that $|a(u,v)|\leq ||| u |||~||| v |||$ and $|L(v)| \leq ||g||_{L^2(\gamma)} ~||| v |||, ~\forall u,v \in V({\mathcal T}_h)$. Thus, a direct application of the Risez representation theorem \cite{FO} gives rise to the existence of a unique $u\in V({\mathcal T}_h)$ satisfying (\ref{pwlsvar}).



\section{Element-wise plane wave neural networks} \setcounter{equation}{0}

An element-wise neural network consists of a single hidden layer of $n\in \mathbb{R}$ neurons on each element $\Omega_k\in {\cal T}_h$, defining a function $\varphi^{NN}: \mathbb{R}^d \rightarrow \mathbb{C}$ as follows:
\be \label{nn_app}
\varphi^{NN}(x; \theta)|_{\Omega_k} = \sum_{j=1}^n c_j^{(k)} \sigma(W_j^{(k)} \cdot {\bf x} + b_j^{(k)}), ~\forall \Omega_k \in {\cal T}_h,
\en
where $n$ is called the width of the network. Nonlinear coefficients $W_j^{(k)}\in \mathbb{R}^d, b_j^{(k)}\in \mathbb{R}$ are defined element-wisely, $c_j^{(k)} \in \mathbb{C}$ are element-wise linear coefficients and $\sigma: \mathbb{R}\rightarrow \mathbb{C}$ is a bounded, elementwise smooth activation function. For simplicity, we shall use the elementwise notation $W|_{\Omega_k}=[W^{(k)}_1, \cdots, W^{(k)}_n]\in \mathbb{R}^{d\times n}$ referred as weights, $b|_{\Omega_k}=[b^{(k)}_1, \cdots, b^{(k)}_n]^T\in \mathbb{R}^{n}$ referred as biases, and $c|_{\Omega_k}=[c^{(k)}_1, \cdots, c^{(k)}_n]^T\in \mathbb{C}^n$. A set of nonlinear and linear parameters over a collection of all elements is denoted by $\theta=\{W; b; c\}$.


As pointed out in Introduction, the plane wave method, which is based on the Trefftz approximation space made of plane wave basis functions, has a related advantage over Lagrange finite elements for discretization of the time-harmonic wave equations \cite{pwdg}. Thus we choose plane wave complex-valued functions denoted by $\text{e}^{\text{i}\omega P\cdot x}~(P\in \mathbb{R}^d, |P|=1)$ as the used activation function $\sigma(x)$.  Notice that $\text{e}^{\text{i}\omega b^{(k)}_j}$ in (\ref{nn_app}) can be equivalently learned by the linear coefficient $c^{(k)}_j$, we can ignore $b^{(k)}_j$ in (\ref{nn_app}) in the situation to obtain an equivalent neural network on each element $\Omega_k\in {\cal T}_h$
\be \label{nn_app2_general}
\varphi^{NN}(x; \theta)|_{\Omega_k} = \sum_{j=1}^n c^{(k)}_j \text{e}^{\text{i} \omega W^{(k)}_j \cdot x}, ~~|W^{(k)}_j| = 1,
\en
with the training parameters $\theta^{(k)}=\{D^{(k)}; c^{(k)}\}$ defined below. Denote by $V_n^{\sigma}$ a discontinuous plane wave neural network space consisting of all functions of the form (\ref{nn_app2_general}), which is defined by
\be \label{vsign_general}
V_n^{\sigma}({\mathcal T}_h) := \bigg\{v: ~v|_{\Omega_k}=\sum_{j=1}^n c^{(k)}_j \text{e}^{\text{i} \omega W^{(k)}_j \cdot x}, ~c^{(k)}_j \in \mathbb{C}, ~W^{(k)}_j\in \mathbb{R}^d, ~|W^{(k)}_j| = 1, ~\forall \Omega_k\in {\mathcal T}_h  \bigg\}.
\en
Note that $V_n^{\sigma}$ is not a linear space owing to the existence of nonlinear coefficients, namely, the set $W$ of plane wave propagation directions.

For the 2D case, element-wise plane wave propagation directions $W^{(k)}_{l}$ ($l=1,\cdots,n$) are defined by
\be \label{nn_app2}
W^{(k)}_l = [cos(d^{(k)}_l), sin(d^{(k)}_l)]^T, ~d^{(k)}_l\in (-\pi, \pi],
\en
 where a set of plane wave propagation angles and linear parameters defined on every element are represented in the form $D|_{\Omega_k}=[d^{(k)}_1, \cdots, d^{(k)}_n]^T$, $c|_{\Omega_k}=[c^{(k)}_1, \cdots, c^{(k)}_n]^T\in \mathbb{C}^n$, respectively.

 For the 3D case, plane wave propagation directions $W^{(k)}_{l}$ ($l=1,\cdots,n$) based on spherical coordinates are defined by
\begin{equation} \label{3ddirec}
W^{(k)}_l = \left ( {\begin{array}{c}
sin(\zeta_m^{(k)}) cos(\vartheta_t^{(k)})  \\
sin(\zeta_m^{(k)}) sin(\vartheta_t^{(k)})  \\
cos(\zeta_m^{(k)})
\end{array}}
\right ),  ~\zeta_m^{(k)}\in [0, \pi],~m =1,\cdots, m^{\ast}; ~\vartheta_t^{(k)}\in (-\pi, \pi], ~t=1,\cdots, t^{\ast},
\end{equation}
where $m^{\ast}$ and $t^{\ast}$ are positive integers, which satisfy $t^{\ast}\approx 2m^{\ast}$ and $n=m^{\ast}t^{\ast}$. We will use the element-wise notation $\zeta|_{\Omega_k}=\zeta^{(k)}=[\zeta_1^{(k)}, \cdots, \zeta_{m^{\ast}}^{(k)}]^T\in \mathbb{R}^{m^{\ast}}$, $\vartheta|_{\Omega_k}=\vartheta^{(k)}=[\vartheta_1^{(k)}, \cdots, \vartheta_{t^{\ast}}^{(k)}]^T \in \mathbb{R}^{t^{\ast}}$, and $c|_{\Omega_k}=c^{(k)}=[c_1^{(k)}, \cdots, c_{n}^{(k)}]^T\in \mathbb{C}^{n}$. Particularly, a set of nonlinear parameters, namely, a set of plane wave propagation angles defined on every element, is collectively denoted by $D=\{\zeta; \vartheta\}$.


\subsection{A universal approximation of discontinuous neural networks}

As previously described, the plane wave basis functions can generate the associated plane wave neural network. We first give the existing approximation results about $hp-$approximation estimates for homogeneous Helmholtz solutions in $H^s(\Omega_k)$ using plane wave neural networks, based on \cite[Lemma 4.4]{hy3} and \cite[Corollary 5.5]{mhp}.


 Let \(s\) and \(m\) be given positive integers
satisfying \(m\geq 2(s-1)+1=2s-1\). Let the
number \(n\) of plane wave propagation directions be chosen as
\(n=2m+1\) in 2D and \(n=(m+1)^2\) in 3D, respectively.

\begin{lemma}  \label{elementappro}
Let $2\leq s\leq {m+1\over 2}$ with a sufficiently large $m$.
Let $u\in H^{s}(\Omega_k)$ be a solution of the homogeneous Helmholtz
equation for each element $\Omega_k$. Then, there exists $[c_1^{(k)},c_2^{(k)},\cdots,c_n^{(k)}] \in \mathbb{C}^n$ and $[W_1^{(k)}, W_2^{(k)},\cdots, W_n^{(k)}]\in \mathbb{R}^{d\times n}$,  such that
\be \label{approx}
||u - \sum_{l=1}^n c_l^{(k)} \text{e}^{\text{i} \omega W_l^{(k)} \cdot{\bf x}} ||_{j,\Omega_k}\leq C h^{s-j}m^{-\lambda(s-j-\varepsilon)}||u||_{s,\omega,\Omega_k}~~~~(0\leq j\leq s; k=1,\cdots,N),
\en
where $\lambda>0$ is a constant depending only on the shape of the elements (in particular, $\lambda=1$ for the
 case of two dimensions), $\varepsilon=\varepsilon(m)>0$ satisfies $\varepsilon(m)\rightarrow 0$ when $m\rightarrow\infty$, and the $\omega-$weighted Sobolev norm is defined by
$||v||_{s,\omega,D} = \sum_{j=0}^{s}\omega^{2(s-j)}|v|_{j,D}^2.$
\end{lemma}

\begin{remark}\label{univerappro}
The universal approximation estimates (\ref{approx}) are classical plane wave approximation estimates, which are usually stated for fixed, prescribed, or suitably distributed propagation directions. They justify approximation in a classical plane wave space with appropriate directions.
\end{remark}

Next, the following universal approximation property of discontinuous plane wave neural networks can be obtained, which can also be achieved by the classical approximate procedure \cite{HO}.




\begin{theorem} \label{aninappro}
Suppose that $s \geq 2$, and that $\Omega\subset \mathbb{R}^d$ is compact, 
 then $V_n^{\sigma}({\mathcal T}_h)$ is dense in the Sobolev space
\be\nonumber
H^{s}({\cal T}_h) := \bigg\{v\in L^2(\Omega): v \in  V(\Omega_k) ~\text{and}~D^{\alpha}v\in  L^2(\Omega_k)~~\forall |\alpha|\leq s, k=1,\cdots,N   \bigg\}.
\en
\end{theorem}

In particular, the universal approximation theorem means that, given a function
$f \in  H^{s}({\cal T}_h)$ and $\tau  > 0$, there exists $n(\tau , f) \in  \mathbb{N}$, and $\tilde{f} \in  V^{\sigma}_{n(\tau ,f )}({\cal T}_{h})$ such that $| | f-\tilde{f}| |_{H^{s}({\cal T}_h)} < \tau$.

\section{A discontinuous Galerkin plane wave neural network method}\label{helmdis} \setcounter{equation}{0}

In this section, we present an iterative algorithm to generate element-wise plane wave neural network solutions that approximately satisfy the minimization problem (\ref{minmum}).



\subsection{An overview of DGPWNN algorithm}


Suppose $u_{i-1}\in { V}({\cal T}_h)$ is any given approximation to the solution of (\ref{pwlsvar}). Define a new residual functional $r(u_{i-1}): V({\cal T}_h) \rightarrow \mathbb{R}$ given by the formula
\be \label{resiformulaiter}
\langle r(u_{i-1}), v \rangle = \text{Re} \{ L(v) - a(u_{i-1}, v) \} = \text{Re} \{ a(u-u_{i-1}, v) \},
\en
where $\text{Re} \{ \cdot \}$ denotes the real part of quantity $\cdot$.

Note that owing to the injection of complex Euclidean space in the considered model, the current residual $\langle r(u_{i-1}), v \rangle$ is the real part of the pre-residue $\langle \xi(u_{i-1}), v\rangle=L(v) - a(u_{i-1}, v)$, which is completely different from the definition of  $\langle r(u_{i-1}), v \rangle$ in \cite{AD}. 

Set
\be \label{defivarphi2}
\varphi_i=(u-u_{i-1})/|||u-u_{i-1}|||.
\en
Then by \cite[Proposition 2.3]{AD}, we can obtain
\be \label{conmax}
\langle r(u_{i-1}), \varphi_i \rangle = \max\limits_{v\in B} ~\langle r(u_{i-1}), v\rangle = |||u-u_{i-1}|||,
\en
 where $B$ is the closed unit ball in $V({\cal T}_h)$. It shows that the maximizer $\varphi_i \in  B$ of the residual is proportional to the error $e_{i-1} = u -u_{i-1}$.

Unfortunately, we can not compute the maximizer $\varphi_i$ explicitly, which would be equivalent to solving the variational problem (\ref{pwlsvar}) exactly. Instead, we compute an approximate maximizer $\varphi_i^{NN} \approx \varphi_i$ in an appropriate way.

For a positive integer $n_i$ defined later,  when aiming to maximize $\langle r(u_{i-1}), \cdot \rangle$ in the subset $V_{n_i}^{\sigma}\cap B$, a significant complication is that it  is not a linear subspace and it is not closed (topologically) in $B$. Hence, even though $\langle r(u_{i-1}), \cdot \rangle$ have an supremum on $V_{n_i}^{\sigma}\cap B$, there may not be a maximizer in $V_{n_i}^{\sigma}\cap B$. Therefore, one should not aim to completely maximize $\langle r(u_{i-1}), \cdot \rangle$, but instead use a relaxed notion of quasi-maximization as used by \cite{Karniadakis} (for which the existence of an infimum implies the existence of a quasi-minimizer):

 {\bf Definition 4.1.} ~~{\it Let $\epsilon_i$ ($0<\epsilon_i<1$) be a sufficiently small positive parameter, and let $V_{n_i}^{\sigma}\cap B$ be a subset of $B$. A function $\varphi_i^{NN} \in V_{n_i}^{\sigma}\cap B$ is said to be a quasi-maximizer of
 $\langle r(u_{i-1}), \cdot \rangle$ associated with $\epsilon_i$ if the following inequality holds:
 \be \label{iteradisminmum}
 \langle r(u_{i-1}), \varphi_i^{NN} \rangle \geq 
 \sup\limits_{v_i^{NN}\in V_{n_i}^{\sigma}\cap B} \langle r(u_{i-1}), v_i^{NN} \rangle - \epsilon_i |||u-u_{i-1}|||.
 \en
 }
It is clear that the value of the functional $\langle r(u_{i-1}), \cdot \rangle$ at the quasi-maximizer $\varphi_i^{NN}$ is a good approximation of the supremum.
Besides, the quantity $\eta$ is defined by
\be \label{etaquan}
\eta(u_{i-1},v) := \langle r(u_{i-1}), v\rangle /|||v|||, \quad \forall v\in V_{n_i}^{\sigma}.
\en
Then it is obvious that
\be \label{etaequa}
\eta(u_{i-1},\varphi_i) =  \langle r(u_{i-1}), \varphi_i \rangle = \max\limits_{v\in V({\cal T}_h)} \eta(u_{i-1},v), \quad
\eta(u_{i-1}, \varphi_i^{NN} ) \geq  \sup\limits_{v_i^{NN}\in V^{\sigma}_{n_i}} \langle \eta(u_{i-1}, v_i^{NN} ) - \epsilon_i.
\en

After solving $\varphi_i^{NN}$, we use it to further augment the discretization subspace by setting $\Phi_i^{NN}=\text{span}\{ \varphi_1^{NN}, \cdots, \varphi_i^{NN}\}$.

Then a better approximation $u_i$ to $u$ can be computed by solving the discrete variational problem:
 \begin{eqnarray}\label{dgappro}
\left\{\begin{array}{ll} \text{Find} \ u_i\in \Phi_i^{NN}~~s.t. \\
a(u_i,v) = L(v) \quad \forall v\in \Phi_i^{NN}.
\end{array}\right.
\end{eqnarray}

Algorithm 1 summarizes this recursive approach.

\begin{algorithm}[H]
\caption{Discontinuous Galerkin plane wave neural networks}

\quad Set $i=1$ and $u_0=0$.

\quad $(\varphi_1^{NN}, \eta) \leftarrow $ AugmentBasis($u_0$).

\quad {\bf while } $\eta > tol$ do

\quad \quad  $u_{i} \leftarrow $ dGSolve($\{\varphi_1^{NN}, \cdots, \varphi_{i}^{NN}\}$, $a, L$).

\quad \quad $(\varphi_{i+1}^{NN}, \eta) \leftarrow $ AugmentBasis($u_i$).

\quad \quad Set $i \leftarrow i+1$.

\quad  {\bf end }

\quad Return $N=i-1, ~u_N,$ and $\{\varphi_{j}^{NN}\}_{j=1}^N$.
\end{algorithm}

Here, Function dGSolve($\{\varphi_1^{NN}, \cdots, \varphi_{i}^{NN}\}$, $a, L$) returns the projection of $u$ onto the space $\Phi_{i}^{NN}$ with respect to the variational equation (\ref{dgappro}) associated with the bilinear $a$ and data $L$.

The training procedure for the approximate maximizer $\varphi_i^{NN}$ of  the weak residual functional $r(u_{i-1})$ is summarized in Algorithm 2.


\begin{algorithm}[H]
\caption{Plane wave basis function generation.} 

\quad {\bf Function} AugmentBasis($u_{i-1}$):

\quad\quad  Initialize hidden parameters $D^{(i)}\in \mathbb{R}^{n_i}$ by (\ref{nn_app2}) or (\ref{3ddirec}).

\quad \quad Compute corresponding activation coefficients:

\quad \quad $c^{(i)} \leftarrow$ dGLSQ-R$(D^{(i)}, \sigma, a, L(\cdot)-a(u_{i-1}, \cdot))$.

\quad \quad {\bf for} each training epoch {\bf do}

\quad \quad  \quad Compute the update propagation wave directions for $D^{(i)}$ by an randomly shuffled Adam optimizer \cite{luo}.

\quad \quad  \quad $\{\zeta^{(i)}_m\}$ corrections: correct associated angles $\{\zeta^{(i)}_m\}$ that satisfy $|sin(\zeta^{(i)}_m)|< 1.0e-6$ by adding suitable disturbance term $\varepsilon_{disturb}$ to guarantee the linear independence of the resulting basis functions $\{e^{\text{i}\omega W^{(i)}_{l}\cdot x}\}$ for different $\vartheta_t^{(i)}$, and generally set $\varepsilon_{disturb}=\frac{\pi}{1000m^{\ast}_{i}}$ and $\zeta^{(i)}_m = \zeta^{(i)}_m + \varepsilon_{disturb}$.  

\quad \quad  \quad Compute corresponding activation coefficients:

\quad \quad  \quad $c^{(i)} \leftarrow$ dGLSQ-R$(D^{(i)}, \sigma, a, L(\cdot)-a(u_{i-1}, \cdot))$.

\quad \quad  {\bf end for}

\quad \quad  Return $\varphi_i^{NN}=\frac{v(\theta^{(i)})}{|||v(\theta^{(i)})|||}$ and $\eta(u_{i-1},\varphi_i^{NN})$.

\quad  {\bf end function}

\end{algorithm}

Here Function dGLSQ-R$(D^{(i)}, \sigma, a, L(\cdot)-a(u_{i-1}, \cdot))$ computes the expansion coefficients of the projection of the error $u-u_{i-1}$ onto the subspace $V_{n_i}^{\sigma}({\mathcal T}_h)$ with respect to $a$ and the pre-residue $\langle \xi(u_{i-1}), v\rangle = L(v) - a(u_{i-1}, v)$ provided that hidden parameters $D^{(i)}$ are given.

\begin{remark}
The ``While" loop in Algorithm 1 is called an outer iteration, which expands the discrete space $\Phi^{NN}$ by finding a quasi-maximizer $\varphi^{NN}$ and solves the least squares problem (\ref{dgappro}) on it to find a better numerical approximation solution. Further, the training procedure for searching the approximate maximizer $\varphi^{NN}$ of the weak residual functional $r(u_{i-1})$ in Algorithm 2 is called an inner iteration, which appears to be alternate between: a fixed directions linear least-squares/Galerkin coefficient solve, and, the direction update with nonlinear optimization step. Thus the neural-network component is primarily a nonlinear basis-search mechanism, not a neural solver map from data to solutions.
\end{remark}

\begin{remark}\label{diffongrad}
The method described in this section is different from that developed in \cite{AD}: the former is naturally derived from the quasi-maximization problem (\ref{iteradisminmum}), but the latter was derived from a maximization problem. Furthermore, we would like to emphasize that, due to the injection of complex Euclidean space in the considered model and plane wave activation functions, the residual loss functional $r(u_{i-1})$ defined by (\ref{resiformulaiter}) is the real part of the pre-residue $\langle \xi(u_{i-1}), v\rangle =L(v) - a(u_{i-1}, v)$, which is a visible difference from that of \cite{AD}.
 \end{remark}

\begin{remark}
The generation of the sequence of discretized discontinuous Galerkin subspaces is completely from that of \cite{AD}.  The discontinuous Galerkin subspace $\Phi_i^{NN}$ is augmented by the current approximate maximizer $\varphi_i^{NN}$ of the new residual $\langle r(u_{i-1}), v \rangle$. While each discontinuous basis function $\varphi_i^{NN}$ may be expressed as the realization of the {\it discontinuous} Galerkin plane wave neural network space $V_{n_i}^{\sigma}$ with a single hidden layer by learning the quasi-maximization representation of the weak residual functional $r(u_{i-1}): V({\cal T}_h) \rightarrow \mathbb{R}$.
\end{remark}

\begin{remark}\label{alterupdate}
It is worth noting that, the alternating updating of plane wave direction angles and activation coefficients in Algorithm 2 is necessary. On the contrary, if a standard gradient-based Adam optimizer \cite{luo} with full-batch mode is used to update the plane wave direction angles and activation coefficients as a whole, the existing numerical results indicate that the DGPWNN method is likely to diverge.
\end{remark}

\subsection{Implementation on Function dGLSQ-R}

Provided that the nonlinear parameters ${D^{(i)}}$ are fixed, we can form the discontinuous Galerkin plane wave subspace $V_{n_i}^{\sigma}({\mathcal T}_h)$ consisting of $N\times n_i$ neural network basis functions $\{\psi_{kj}\}$, which satisfy
\begin{eqnarray}
\psi_{kj}(x)=\left\{\begin{array}{ll}
 e^{\text{i} \omega W_j\cdot x},~~x\in\Omega_k,\\
 0,~~x\in\Omega_l~~\mbox{satisfying}~~l\neq k
\end{array}\right.~~(k,l=1,\cdots,N;~j=1,\cdots, n_i).
\label{neubasis}
\end{eqnarray}
 Then a discontinuous Galerkin PWLS approximation associated with (\ref{iteradisminmum}) can be calculated from $V_{n_i}^{\sigma}({\mathcal T}_h)$ by choosing the linear parameters $c^{(i)}$ as follows:

\begin{eqnarray}
\left\{\begin{array}{ll}  {\cal A}c^{(i)} = b, \\
{\cal A}^{k,l}_{r,j} = a(\psi_{lj},\psi_{kr}), ~~b_k^r=L(\psi_{kr}) - a(u_{i-1},\psi_{kr}).
\end{array}\right.
\label{galerlsq}
\end{eqnarray}
This procedure defines a function $c^{(i)}=$dGLSQ-R$(D^{(i)}, \sigma, a, L(\cdot)-a(u_{i-1}, \cdot))$ which computes the expansion coefficients of the projection of $u-u_{i-1}$ onto the subspace $V_{n_i}^{\sigma}({\mathcal T}_h)$ with respect to the sesquilinear form $a(\cdot,\cdot)$ and the pre-residue righthand $\langle \xi(u_{i-1}), v\rangle=L(v) - a(u_{i-1}, v)$.

Once the training epoch in Algorithm 2 is completed, we can assemble the ultimate plane wave basis functions $\{\psi_{kj}\}$ and the corresponding linear coefficients $c^{(i)}$ in a proper way, and finally get the current approximate maximizer $\varphi_i^{NN}$ of the residual functional $\langle r(u_{i-1}), \cdot \rangle$ on the set $V_{n_i}^{\sigma}\cap B$.

\subsection{Calculation of the training parameters $D$ for DGPWNN} \label{hlempara}

To compute optimal values of the objective function (\ref{iteradisminmum}) associated with discontinuous Galerkin neural networks, the modified gradient-based Adam optimizer \cite{luo} with full-batch mode is employed to update the hidden parameters ${D}$. On the first training iteration, hidden parameters ${D}$ are initialized uniformly (details specified in numerical experiments). For subsequent iterations, ${D}$ is updated using a randomly shuffled Adam optimizer. Particularly, only one optimizer step is performed per update, and the loss function is computed over $N_G$ quadrature nodes, which serves as the training dataset. Besides, we choose $\eta=\frac{\eta_{\omega,i}}{\sqrt{N_G}}$ as a learning rate scheduler to reduce the learning rate during the process. In particular, $\eta_{\omega,i}$ is the learning rate of decay depending on the frequency $\omega$ and $i$. In our approach, for a fixed $\omega$, we only decrease the learning rate of decay between each training epoch to compensate for the increasing widths $n_i$.

As aforementioned in Algorithm 2, 
 we would like to emphasize the computation of derivative of $\eta(u_{i-1},v)$ with respect to the training parameters $D$, which will be used to update the hidden parameters.



By the direct reduction, we obtain
\be
\nabla_D \eta(u_{i-1},v) = \frac{\nabla_D \langle r(u_{i-1}), v \rangle }{|||v|||} + \langle r(u_{i-1}), v \rangle ~\nabla_D \big(\frac{1}{|||v|||}\big),
\en
where
\be
\nabla_D \langle r(u_{i-1}), v \rangle = \text{Re}\big\{ \langle \xi(u_{i-1}), \nabla_D v\rangle \big\} 
\en
and
\be
\nabla_D \big(\frac{1}{|||v|||}\big) = -|||v|||^{-3}\text{Re}\{a(v, \nabla_D v)\}.
\en

Exactly, taking as an example
$v =  \sum_{l=1}^n c^{(k)}_l \text{e}^{\text{i} \omega W^{(k)}_l \cdot x}, ~x\in \Omega_k \subset {\cal T}_h$ in two dimensions, we get $\nabla_{D^{(k)}} v =(\partial_{d^{(k)}_1}v, \cdots, \partial_{d^{(k)}_n}v)^T$, where
\be \nonumber
\partial_{d^{(k)}_m}v = i\omega  c^{(k)}_m \big(-x_1 sin(d^{(k)}_m) + x_2 cos(d^{(k)}_m)\big)~e^{\text{i}\omega W^{(k)}_{m}\cdot x}.
\en
In three dimensions, taking $v = \sum_{m=1}^{m^{\ast}}\sum_{t=1}^{t^{\ast}}  c_l^{(k)} ~e^{\text{i}\omega
{\bf d}_{l}^{(k)}\cdot {\bf x}}, ~l=(m-1)k^{\ast}+t, ~{\bf x}=(x_1, x_2, x_3)^T \in \Omega_k \in {\cal T}_h$, we get $\nabla_{D^{(k)}} v = (\partial_{\vartheta^{(k)}_1}v, \cdots, \partial_{\vartheta^{(k)}_{t^{\ast}}}v, \partial_{\zeta^{(k)}_1}v, \cdots, \partial_{\zeta^{(k)}_{m^{\ast}}}v)^T$, where
\be \nonumber
\partial_{\vartheta_t^{(k)}}v = i\omega \sum_{m=1}^{m^{\ast}} c_l^{(k)} \big(-x_1 sin(\zeta_m^{(k)})sin(\vartheta_t^{(k)}) + x_2 sin(\zeta_m^{(k)})cos(\vartheta_t^{(k)})\big)~e^{\text{i}\omega
W_{l}^{(k)}\cdot {\bf {\bf x}}},
\en
and
\be\nonumber
\partial_{\zeta_m^{(k)}}v = i\omega \sum_{t=1}^{t^{\ast}} c_l^{(k)}  \big(x_1 cos(\zeta_m^{(k)})cos(\vartheta_t^{(k)}) + x_2 cos(\zeta_m^{(k)})sin(\vartheta_t^{(k)})-x_3sin(\zeta_m^{(k)})\big)~e^{\text{i}\omega
W_{l}^{(k)}\cdot {\bf x}}.
\en

\begin{remark}
We would like to clearly separate the following three sets. Once the nonlinear parameters $D$ are fixed, the discontinuous Galerkin plane wave neural network set under such restricted meaning transforms into the subspace denoted as $V_{n,D}^{\sigma}({\mathcal T}_h)$, whose basis functions satisfy (\ref{neubasis}). Of course, the set $V_{n}^{\sigma}({\mathcal T}_h)$ defined in (\ref{vsign_general}) with the nonlinear dictionary over all admissible directions is composed of the union of such subsets, $V_{n}^{\sigma}({\mathcal T}_h)=\bigcup_D V_{n,D}^{\sigma}({\mathcal T}_h)$. 
\end{remark}

\subsection{The convergence analysis}

The following result demonstrates that the quasi-maximizer $\varphi_i^{NN} \in V_{n_i}^{\sigma}\cap B$ of
 $\langle r(u_{i-1}), \cdot \rangle$ associated with $\epsilon_i$, satisfying {\bf Definition 4.1} exists on the algorithmic level, and is a good approximation to $\varphi_i$. It also yields that $\langle r(u_{i-1}), \varphi_i^{NN} \rangle$ is equivalent to the true error $|||u-u_{i-1}|||$.

\begin{lemma} \label{kerlem}
 There exists $n(\epsilon_i, u_{i-1})\in \mathbb{N}$ such that if $n_i = n(\epsilon_i, u_{i-1})$, $\varphi_i^{NN}\in V_{n_i}^{\sigma}\cap B$ defined in {\bf Definition 4.1} exists, and is computable. It also satisfies the approximation:
\be \label{direxist}
||| \varphi_i^{NN} - \varphi_i ||| \leq \epsilon_i.
\en
Moreover, $\langle r(u_{i-1}), \varphi_i^{NN} \rangle$ is equivalent to the true error $|||u-u_{i-1}|||$, 
\be \label{equiestimator}
\frac{\langle r(u_{i-1}), \varphi_i^{NN} \rangle}{1+\epsilon_i}  \leq |||u-u_{i-1}||| \leq \frac{\langle r(u_{i-1}), \varphi_i^{NN} \rangle}{1-\epsilon_i}.
\en
\end{lemma}


{\it Proof.} For seeking a good initial value $\tilde{\varphi}_i$ that approaches $\varphi_i$, we find it by the PWLS method \cite{hy2}, such that
\be
a(\tilde{\varphi}_i, v_h) = a(u-u_{i-1},v_h) = L(v_h)- a(u_{i-1},v_h), \quad \forall v_h\in V_{n, D}^{\sigma}({\mathcal T}_h),
\en
where $n$ satisfies the following approximate estimate. By the universal approximation property (\ref{approx}) and further Remark \ref{univerappro},
there exists suitably large $n(\epsilon_i, u_{i-1})\in \mathbb{N}$ and $\tilde{\varphi}_i \in V_{n(\epsilon_i, u_{i-1}), D}^{\sigma}$ such that
$||| \varphi_i - \tilde{\varphi}_i ||| < \frac{\epsilon_i}{\epsilon_i+2}$,
 which indicates that, by the triangle inequality,
  $1- \frac{\epsilon_i}{\epsilon_i+2}  < | | |  \tilde{\varphi}_i | | |  < 1 + \frac{\epsilon_i}{\epsilon_i+2}$.
   Let $\hat{\varphi}_i  = \tilde{\varphi}_i/|||\tilde{\varphi}_i|||\in V_{n(\epsilon_i, u_{i-1}), D}^{\sigma}\cap B$.
  Then
$$|||\varphi_i-\hat{\varphi}_i|||= \frac{ ||| \varphi_i ||| \tilde{\varphi}_i ||| - \tilde{\varphi}_i ||| }{||| \tilde{\varphi}_i|||}
= \frac{ ||| \varphi_i ||| \tilde{\varphi}_i ||| - \varphi_i + \varphi_i - \tilde{\varphi}_i ||| }{||| \tilde{\varphi}_i|||}
 \leq \frac{ ||| \varphi_i ||| \tilde{\varphi}_i ||| - \varphi_i ||| + ||| \varphi_i - \tilde{\varphi}_i ||| }{||| \tilde{\varphi}_i|||}
 \leq \epsilon_i.  $$

Since $\varphi_i^{NN}$ is updated for objective functional $\langle r(u_{i-1}), \cdot \rangle$ using a randomly shuffled Adam optimizer from the initial point $\hat{\varphi}_i$, there must be
\be
\langle r(u_{i-1}), \varphi_i^{NN} \rangle \geq \langle r(u_{i-1}), \hat{\varphi}_i \rangle.
\en
Thus
\beq
\sup\limits_{v_i^{NN}\in V_{n_i}^{\sigma}\cap B} \langle r(u_{i-1}), v_i^{NN} \rangle - \langle r(u_{i-1}), \varphi_i^{NN} \rangle
& \leq &  
 \sup\limits_{v_i^{NN}\in V_{n_i}^{\sigma}\cap B} \langle r(u_{i-1}), v_i^{NN} \rangle -  \langle r(u_{i-1}), \hat{\varphi}_i \rangle
 \cr
 & \leq &  \max\limits_{v\in B} ~\langle r(u_{i-1}), v\rangle -  \langle r(u_{i-1}), \hat{\varphi}_i \rangle
 \cr
 & =  &   \langle r(u_{i-1}), \varphi_i - \hat{\varphi}_i \rangle = \text{Re} \{ a(u-u_{i-1}, \varphi_i - \hat{\varphi}_i) \}
  \cr
 & \leq &  |||u-u_{i-1}||| ~|||\varphi_i - \hat{\varphi}_i)  ||| \leq \epsilon_i |||u-u_{i-1}|||,
\eq
which coincides the numerical existence of the quasi-maximizer.

Besides, we have the approximation for the quasi-maximizer,
$$||| \varphi_i - \varphi_i^{NN}  |||^2 = 2 - 2\text{Re} \{ a(\varphi_i, \varphi_i^{NN}) \} \leq  2 - 2\text{Re} \{ a(\varphi_i, \hat{\varphi}_i ) \} = ||| \varphi_i - \hat{\varphi}_i  |||^2  \leq \epsilon_i^2.$$

Since $$\langle r(u_{i-1}), \varphi_i \rangle = \text{Re} \{ L(v) - a(u_{i-1}, \varphi_i) \} = \text{Re} \{ a(u-u_{i-1}, \varphi_i) \} = |||u-u_{i-1}|||,$$
we have, 
\be
| \langle r(u_{i-1}), \varphi_i^{NN} \rangle - |||u-u_{i-1}||| ~| = | \langle r(u_{i-1}), \varphi_i^{NN}- \varphi_i \rangle | = | \text{Re} \{ a(u-u_{i-1}, \varphi_i^{NN}- \varphi_i ) \}| \leq |||u-u_{i-1}||| \epsilon_i.
\en
By the triangle inequality, we obtain
\be
(1 - \epsilon_i) |||u-u_{i-1}||| \leq \langle r(u_{i-1}), \varphi_i^{NN} \rangle \leq (1 + \epsilon_i) |||u-u_{i-1}|||,
\en
namely (\ref{equiestimator}).

$\hfill\Box$

\begin{remark}
The determination on width $n(\epsilon_i, u_{i-1})$ can be expressed as a linear function of Galerkin iterations in the 2D case and a square function of Galerkin iterations in the 3D case, from numerical experiments.
\end{remark}

\begin{remark}
We would like to point out that the key result Lemma \ref{kerlem} was established for a different method from Proposition 2.5 of \cite{AD} and Lemma 4.1 of \cite{yuanhu}. Lemma \ref{kerlem} calculates an approximation $\varphi_i^{NN}$ tending to the true maximizer $\varphi_i$ from specific initial values. Regardless of the convergence speed of the optimizer used, it does not affect the final convergence of the DGPWNN algorithm.
\end{remark}

\begin{remark}
There is another important difference for Lemma \ref{kerlem} from Proposition 2.5 introduced in \cite{AD}. The boundedness assumption of the neural network parameters was implicitly required in Proposition 2.5 of \cite{AD}. However the standard plane wave method is heavily unstable \cite{ref21,HMPsur,HGA}, the coordinate coefficients $\{c_l^{(k)}\}$ in (\ref{approx}) are unbounded when the number of degrees of freedom on each element increases to infinity (see Section 4 of \cite{Par2022}). It means that there is no boundedness requirement on the neural network parameters in the proposed DGPWNN method.
 \end{remark}

The sequential enhancements can be shown to provide further an alternative indicator $\langle r(u_{i-1}), \varphi_i^{NN} \rangle \cdot $ $||\varphi_i^{NN}||_{L^2}$ for the errors in the $L^2$ norm that can be used as a useful criterion for terminating the sequential updates.

\begin{corollary}
We have the following equivalence estimates:
\beq \label{l2indi}
||u-u_{i-1}||_{L^2} \geq \frac{\langle r(u_{i-1}), \varphi_i^{NN} \rangle}{1+\epsilon_i}
~\big(||\varphi_i^{NN}||_{L^2} - \epsilon_i \big),
\cr
||u-u_{i-1}||_{L^2} \leq  \frac{\langle r(u_{i-1}), \varphi_i^{NN} \rangle}{1-\epsilon_i}
~\big(||\varphi_i^{NN}||_{L^2} + \epsilon_i  \big).
\eq
\end{corollary}

{\it Proof.} By the definition (\ref{defivarphi2}) of $\varphi_i$, and the triangle inequality, we have
\beq \label{l2inequ2}
||| u-u_{i-1} |||  ~(||\varphi_i^{NN}||_{L^2} - ||\varphi_i^{NN} - \varphi_i||_{L^2})
&\leq &  ||u-u_{i-1}||_{L^2}=  ||| u-u_{i-1} ||| ~||\varphi_i||_{L^2}
\cr & \leq &
||| u-u_{i-1} |||  ~(||\varphi_i^{NN}||_{L^2} + ||\varphi_i^{NN} - \varphi_i||_{L^2}).
 \eq
Besides, by the approximation property (\ref{direxist}) of $\varphi_i^{NN}$ and the error estimator $\langle r(u_{i-1}), \varphi_i^{NN} \rangle$, we get the desired equivalent formula for the error $||u-u_{i-1}||_{L^2}$.

$\hfill\Box$

The above estimates indicate that provided that $\epsilon_i$ is small enough, the numerical value of $||\varphi_i^{NN}||_{L^2} \cdot $ $\langle r(u_{i-1}), \varphi_i^{NN} \rangle $ is indeed equivalent to the true error $||u-u_{i-1}||_{L^2}$. The numerical tests will demonstrate the robustness of the error indicator.

The next result addresses the convergence of the proposed recursive algorithm.


\begin{theorem} \label{helmtheom}
Assume that, for $j=1,\cdots,i$ and $|||u-u_{j-1}|||>0$. Let $n(\epsilon_j, u_{j-1})\in \mathbb{N}$ and
$\varphi_j^{NN}\in V_{n_j}^{\sigma}\cap B$ with $n_j = n(\epsilon_j, u_{j-1})$ be determined as in Lemma \ref{kerlem} ($j=1,\cdots,i$).
Then the resulting approximate solution $u_r$ has the estimate:
\be \label{conver}
|||u-u_i||| < |||u-u_0|||~\prod_{j=1}^i  \epsilon_j.
\en
\end{theorem}

{\it Proof.} Choosing $v=u_{i-1} + |||u-u_{i-1}|||\varphi_{i}^{NN}$, we obtain
\beq \label{iteresti}
|||u-u_i||| & \leq &  ||| u-u_{i-1} - |||u-u_{i-1}|||\varphi_{i}^{NN} |||
= ||| ~\varphi_{i} |||u-u_{i-1}||| - |||u-u_{i-1}||| \varphi_{i}^{NN} ~||| \cr
 & = & |||u-u_{i-1}||| \cdot ||| \varphi_{i} - \varphi_{i}^{NN} |||
 \leq  |||u-u_{i-1}||| \cdot \epsilon_i.
\eq
Thus, we have
$$|||u-u_i|||  \leq |||u-u_{i-1}||| \cdot \epsilon_i.$$
Iterating this process forward gives the desired result.

$\hfill\Box$

Particularly, (\ref{equiestimator}) provides a computable estimate of the error $|||u-u_{i-1}|||$ that can be used as a stopping criterion for the recursive solution procedure as follows: Given a tolerance $tol > 0$, if $\eta(u_{i-1}, \varphi_i^{NN}) > tol$, then we augment the discontinuous Galerkin subspace $\Phi_{i-1}^{NN}$ to $\Phi_i^{NN}=\Phi_{i-1}^{NN}\oplus \{\varphi_i^{NN}\}$  and compute an updated approximation $u_i$; otherwise, the approximation $u_{i-1}$ is deemed satisfactory, and the procedure terminates. Additionally, owing to augmenting the discontinuous Galerkin subspace, (\ref{conver}) produces the convergence of the recursive procedure.

\begin{remark}\label{condicover}
It should be strengthened that, the above convergence result is unconditional: as long as the initial directions are chosen as suitably distributed propagation directions, regardless of which optimizer is used, the least squares solution converges for the discrete problem (\ref{dgappro}) defined in the discontinuous Galerkin space generated by the quasi-maximizers.
\end{remark}


The following estimate of condition numbers can be viewed as the extensions to the discontinuous Galerkin plane wave neural networks of the existing results of \cite[Proposition 2.11]{AD}. We omit the proof.

\begin{lemma}
 Denote by $K^{(j)}\hat u = F^{(j)}$ the linear system (\ref{dgappro}). Then the functions $\{\varphi_i^{NN}\}_{i=1}^j$ are linearly independent. Moreover, if $\gamma = \sum\limits_{i=1}^j \epsilon_i <1$, the condition number of the stiffness matrix $K^{(j)}$ satisfies $\text{cond}(K^{(j)})< \frac{1+\gamma }{1-\gamma }$.
\end{lemma}

{\it Proof.} see \cite[Proposition 2.11]{AD}.

$\hfill\Box$



\subsection{Discussions on the algorithm} \label{sectionnr}

On the choice of the width $n_i$ of the network as $i$ increases, Theorem \ref{helmtheom} indicates that if $0 < \epsilon < 1$, the mesh ${\cal T}_h$ are fixed, and $\epsilon_i=\epsilon$ for all $\epsilon$, then the error in $u_i$ is a factor of $\epsilon$ smaller than the error in $u_{i-1}$. However, as pointed out in \cite[Section 2.3.2]{AD}, if $n_i$ is fixed, the rate of convergence will become arbitrarily slow due to the decreased ability of a fixed-width network $\varphi_i^{NN}\in V^{\sigma}_{n(\epsilon, u_{i-1})}$ to capture higher resolution features of the error as $i$ increases.
In order for $\varphi_i^{NN}\in V^{\sigma}_{n(\epsilon, u_{i-1})}$ to capture the higher resolution features of $\varphi_i$, we necessarily require that $n_i=n(\epsilon, u_{i-1})>n_{i-1}=n(\epsilon, u_{i-2})$.


On the other hand, motivated by the $hp$-convergence (\ref{approx}) of the aforementioned discontinuous Galerkin method, it is expected that,
 if we decrease the meshwidth $h$ while increasing $n_i$ as $i$ increases, then the error in $u_i$ will not only reach a factor of $\epsilon$ of the error in $u_{i-1}$, but also decrease more sharply than the error generated by the case of fixed $h$. That is, a judicious $hp-$refinement strategy will be the most attractive option to solve the well-known the problem of large wavenumbers.

\begin{remark} \label{requon_p}
Owing to the employed strategy of $hp-$refinement, to achieve the same $\epsilon$, the requirement above on the width $n$ of the network for the proposed discontinuous Galerkin neural network is to a certain extent decreasing, and is obviously weaker than that of \cite{AD} for the continuous Galerkin neural network.
\end{remark}

Numerical results validate that, when choosing $h\approx \mathcal{O}(\frac{\pi}{\omega})$ and gradually increasing the width $n_i$ of the network at each discontinuous Galerkin iteration, the resulting approximation generated by the discontinuous Galerkin plane wave neural networks can reach the accuracy of a given approximation. While for the current situation, compared with the standard plane wave method, the condition number of discrete system (\ref{dgappro}) satisfies $\text{cond}(K^{(j)})=\mathcal{O}(1)$.
Additionally, the performance of the algorithms is independent of wavenumbers. In particular, only no more than ten outer discontinuous Galerkin iterations may guarantee the convergence of our plane wave neural network algorithm, while the total degrees of freedom and computational cost still increase with $\omega$.

\section{A DGPWNN algorithm for underlying Maxwell's equations} \setcounter{equation}{0}

In this section we recall the second order system of Maxwell's equations, and derive a discontinuous Galerkin plane wave neural network algorithm.

We want to compute a numerical approximation $E$ of the second order homogeneous Maxwell equations with the lowest order absorbing boundary condition in three space dimensions:
\begin{equation} \label{maxeq}
\left\{ \begin{aligned}
     &  \nabla\times(\frac{1}{i\omega\mu}\nabla\times E) + i\omega\varepsilon E=0 & \text{in}\quad
\Omega,\\
    & -E\times n+\frac{\varsigma}{i\omega\mu}((\nabla\times E)\times  n)\times
n=g & \text{on}\quad\gamma.
                          \end{aligned} \right.
                          \end{equation}
Here $\omega>0$ is the temporal frequency of the field, and $
g\in L_T^2(\partial\Omega)^3$. The material coefficients
$\varepsilon,\mu$ and $\varsigma$ are assumed to be piecewise constant in the
whole domain. In particular, if $\varepsilon$ takes complex-valued,
then the imaginary part typically represents absorption; else the material is non-absorbing (see \cite{hy2}).

\subsection{A minimization method}
For an element $\Omega_k$, denote by $\mbox{
H}(\text{curl};\Omega_k)$ the standard Sobolev space. Set
\begin{equation}
 V(\Omega_k)=\bigg\{E_k\in \mbox{
H}(\text{curl};\Omega_k); ~ E_k \text{~elementally satisfies the
first equation of (\ref{maxeq})}\bigg\}.
\end{equation}
Define
 \be
V({\cal T}_h) = \{ v\in L^2(\Omega): ~v|_{\Omega_k}\in V(\Omega_k) ~\forall \Omega_k\in {\cal T}_h  \}.
\en

 Let $\rho_1$ and $\rho_2$ be two given relaxed constant. The tangential jump of a quantity across an interface $\Gamma_{kj}$ $(k<j)$ between two subdomains $\Omega_k$
  and $\Omega_j$ is typically defined as follows.
\be
\llbracket F\times  n \rrbracket  =  F_k\times n_k + F_j\times n_j.
 \en

Define the functional
\beq \label{max_mini_func}
J(F)&=&\sum_{k=1}^N\int _{\gamma_k}| - F_k\times  n_k
 + \frac{\varsigma}{i\omega\mu}((\nabla\times  F_k)\times
n_k) \times n_k - g|^2~ds \cr&+&
 \sum_{k<j} \Bigg(\int_{\Gamma_{kj}}
\rho_1 |\llbracket{ F}\times{ n} \rrbracket|^2~~ds
+
 \int_{\Gamma_{kj}} \rho_2 |\llbracket {1\over i\omega\mu}(\nabla\times{F})\times n \rrbracket|^2 ~ds \bigg), \quad F\in V({\cal T}_h).
\eq

Consider the minimization problem: find $E \in V({\cal T}_h)$ such that
\begin{equation} \label{maxminmum}
J({E})=\min\limits_{{ F}\in  V({\cal T}_h)}~J( F).
\end{equation}

The variational problem related to the minimization problem (\ref{maxminmum}) typically involves finding a function $E \in V({\mathcal T}_h)$ such that (refer to \cite{hy2}):

\begin{equation}
\begin{split}
 &\sum_{k=1}^N\int _{\gamma_k}\bigg(-E_k\times  n_k + \frac{\varsigma}{i\omega\mu}\big((\nabla\times  E_k)\times
n_k\big) \times n_k - g\bigg)\cdot\overline{-F_k\times n_k + \frac{\sigma}{i\omega\mu}\big((\nabla\times  F_k)\times
n_k\big) \times n_k}~ds
\\ & +
 \sum_{k<j}\int_{\Gamma_{kj}}\Bigg(\rho_1 \llbracket  E\times  n
  \rrbracket \cdot \overline{\llbracket  F\times  n \rrbracket}
+ \rho_2
\llbracket {1\over i\omega\mu}(\nabla\times E)\times n \rrbracket \cdot \overline{ \llbracket {1\over i\omega\mu}(\nabla\times F)\times n \rrbracket }\bigg) ~ds =0.
\end{split}
\label{MAXVAR}
\end{equation}

Define the sesquilinear form $a(\cdot, \cdot)$ by
\beq \nonumber
 a(E,F) =  \sum_{k<j}\int_{\Gamma_{kj}}\Bigg(\rho_1 \llbracket  E\times  n
  \rrbracket \cdot \overline{\llbracket  F\times  n \rrbracket}
+ \rho_2
\llbracket {1\over i\omega\mu}(\nabla\times E)\times n \rrbracket \cdot \overline{ \llbracket {1\over i\omega\mu}(\nabla\times F)\times n \rrbracket }\bigg) ~ds
 \cr
 +\sum_{k=1}^N\int _{\gamma_k} \bigg(-E_k\times  n_k + \frac{\varsigma}{i\omega\mu}\big((\nabla\times  E_k)\times
n_k\big) \times n_k \bigg) \cdot \overline{-F_k\times n_k + \frac{\sigma}{i\omega\mu}\big((\nabla\times  F_k)\times
n_k\big) \times n_k}~ds,
\eq
and semilinear form $L(F)$ by
\be \nonumber
L(F)=\sum_{k=1}^N\int _{\gamma_k} g \cdot ~\overline{-F_k\times n_k + \frac{\varsigma}{i\omega\mu}\big((\nabla\times  F_k)\times
n_k\big) \times n_k}~ds.
\en
 Then (\ref{MAXVAR}) can be written as the following variational equation:
 \begin{eqnarray}\label{maxpwlsvar}
\left\{\begin{array}{ll} \text{Find} \ E\in V({\mathcal T}_h) ~~s.t. \\
a(E,F)=L(F) \quad \forall F\in V({\mathcal T}_h).
\end{array}\right.
\end{eqnarray}

We can see that $a(F,F)$ is a norm on $V({\mathcal T}_h)$, and denote it by $||| F |||^2$. Using the Cauchy-Schwarz inequality, it can be proved that $|a(E,F)|\leq ||| E |||~||| F |||$ and $|L(F)| \leq ||g||_{L^2(\gamma)} ~||| F |||, ~\forall E,$ $F \in V({\mathcal T}_h)$. Thus, a direct application of the Riesz representation theorem \cite{FO} gives rise to the existence of a unique $E \in V({\mathcal T}_h)$ satisfying (\ref{maxpwlsvar}).


\subsection{Element-wise plane wave neural networks}

A family of plane wave solutions to the constant-coefficient Maxwell equations are generated on every $\Omega_k \in {\cal T}_h$ by selecting $n$ unit propagation directions $W^{(k)}_{l} ~(l=1,\cdots,n)$ and the associated unit real polarization vectors $G^{(k)}_{l}$ orthogonal to $ W^{(k)}_{l}$. The complex polarization vectors $F^{(k)}_{l}$ and $F^{(k)}_{l+n}$ are then constructed from $W^{(k)}_{l}$ and $G^{(k)}_{l}$. These typically take the form:
\be \label{polarizavec}
 F^{(k)}_{l}=G^{(k)}_{l}, \quad F^{(k)}_{l+n}=G^{(k)}_{l}\times
W^{(k)}_{l}~~~(l=1,\cdots,n).
\en
It yields two pairs of functions $E^{(k)}_{l}$ (refer to \cite{pwdg}) satisfying the Maxwell system (\ref{maxeq}), respectively,
\begin{equation}
E^{(k)}_{l}=\sqrt{\mu}~F^{(k)}_{l}~e^{\text{i}\kappa
W^{(k)}_{l}\cdot x}, \quad
 E^{(k)}_{l+n}=\sqrt{\mu}~ F^{(k)}_{l+n}~e^{\text{i}\kappa
W^{(k)}_{l}\cdot x}~~~
(l=1,\cdots, n),\label{maxeq17}
\end{equation}
where $\kappa=\omega\sqrt{\mu\varepsilon}$. 

Plane wave propagation directions $\{W^{(k)}_{l}\}$ are defined in spherical coordinates, similarly to section 3. However, for a given direction $W^{(k)}_l$, there are multiple orthogonal vectors. For now, $W^{(k)}_l$ is a unit vector and set as $W^{(k)}_l=(a_l  ~~b_l  ~~c_l)^T$. Without loss of generality, we assume $|b_l|<1$, then choose $G^{(k)}_l$ as follows
$$G^{(k)}_l = \bigg(\frac{a_lb_l}{\sqrt{1-b_l^2}} ~~~ -\sqrt{1-b_l^2}  ~~~ \frac{-a_l^2b_l+b_l(1-b_l^2)}{c_l\sqrt{1-b_l^2}} \bigg)^T.$$

Subsequently, a hybrid approach combining neural networks with the discontinuous Galerkin (DG) finite element framework generates a vector discontinuous Galerkin neural network consisting of a single hidden layer of $2n\in \mathbb{R}$ neurons on each element $\Omega_k\in {\cal T}_h$ and defining an output function $\varphi^{NN}: \mathbb{R}^3 \rightarrow \mathbb{C}$ as follows:
\be \label{max_nn_app}
\varphi^{NN}(x; \theta)|_{\Omega_k} =  \sum_{m=1}^{m^{\ast}}\sum_{t=1}^{t^{\ast}}
\big( c_l^{(k)} E_{l}^{(k)} + c_{l+n}^{(k)} E_{l+n}^{(k)} \big), ~l=(m-1)k^{\ast}+t.
\en

For consistency, we shall use the elementwise notation $\zeta|_{\Omega_k}=[\zeta_1^{(k)}, \cdots, \zeta^{(k)}_{m^{\ast}}]^T \in \mathbb{R}^{m^{\ast}}$, $\vartheta|_{\Omega_k}=[\vartheta_1^{(k)}, \cdots, \vartheta^{(k)}_{t^{\ast}}]^T \in \mathbb{R}^{t^{\ast}}$, and $c|_{\Omega_k}=[c^{(k)}_1, \cdots, c^{(k)}_{2n}]^T\in \mathbb{C}^{2n}$. A set of nonlinear and linear parameters over a collection of all elements is denoted by $\theta=\{\zeta; \vartheta; c\}$. Particularly, a set of nonlinear parameters defined on every element, is collectively denoted by $D=\{\zeta; \vartheta\}$.

Denote by $V_n^{\sigma}({\cal T}_h)$ the set of all functions of the form (\ref{max_nn_app}), which is defined by
\be \label{maxvsign}
V_n^{\sigma}({\cal T}_h) := \bigg\{\varphi: ~\varphi(x)|_{\Omega_k} =  \sum_{m=1}^{m^{\ast}}\sum_{t=1}^{t^{\ast}}
\big( c_l^{(k)} E_{l}^{(k)} + c_{l+n}^{(k)} E_{l+n}^{(k)} \big), ~l=(m-1)k^{\ast}+t,~\forall \Omega_k\in {\cal T}_h  \bigg\}.
\en

\subsection{Discontinuous Galerkin plane wave neural network discretization}
As in section \ref{helmdis}, we can design iterative algorithms to generate discontinuous Galerkin plane wave neural network solutions that approximately satisfy the minimization problem (\ref{maxminmum}).
We would like to emphasize that updating hidden parameters $D$ requires calculating derivative $\nabla_D \eta(E_{i-1},\varphi)$ of the quantity $\eta(E_{i-1},\varphi)$ similarly defined in (\ref{etaquan}). By the direct reduction, we get
\be
\nabla_D \eta(E_{i-1}, \varphi) = \frac{\nabla_D \langle r(E_{i-1}), \varphi \rangle }{|||\varphi|||} + \langle r(E_{i-1}), \varphi \rangle~ \nabla_D \frac{1}{|||\varphi|||},
\en
where
\be
\nabla_D \langle r(E_{i-1}), \varphi  \rangle = \text{Re}\big\{ \langle \xi(E_{i-1}), \nabla_D \varphi \rangle \big\}
\en
and
\be
\nabla_D \big(\frac{1}{|||\varphi|||}\big) = -|||\varphi|||^{-3}\text{Re}\{a(\varphi, \nabla_D \varphi)\}.
\en

Exactly, taking as an example $\varphi(x)=  \sum_{m=1}^{m^{\ast}}\sum_{t=1}^{t^{\ast}}
\big( c_l^{(k)} {\bf E}_{l}^{(k)} + c_{l+n}^{(k)} {\bf E}_{l+n}^{(k)} \big), ~l=(m-1)k^{\ast}+t;~~{\bf x}\in \Omega_k \in {\cal T}_h$, we have
\beq \nonumber 
\nabla_{D^{(k)}} \varphi(x)= \sum_{m=1}^{m^{\ast}}\sum_{t=1}^{t^{\ast}}  \bigg(c^{(k)}_l \sqrt{\mu}~ \big(\nabla_{D^{(k)}} F^{(k)}_l + i\kappa F^{(k)}_l x^T \nabla_{D^{(k)}} W^{(k)}_l \big)
\\
+ c^{(k)}_{l+n} \sqrt{\mu}~  \big(\nabla_{D^{(k)}} F^{(k)}_{l+n} + i\kappa F^{(k)}_{l+n} x^T \nabla_{D^{(k)}} W^{(k)}_l \big) \bigg) ~e^{\text{i}\kappa
W^{(k)}_{l}\cdot x},
\eq
 where
 $\nabla_{D^{(k)}} W^{(k)}_l=(\nabla_{D^{(k)}} a_l ~\nabla_{D^{(k)}} b_l ~\nabla_{D^{(k)}} c_l)^T$, and the calculation of $\nabla_{D^{(k)}} F^{(k)}_l$ and $\nabla_{D^{(k)}} F^{(k)}_{l+n}$ is similarly to $\nabla_{D^{(k)}} W^{(k)}_l$'s. 

The following convergence estimate can be proved as in the proof of Theorem \ref{helmtheom}.



\begin{theorem} Assume that, for $j=1,\cdots,i$ and $|||E-E_{j-1}|||>0$. Let $n(\epsilon_j, E_{j-1})\in \mathbb{N}$ and
$\varphi_j^{NN}\in V_{n_j}^{\sigma}\cap B$ with $n_j = n(\epsilon_j, E_{j-1})$ be determined as in Lemma \ref{kerlem} ($j=1,\cdots,i$).
Then the resulting approximate solution $E_i$ has the estimate:
\be
|||E-E_i||| < |||E-E_0|||~ \prod_{j=1}^i \epsilon_j.
\en
\end{theorem}

\section{Numerical experiments} \setcounter{equation}{0}

In this section we report some numerical results to validate the effectiveness of the discontinuous Galerkin plane wave neural network. For the case of Helmholtz equation, as pointed out in Section 3 of \cite{hy}, we choose the weighted parameters $\alpha$ and $\beta$ in the variational
problem (\ref{pwlsvar}) as $\alpha=\omega^2$ and $\beta=1$.
Besides, for Maxwell's equations, we assume that \(\varepsilon=1+i\), and apply the PWLS method with
\(\rho_1= \rho_2=1\).

A uniform triangulation \(\mathcal {T}_h\) for the domain \(\Omega\) is employed in the examples tested in this section. In the two-dimensional case, \(\Omega\) is divided into equal-sized rectangles. In the three-dimensional case, \(\Omega\) is divided into small cubes of equal size.

The initial iterative solutions $u_0$ and $E_0$ are both uniformly set to be $0$. We uniformly set the tolerance $tol$ to be $tol = 10^{-6}$, and adaptively construct a discontinuous Galerkin subspace that achieves an approximation with $\eta (u_{i-1}, \varphi^{NN}_i) < tol$. Generally, set the number of training epoches as 50 in Algorithm 2. Moreover, in Algorithm 2, the termination condition for training epochs is that $||\nabla_D \eta(u_0,v)||_{\infty} < \varepsilon$, which will be determined in each numerical test.


A set of plane wave propagation angles $D$ introduced in subsection \ref{hlempara} is uniformly initialized as follows. For the 2D case, the hidden parameters are initialized such that
\be \nonumber
d_j^{(i)}=-\pi+\frac{2\pi}{n_i}j, ~j=1,\cdots,n_i.
\en
 For the 3D case, the hidden parameters are typically initialized such that
 \be  \nonumber
\zeta_m^{(i)}= \frac{\pi-2\varepsilon_{disturb}}{m^{\ast}_{i}-1}(m-1)+\varepsilon_{disturb}, ~m=1,\cdots,m^{\ast}_{i}
 \en
 and
 \be  \nonumber
 \vartheta_t^{(i)}=-\pi+\frac{2\pi}{t^{\ast}_{i}}t, ~t=1,\cdots,t^{\ast}_{i}=2m^{\ast}_{i},
 \en
where $m^{\ast}_{i}$ and $t^{\ast}_{i}$ denote the widths $m^{\ast}$ and $t^{\ast}$ (defined in (\ref{3ddirec})) of the network architecture for the $i$-th quasi-maximization iteration, respectively.

For evaluating the objective function and training the discontinuous Galerkin plane wave neural network, we employ a fixed complex Gauss-Legendre quadrature rule to approximate all inner products. We perform all computations on a Dell Precision T5500 graphics workstation (2*Intel Xeon X5650 and 6*12GECC) using homemade \textsc{MATLAB} implementations.

\subsection{An example of two-dimensional Helmholtz equation}
An exact solution (see \cite{HGA}) to the problem can be expressed in the closed form as $$u_{ex}(x,y)=\text{cos}(k\pi y)(A_1e^{-i\omega_x x}+A_2e^{i\omega_x
x})$$ where $\omega_x=\sqrt{\omega^2-(k\pi)^2}$, and coefficients
$A_1$ and $A_2$ satisfy the equation
\begin{equation}
\left( {\begin{array}{cc} \omega_x & -\omega_x \\
(\omega-\omega_x)e^{-2i\omega_x} & (\omega+\omega_x)e^{2i\omega_x}
\end{array} }
\right)
 \left ( {\begin{array}{c} A_1 \\ A_2
\end{array}}
\right ) = \left ( {\begin{array}{c}
-i   \\
0
\end{array}}
\right )
\end{equation}

The solution represents propagating modes and
evanescent modes, respectively, when the mode number $k$ is below the cut-off value
$k\leqslant k_{\text{cut-off}}={\omega\over\pi}$ and up the cut-off
value $k> k_{\text{cut-off}}$. To be simple, we focus on the highest propagating mode with
$k=\omega/\pi-1$ in the following tests.


In order to validate Remark \ref{requon_p}, i.e., the requirement on the width $n$ of the network is to a certain extent decreasing owing to the employed strategy of $hp-$refinement, we consider the case of one subdomain $N=1$ and the case of multi-subdomains $N>1$, respectively. In subsection \ref{largewave}, the case of large wavenumbers is considered. In subsection \ref{comdgnnpwls}, we give a comparison between the DGPWNN and the PWLS method.

\subsubsection{The case of one subdomain $N=1$}
We consider both the case when successive networks have fixed widths $n_i = 14$ for all $i$ and the case when successive networks have increasing widths $n_i = 2i+1$. Set $\omega=2\pi$ and $\varepsilon=10^{-6}$.

\begin{figure}[H]
\begin{center}
\begin{tabular}{cc}
\epsfxsize=0.4\textwidth\epsffile{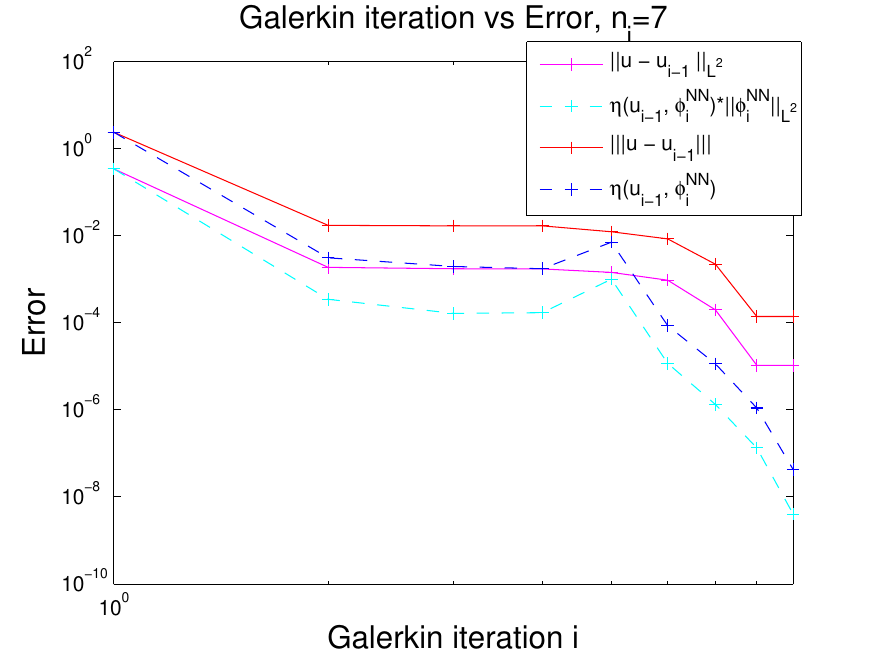}&
\epsfxsize=0.4\textwidth\epsffile{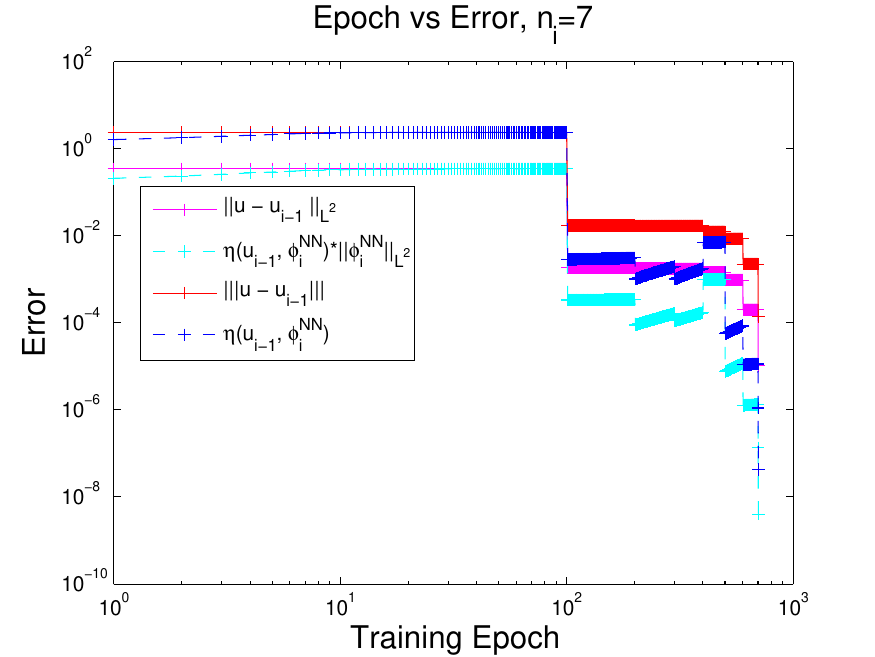}\\
\epsfxsize=0.4\textwidth\epsffile{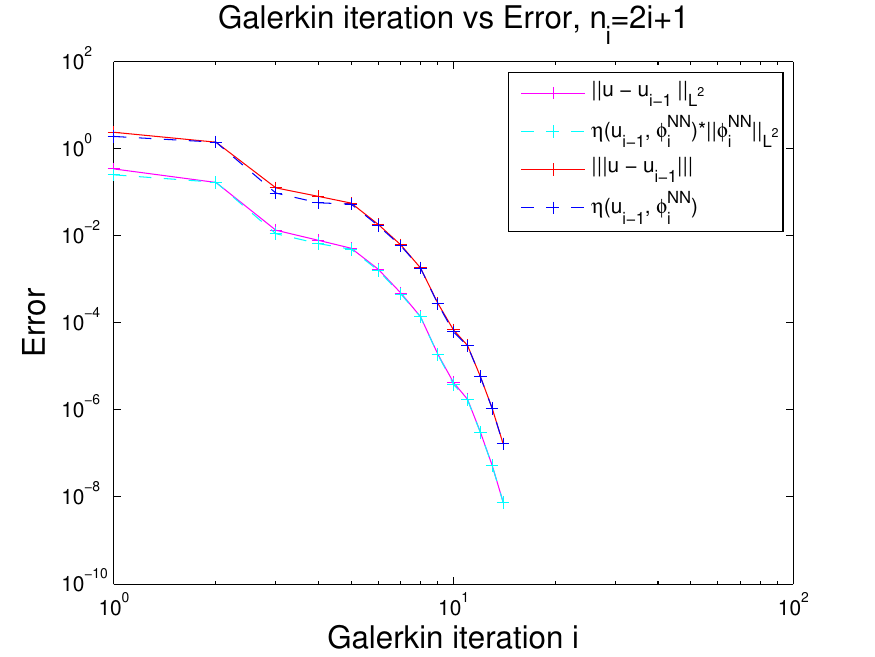}&
\epsfxsize=0.4\textwidth\epsffile{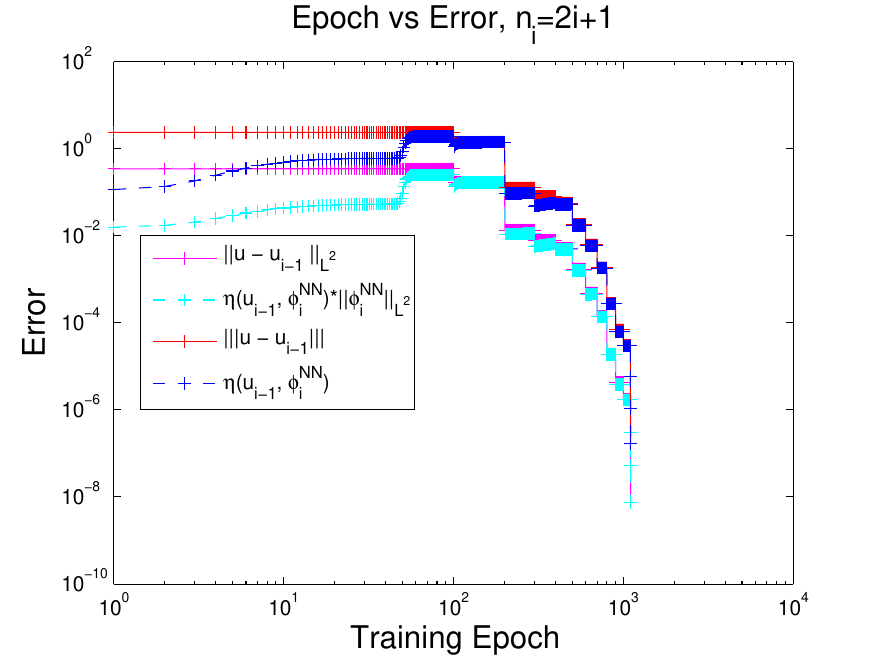}\\
\end{tabular}
\end{center}
 \caption{Displacement of a string (Isotropic Helmholtz equation in two dimension). (Left-Up) Estimated and true error in the $L^2$ and energy norms at each Galerkin iteration, $n_i = 14$ for all i. (Right-Up) The piecewise constant segments (solid lines) denote the true errors for each Galerkin iteration $i$, while the dashed blue line (resp., cyan line) denotes the progress of the loss function (resp., $L^2$ error estimator $\eta(u_{i-1}, \varphi_i^{NN}) ~| | \varphi_i^{NN}| |_{L^2}$) in approximating the true errors within each Galerkin iteration with $n_i=14$. The x-axis thus denotes the cumulative training epoch over all Galerkin iterations. (Left-Bottom) and (Up-Bottom) The analogous results for $n_i =2i+1$. }
\label{2dnn1}
\end{figure}

Figure \ref{2dnn1} shows the true errors $||u-u_{i-1}||_{L^2}$ and $|||u-u_{i-1}|||$ at the end of each discontinuous Galerkin iteration as well as the corresponding error estimates $\eta(u_{i-1}, \varphi_i^{NN}) ~| | \varphi_i^{NN}| |_{L^2}$ and $\eta(u_{i-1}, \varphi_i^{NN})$. We also provide the analogous results after each training epoch.

We observe first that the initial Galerkin iterations reduce the error substantially, with further decreases beginning to stagnate around $0.1-1.0\%$ irrespective of the number of neurons or layers used to define the underlying network architecture, and increasingly poor matching of the estimator $\eta(u_{i-1}, \varphi_i^{NN})$ of the true error as $i$ grows larger whenever $n_i$ is constant, owing to the poor approximations $u_i$ to the analytic solution,  which further, to the contrary, coincides with the validity of (\ref{equiestimator}) and (\ref{l2indi}); stagnation is not an issue when $n_i$ is increased with $i$. In particular, compared with the case that $n_i$ is constant, the case that $n_i$ is increased with $i$ provides more accurate approximations $u_i$ and estimators $\eta(u_{i-1}, \varphi_i^{NN}) ~| | \varphi_i^{NN}| |_{L^2}$ and $\eta(u_{i-1}, \varphi_i^{NN})$ of the true errors, with negligible discrepancies beginning to show in each Galerkin iteration.

Next, Table \ref{2dhelmcond_onedomain} shows the condition number of the matrix $K^{(j)}$ associated with the linear system (\ref{dgappro}) at each discontinuous Galerkin iteration $j$. We observe that the condition number is approximately unity as expected.

\vskip 0.1in
\begin{center}
       \tabcaption{Condition number for the discontinuous Galerkin matrix at each iteration $j$.}
\label{2dhelmcond_onedomain}
       \vskip 0.1in
\begin{tabular}{|c|c|c|c|c|c|c|c|c|c|c|} \hline
   $j$  & 5 & 6 &  7 &  8 &  9 & 10 &  11  &  12 & 13 &14 \\ \hline
$\text{cond}(K^{(j)})$  & 1.01 & 1.03 &  1.05 &  1.09 &  1.13 & 1.17&  1.21  &  1.29 &  1.35 & 1.42
 \\ \hline
   \end{tabular}
     \end{center}


\vskip 0.1in

\begin{figure}[H]
\begin{center}
\begin{tabular}{c}
\epsfxsize=0.5\textwidth\epsffile{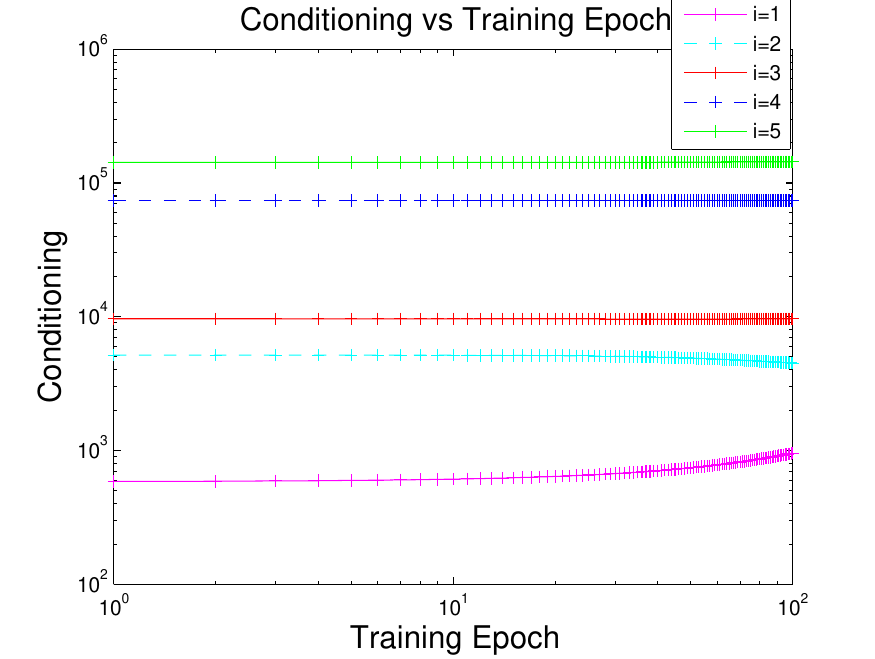}\\
\end{tabular}
\end{center}
 \caption{ Conditioning vs Training Epoch for some outer iterations.}
\label{2dfircond}
\end{figure}

Figure \ref{2dfircond} shows the conditioning of the intermediate systems in inner iterations associated with the linear system (\ref{galerlsq}) for some outer iterations.

It is not hard to see that, the auxiliary plane wave least-squares systems solved during basis construction still suffers from ill-conditioned behavior \cite{ref21,hy2,peng,peng2} for many close directions when successive networks have increasing widths $n_i$.

\begin{figure}[H]
\begin{center}
\begin{tabular}{cc}
\epsfxsize=0.4\textwidth\epsffile{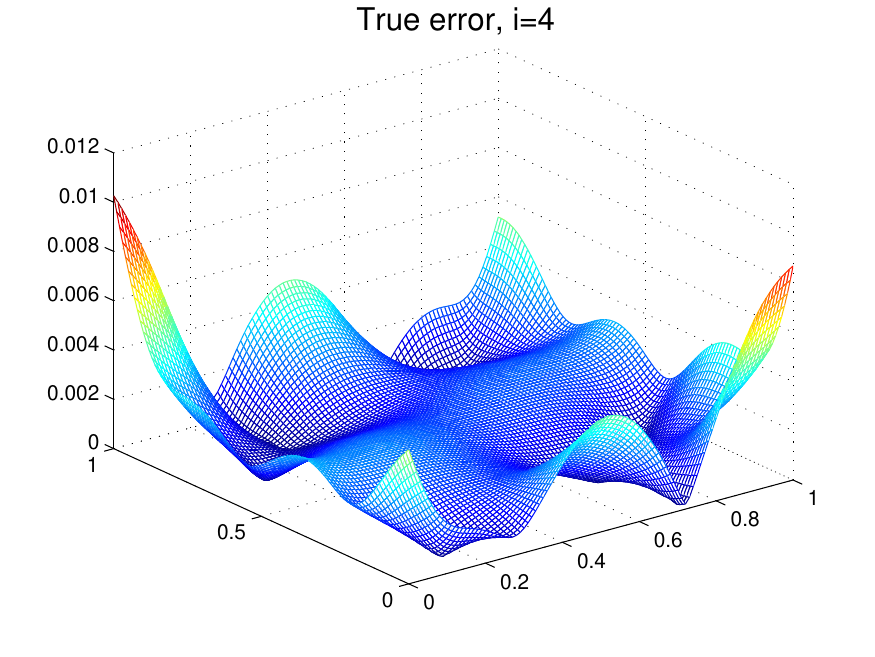}&
\epsfxsize=0.4\textwidth\epsffile{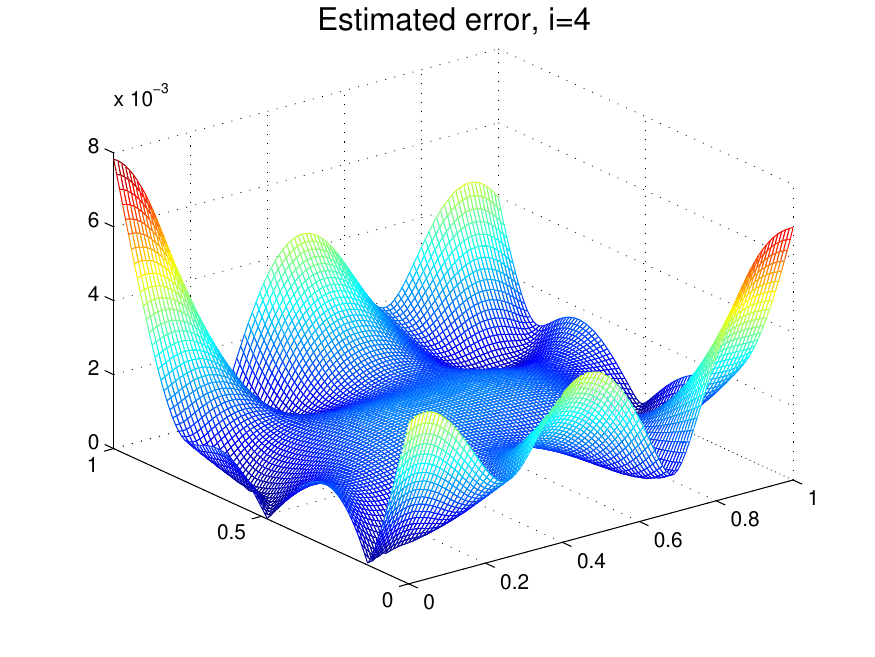}\\
\epsfxsize=0.4\textwidth\epsffile{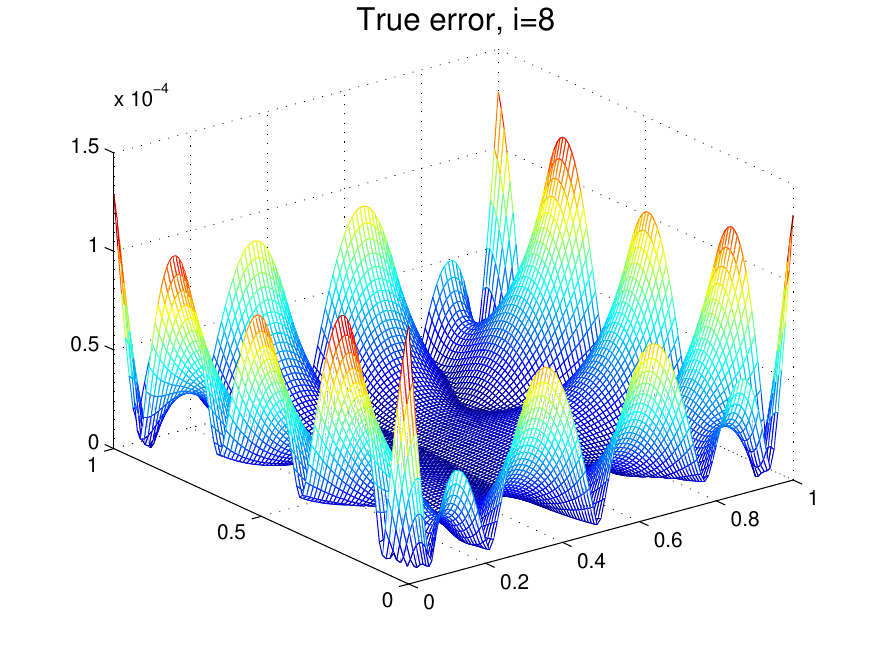}&
\epsfxsize=0.4\textwidth\epsffile{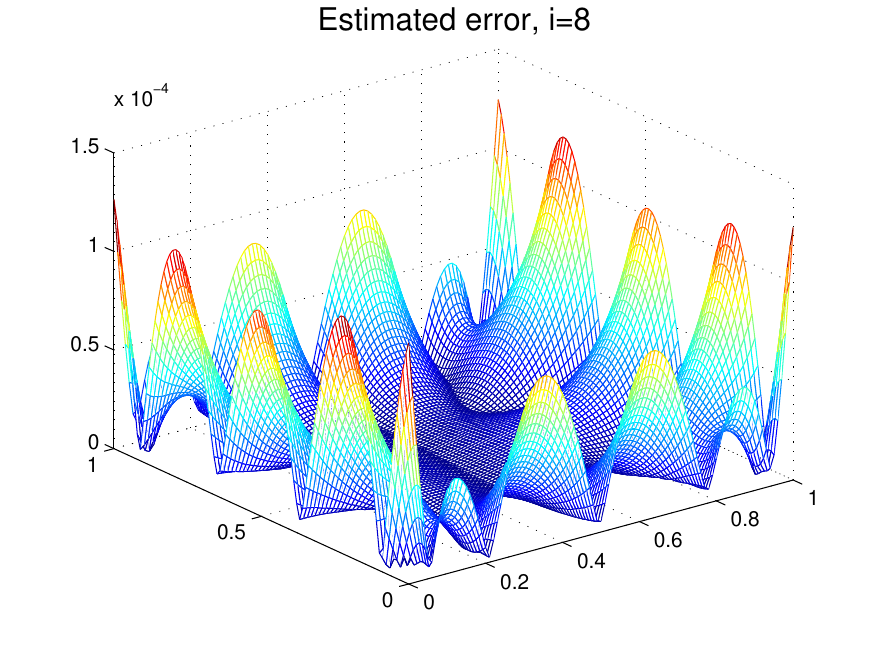}\\
\epsfxsize=0.4\textwidth\epsffile{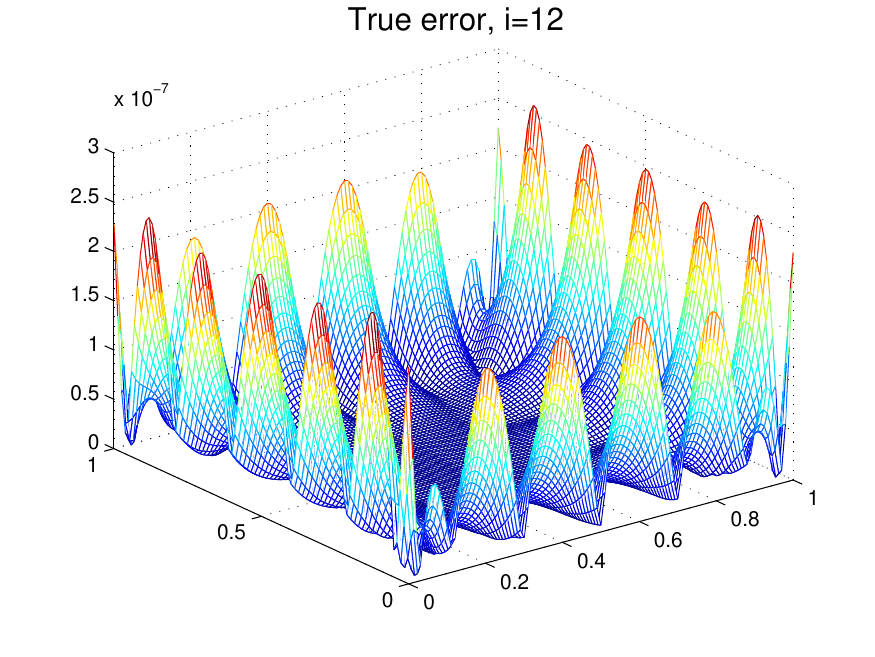}&
\epsfxsize=0.4\textwidth\epsffile{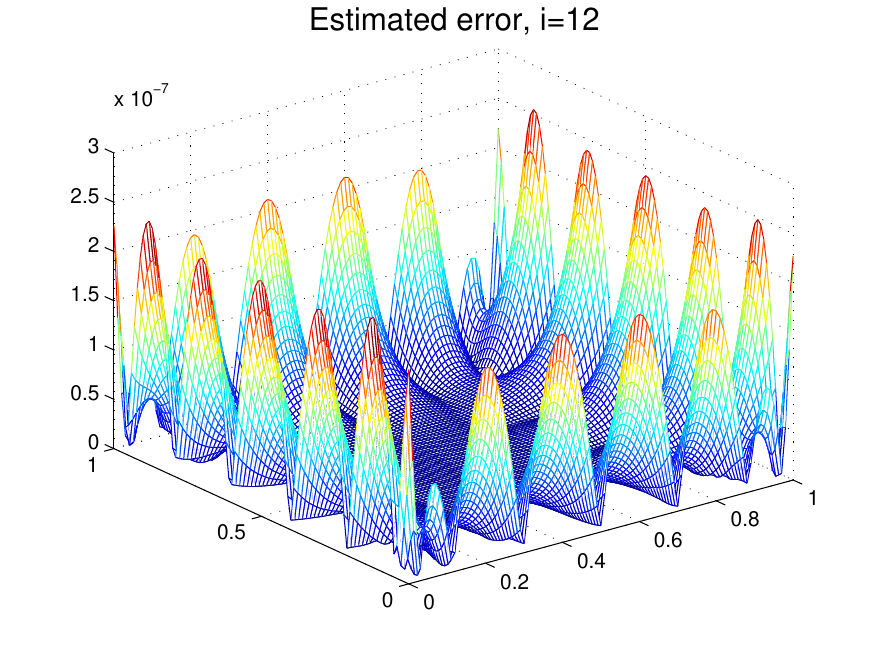}\\
\end{tabular}
\end{center}
 \caption{Displacement of a membrane (Helmholtz equation in two dimensions). Exact error $|u-u_{i-1}|$ (left column) and approximate error $\eta(u_{i-1}, \varphi_i^{NN}) ~|\varphi_i^{NN}|$ (right column) for i = 4, 8, 12. }
\label{2dnn2}
\end{figure}

Last, Figure \ref{2dnn2} shows the exact error $|u-u_{i-1}|$ as well as the maximizer $\eta(u_{i-1}, \varphi_i^{NN}) ~|\varphi_i^{NN}|$ at several stages of the algorithm.

It is apparent that initially, the low frequency components of the error are learned, with later iterations learning the high frequency error components. Besides, the pointwise error is consistent with the Galerkin iteration error in Figure \ref{2dnn1}.

\subsubsection{The case of multi-subdomains $N>1$}
\begin{figure}[H]
\begin{center}
\begin{tabular}{cc}
\epsfxsize=0.4\textwidth\epsffile{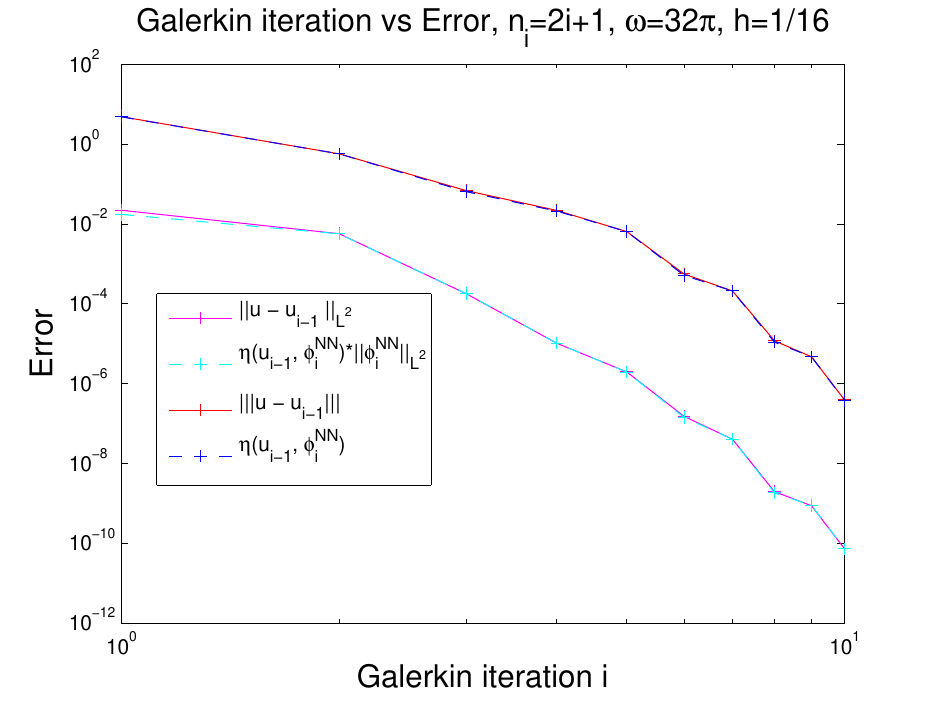}&
\epsfxsize=0.4\textwidth\epsffile{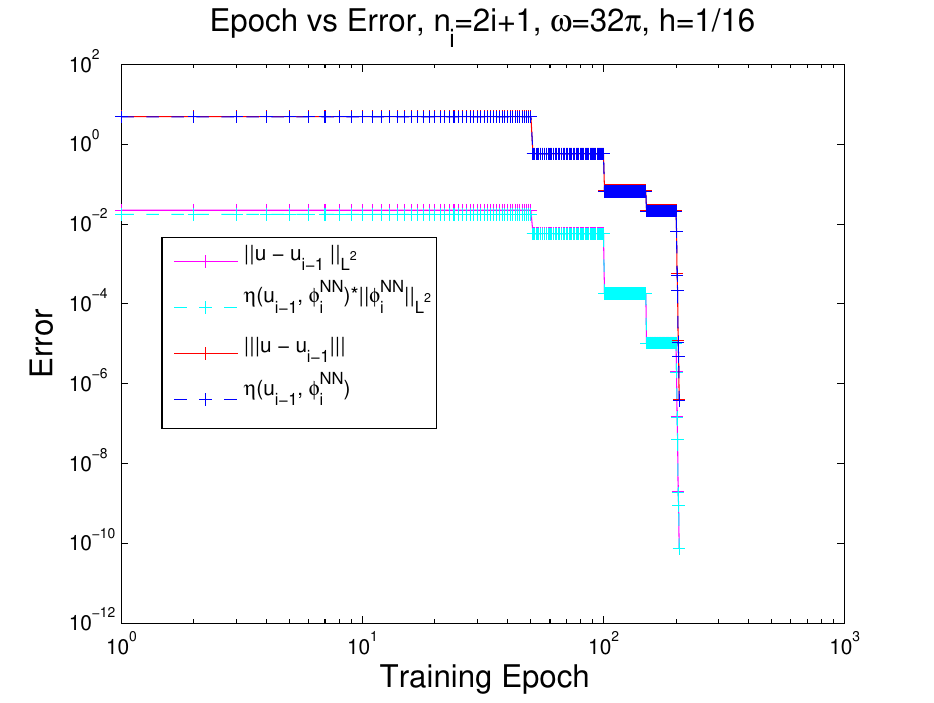}\\
\epsfxsize=0.4\textwidth\epsffile{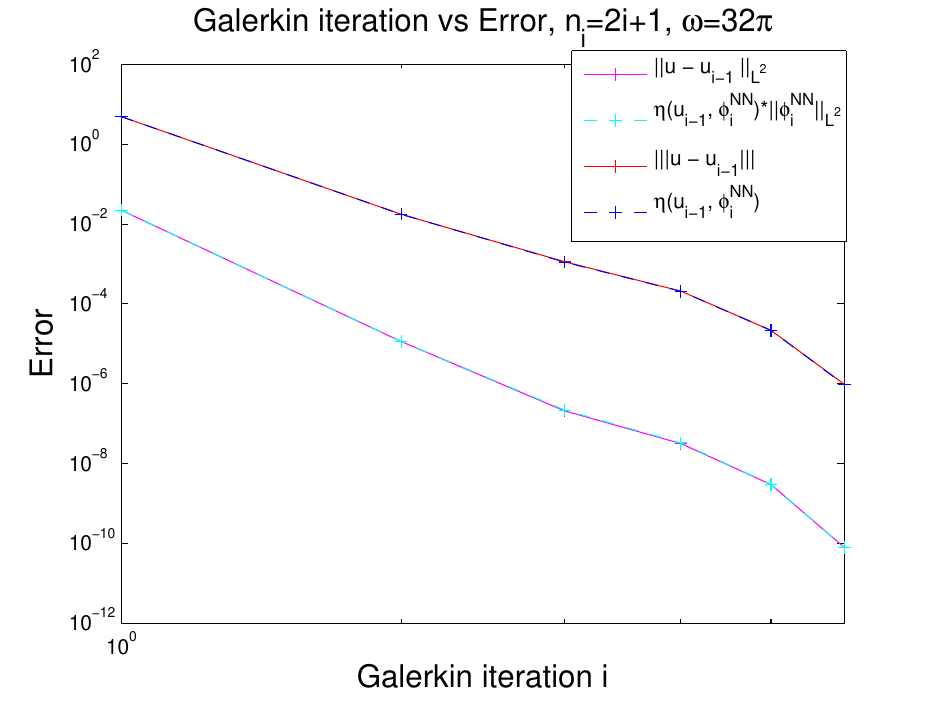}&
\epsfxsize=0.4\textwidth\epsffile{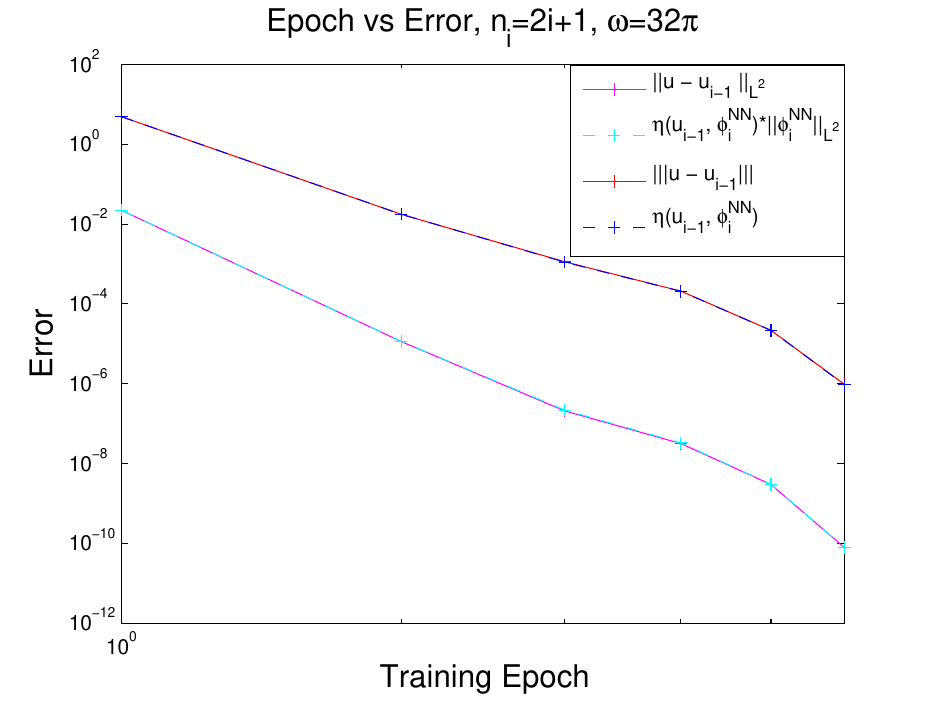}\\
\epsfxsize=0.4\textwidth\epsffile{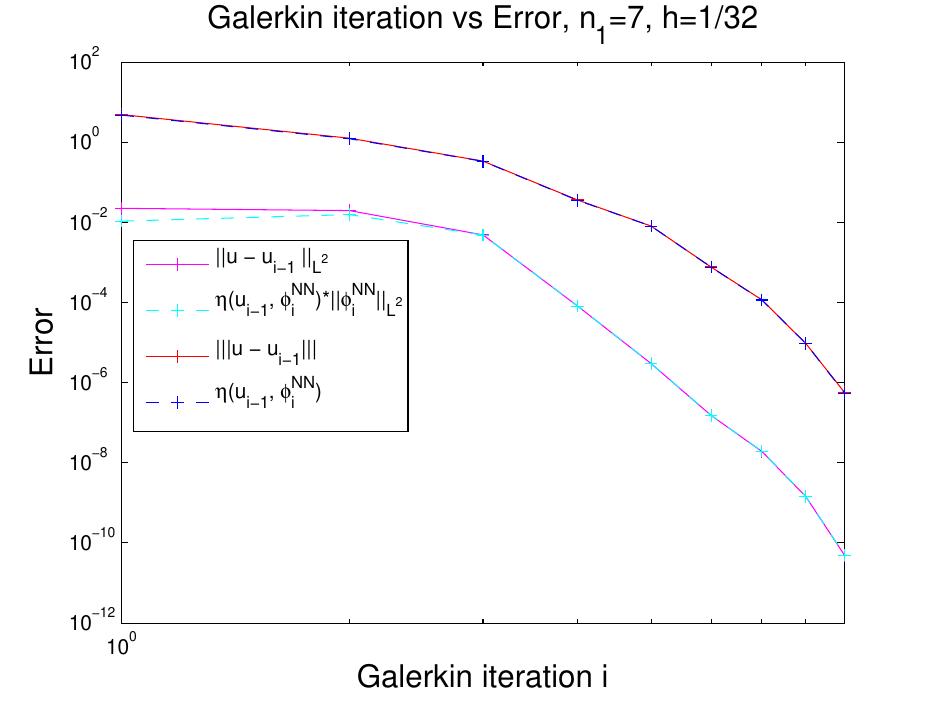}&
\epsfxsize=0.4\textwidth\epsffile{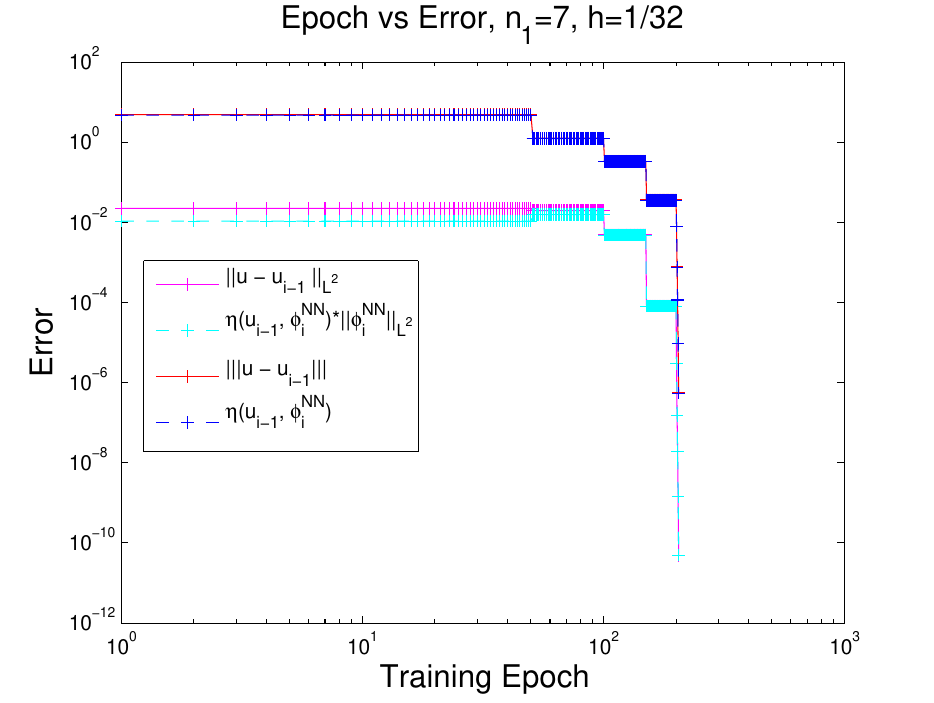}\\
\end{tabular}
\end{center}
 \caption{Displacement of a string (Isotropic Helmholtz equation in two dimension). (Left-Up) Estimated and true error in the $L^2$ and energy norms at each Galerkin iteration, $h=\frac{1}{16}$, $n_i = 2i+1$ with $n_1=13$. (Right-Up) The piecewise constant segments (solid lines) denote the true errors for each Galerkin iteration $i$, while the dashed blue line (resp., cyan line) denotes the progress of the loss function (resp., $L^2$ error estimator $\eta(u_{i-1}, \varphi_i^{NN}) ~| | \varphi_i^{NN}| |_{L^2}$) in approximating the true errors within each Galerkin iteration with $h=\frac{1}{16}$ and $n_i = 2i+1$. The x-axis thus denotes the cumulative training epoch over all Galerkin iterations. (Left-Middle) and (Up-Middle) The analogous results for $h=\frac{1}{32}$ and $n_i =2i+1$ with $n_1=13$. (Left-Bottom) and (Up-Bottom) The analogous results for $h=\frac{1}{32}$ and $n_i =2i+1$ with $n_1=7$.}
\label{2dmul_nn1}
\end{figure}

Set $\omega=32\pi$, $\varepsilon=10^{-6}$. Figure \ref{2dmul_nn1} shows the true errors $||u-u_{i-1}||_{L^2}$ and $|||u-u_{i-1}|||$ at the end of each Galerkin iteration as well as the corresponding error estimates $\eta(u_{i-1}, \varphi_i^{NN}) ~| | \varphi_i^{NN}| |_{L^2}$) and $\eta(u_{i-1}, \varphi_i^{NN})$. We also provide the analogous results after each training epoch.
  For testing $h-$convergence of the error, we choose $n_i=2i+1$, starting from $n_1=13$, and reduce the mesh size $h$ from $1/16$ to $1/32$; see Up and Middle of Figure \ref{2dmul_nn1}.
 Then, in order to validate Remark \ref{requon_p}, i.e., the requirement on the width $n$ of the network is to a certain extent decreasing owing to the employed strategy of $hp-$refinement, we decrease both $h$ and the starting network width $n_1$; see Up and Bottom of Figure \ref{2dmul_nn1}.


We observe first that the initial Galerkin iterations reduce the error substantially, with further sharp decreases to a given precision. Besides, it provides accurate approximations $u_i$ and estimators $\eta(u_{i-1}, \varphi_i^{NN}) ~| | \varphi_i^{NN}| |_{L^2}$ and $\eta(u_{i-1}, \varphi_i^{NN})$ of the true errors, with negligible discrepancies beginning to show in each Galerkin iteration.

In addition, it can be seen from the top and middle row that, $hp-$refinement strategy does accelerate the convergence of the algorithm when $h$ is reduced. Finally, from the top and bottom row, it confirms the validity of Remark \ref{requon_p}, that is, to achieve almost the same approximation accuracy, the requirement above on the width $n$ of the network for the proposed DGPWNN is to a certain extent decreasing owing to the employed strategy of $hp-$refinement.

Table \ref{2dhelmcond_multidomain} shows the condition number of the matrix $K^{(j)}$ associated with the linear system (\ref{dgappro}) at each discontinuous Galerkin iteration $j$ when $h=\frac{1}{32}$ and $n_j =2j+1$ with $n_1=7$. We observe that the condition number is approximately unity as expected.

\vskip 0.1in
\begin{center}
       \tabcaption{Condition number for the discontinuous Galerkin matrix at each iteration $j$.}
\label{2dhelmcond_multidomain}
       \vskip 0.1in
\begin{tabular}{|c|c|c|c|c|c|c|c|c|c|} \hline
   $j$  & 1 & 2 &  3 &  4 &  5 & 6 &  7  &  8 &  9  \\ \hline
$\text{cond}(K^{(j)})$  & 1.00 & 1.00 &  1.00 &  1.00 &  1.01 & 1.03&  1.05  &  1.09 &  1.17
 \\ \hline
   \end{tabular}
     \end{center}

\begin{figure}[H]
\begin{center}
\begin{tabular}{c}
\epsfxsize=0.5\textwidth\epsffile{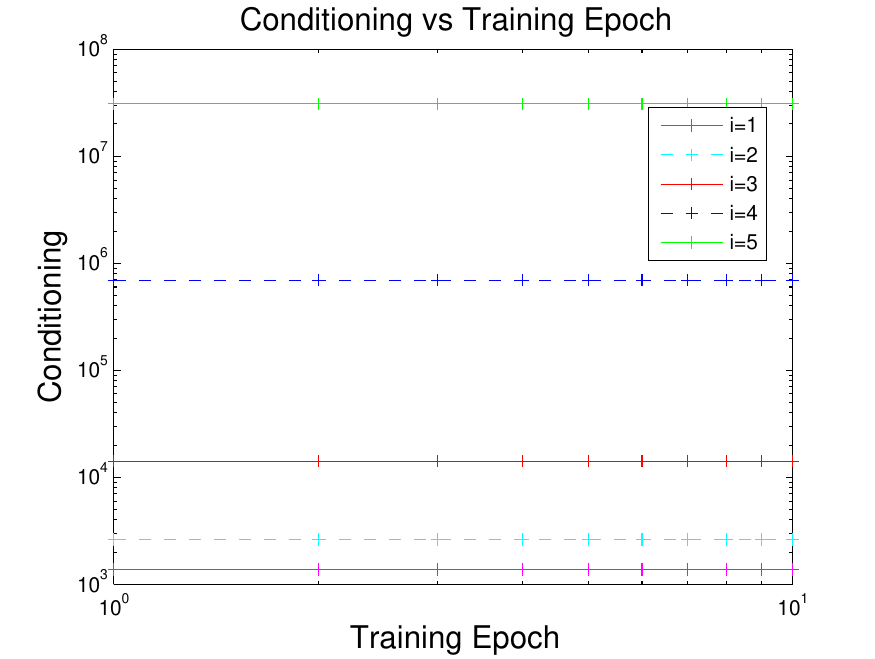}\\
\end{tabular}
\end{center}
 \caption{ Conditioning vs Training Epoch for some outer iterations.}
\label{2dseccond}
\end{figure}

\vskip 0.1in
Figure \ref{2dseccond} shows the conditioning of the intermediate systems in inner iterations associated with the linear system (\ref{galerlsq}) for some outer iterations.

Similarly to Figure  \ref{2dfircond}, the auxiliary plane wave least-squares systems solved during basis construction still suffers from ill-conditioned behavior \cite{ref21,hy2,peng,peng2} for many close directions when successive networks have increasing widths $n_i$.

\begin{figure}[H]
\begin{center}
\begin{tabular}{cc}
\epsfxsize=0.4\textwidth\epsffile{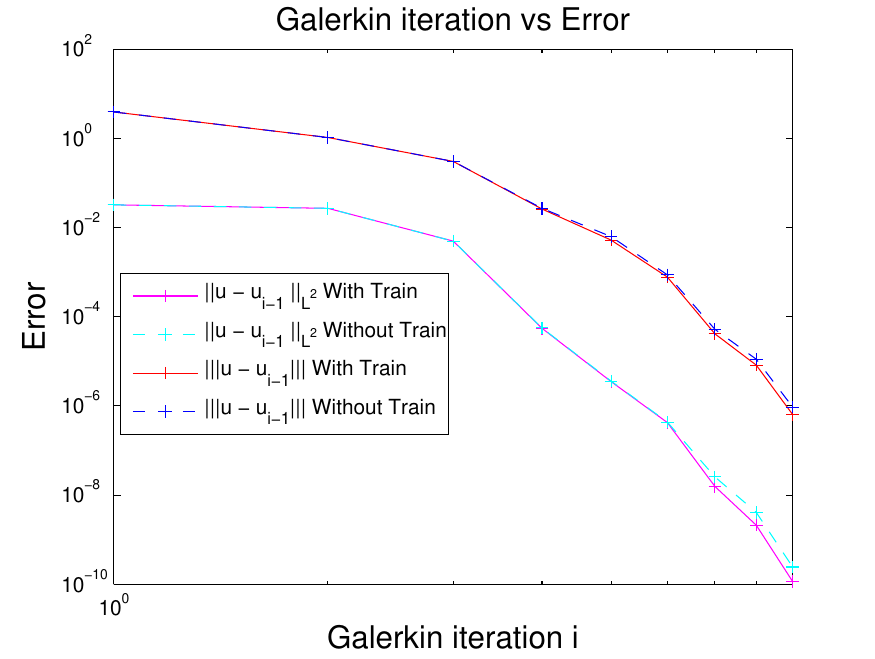} &
\epsfxsize=0.4\textwidth\epsffile{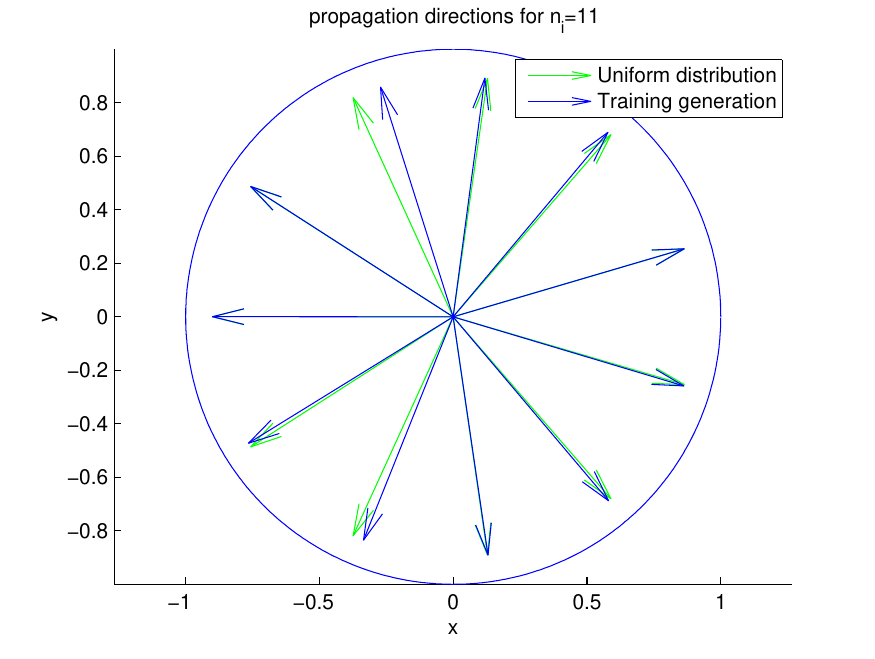} \\
\epsfxsize=0.4\textwidth\epsffile{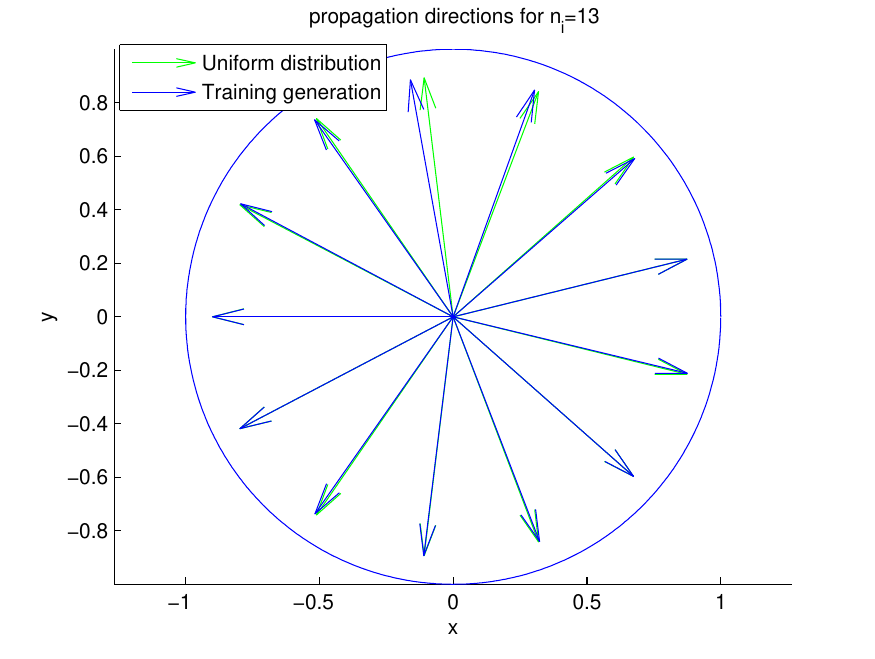} &
\epsfxsize=0.4\textwidth\epsffile{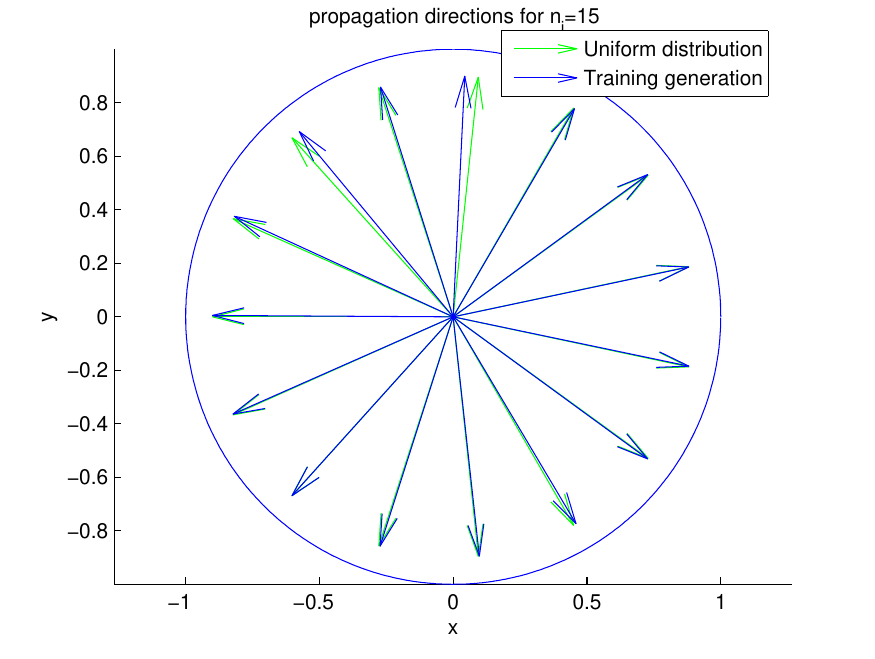} \\
\end{tabular}
\end{center}
 \caption{Left: The comparison of the errors when with training procedures and without the training procedure. Others: the sets of trained directions on the element centred at $(\frac{3}{32},\frac{21}{32})$ for different Galerkin iterations.}
\label{2dhelmcompare_dgf}
\end{figure}

To clarify the contribution of the direction-training step, we augment the discretization subspace without training procedures for the approximate maximizers in Algorithm 2. In the current situation, the proposed DGPWNN method degenerates into the same greedy residual-enrichment algorithm but with fixed (equispaced) directions. Figure \ref{2dhelmcompare_dgf} (Left-Up) shows the comparison of the errors when with training procedures and without the training procedure. Figure \ref{2dhelmcompare_dgf} (Others) shows the sets of trained directions on the element centred at $(\frac{3}{32},\frac{21}{32})$ for different Galerkin iterations.

The results shown in Figure \ref{2dhelmcompare_dgf} (Left) indicate that the accuracy of the approximation generated by the DGPWNN method with the training procedure is twice that of the approximation generated by the proposed method without the training procedure. Besides, we can see that, from Figure \ref{2dhelmcompare_dgf} (Others), the Adam-optimizer did substantially change the propagation directions.


%

\subsubsection{The case of large wavenumbers} \label{largewave}
In order to avoid the notorious ill-conditioning of the linear systems arising from Function dGLSQ-R based on the PWLS approach for the case of large wavenumbers whenever possible, a judicious hp-refinement strategy is employed to overcome the so-called pollution phenomenon \cite{ref21}.

\begin{figure}[H]
\begin{center}
\begin{tabular}{cc}
\epsfxsize=0.4\textwidth\epsffile{32pi_h16_galerkin.eps}&
\epsfxsize=0.4\textwidth\epsffile{32pi_h16_train.eps}\\
\epsfxsize=0.4\textwidth\epsffile{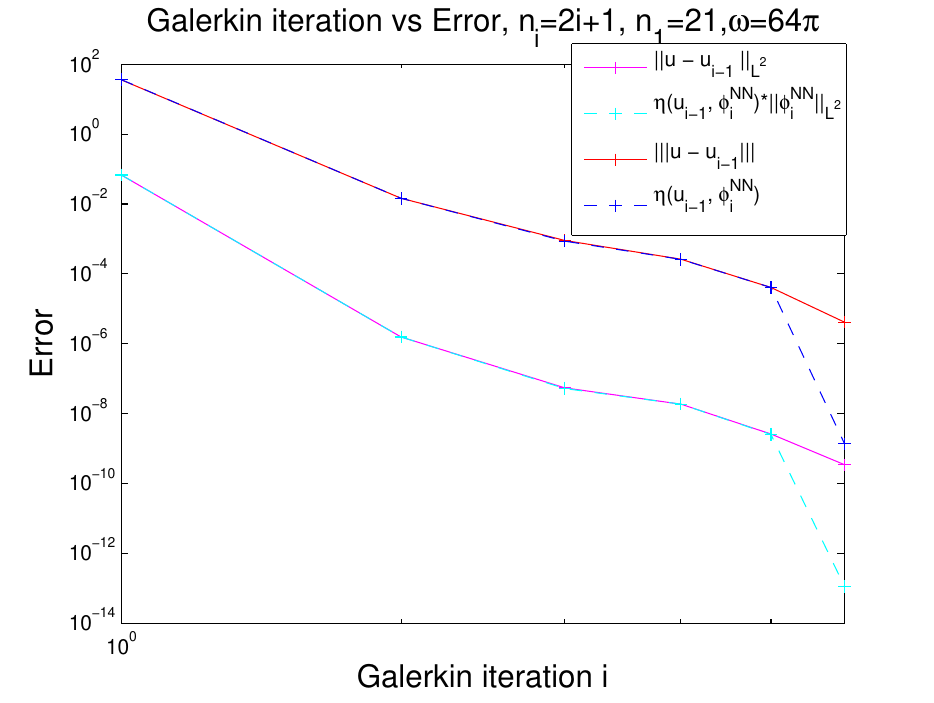}&
\epsfxsize=0.4\textwidth\epsffile{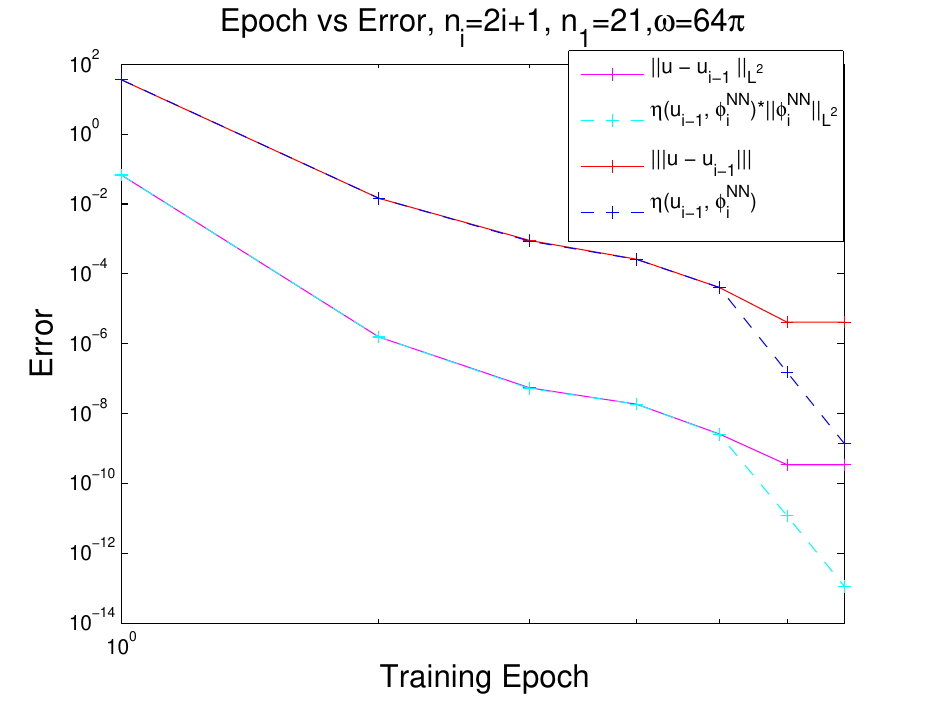}\\
\epsfxsize=0.4\textwidth\epsffile{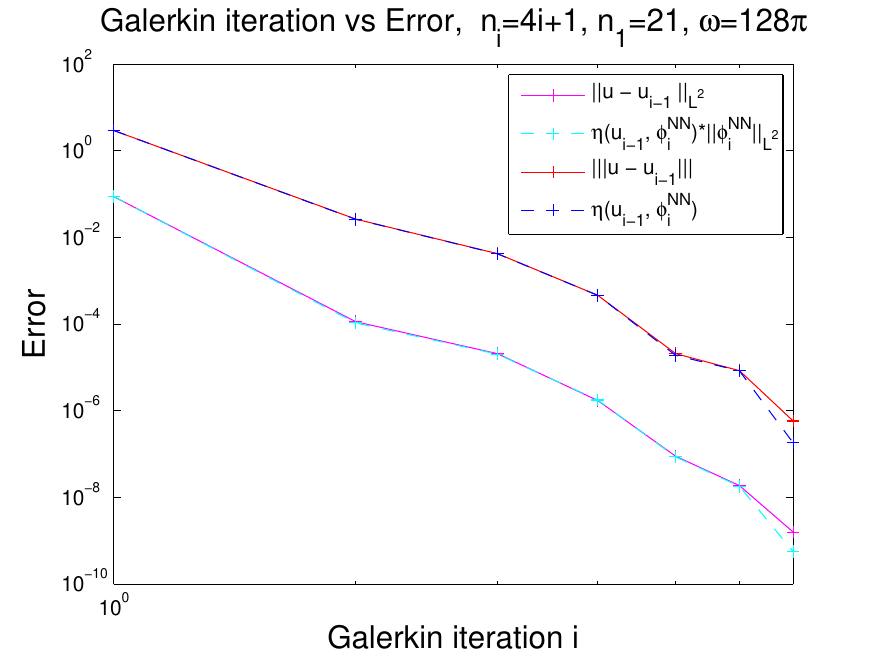}&
\epsfxsize=0.4\textwidth\epsffile{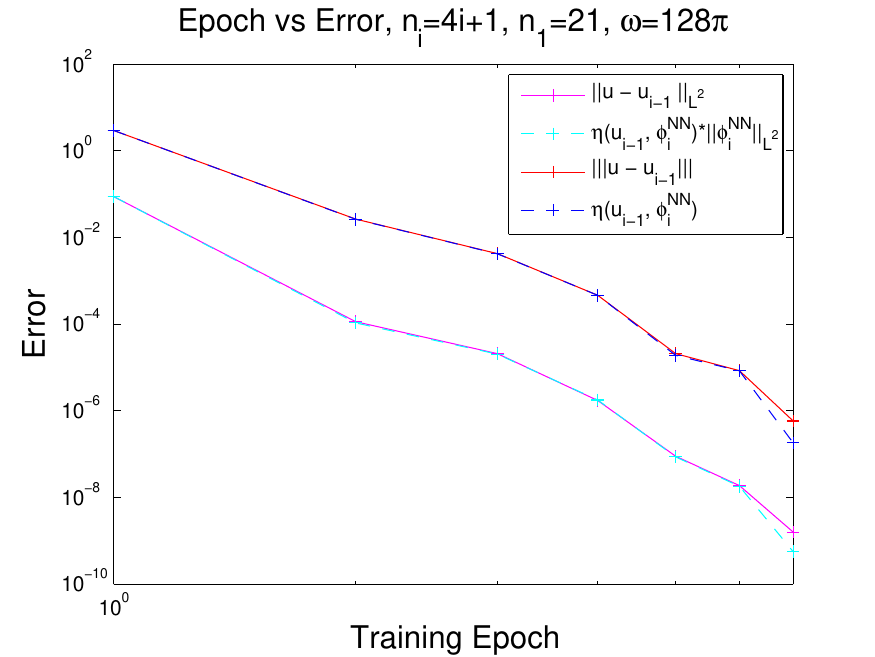}\\
\end{tabular}
\end{center}
 \caption{Displacement of a string (Isotropic Helmholtz equation in two dimension). (Left) Estimated and true error in the $L^2$ and energy norms at each Galerkin iteration, $n_i =2i+1$. (Right) The piecewise constant segments (solid lines) denote the true errors for each Galerkin iteration $i$, while the dashed blue line (resp., cyan line) denotes the progress of the loss function (resp., $L^2$ error estimator $\eta(u_{i-1}, \varphi_i^{NN}) ~| | \varphi_i^{NN}| |_{L^2}$) in approximating the true errors within each Galerkin iteration with $n_i =2i+1$. The x-axis thus denotes the cumulative training epoch over all Galerkin iterations. }
\label{2dmul_larwave_nn1}
\end{figure}

Set $h\approx \mathcal{O}(\frac{2\pi}{\omega})$ and $\varepsilon=10^{-8}$.
Figure \ref{2dmul_larwave_nn1} shows the true errors $||u-u_{i-1}||_{L^2}$ and $|||u-u_{i-1}|||$ at the end of each Galerkin iteration as well as the corresponding error estimates $\eta(u_{i-1}, \varphi_i^{NN}) ~| | \varphi_i^{NN}| |_{L^2}$) and $\eta(u_{i-1}, \varphi_i^{NN})$. We also provide the analogous results after each training epoch.

It can be seen that, only no more than ten outer discontinuous Galerkin iterations may guarantee the convergence of the discontinuous Galerkin plane wave neural network algorithm. Besides, we observe first that the initial Galerkin iterations reduce the error substantially, with further sharp decreases to a given precision in energy-norm approximation error of or $\mathcal{O}(10^{-6})$ and $L^2$ approximation error of or $\mathcal{O}(10^{-10})$. Further, it provides accurate approximations $u_i$ and estimators $\eta(u_{i-1}, \varphi_i^{NN}) ~| | \varphi_i^{NN}| |_{L^2}$ and $\eta(u_{i-1}, \varphi_i^{NN})$ of the true errors, with negligible discrepancies beginning to show in each Galerkin iteration except for the last iteration when $\omega$ is large. For example, for the last iteration when $\omega=64\pi$, the approximation error of $u_5$ in the energy norm and $L^2$ norm can reach the accuracy of $4.13e-6$, and $8.21e-9$, respectively, while two indicators for the errors in two norms are significantly smaller than the actual errors. The existence of this discrepancy is mainly due to the difficult characterization of high-dimensional reconstruction of errors caused by the nature of high-frequency problems of the acoustic equations.

\begin{figure}[H]
\begin{center}
\begin{tabular}{cc}
\epsfxsize=0.4\textwidth\epsffile{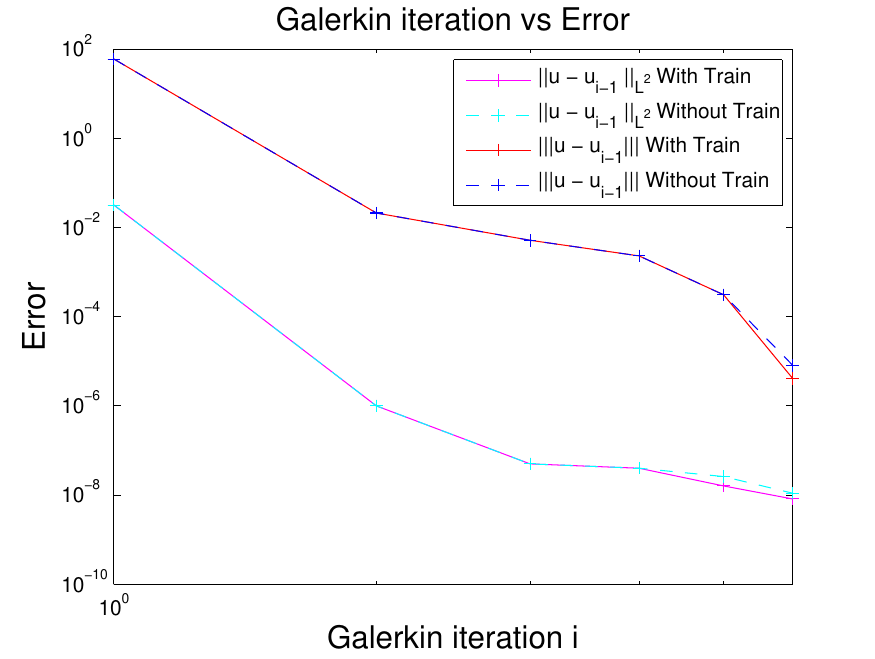} &
\epsfxsize=0.4\textwidth\epsffile{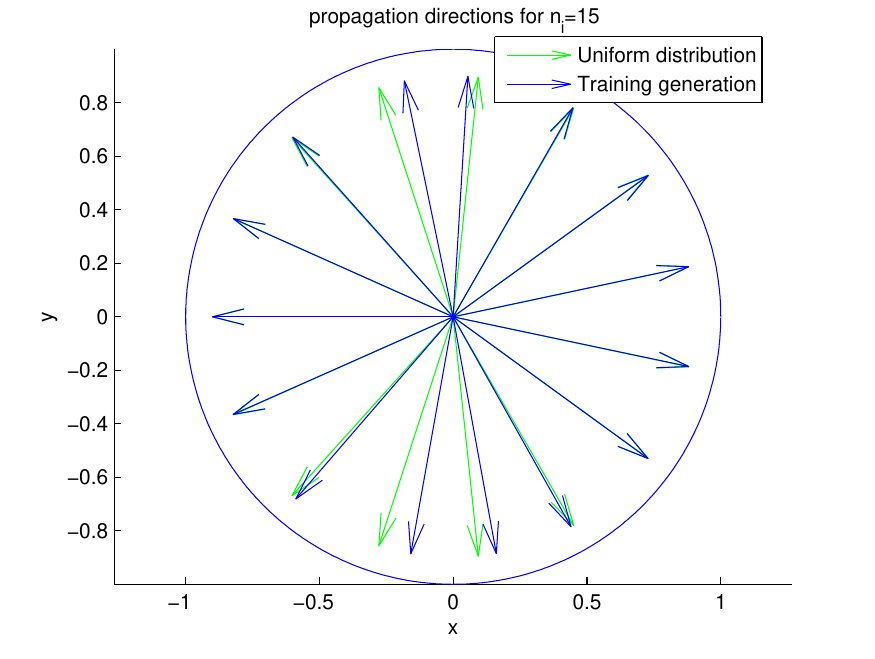} \\
\epsfxsize=0.4\textwidth\epsffile{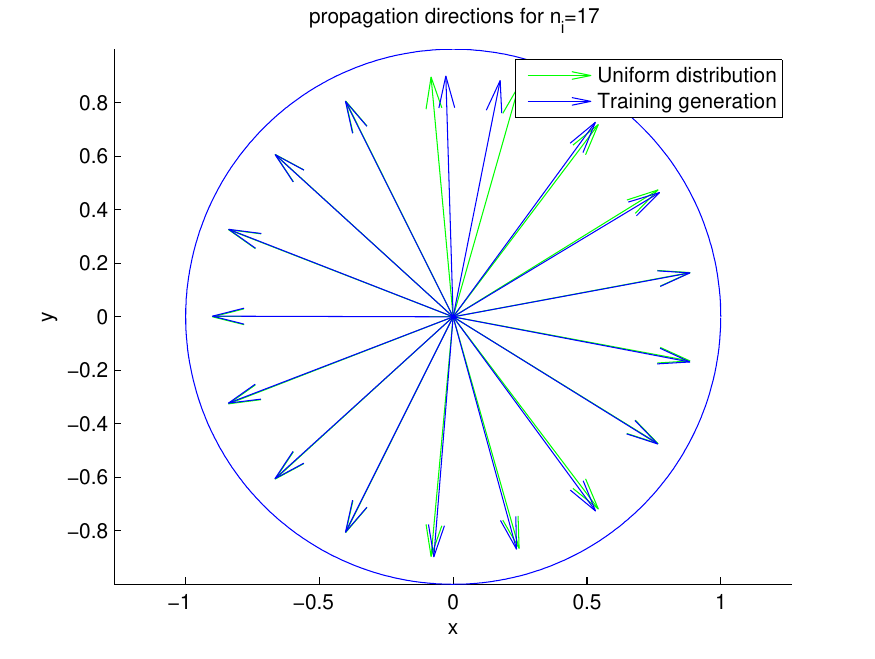} &
\epsfxsize=0.4\textwidth\epsffile{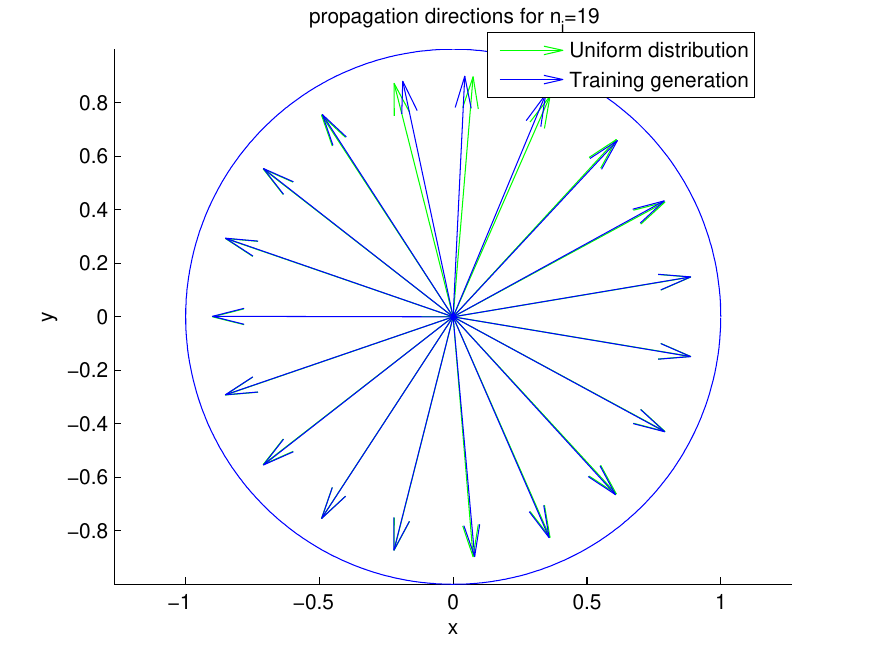} \\
\end{tabular}
\end{center}
 \caption{Left: The comparison of the errors when with training procedures and without the training procedure. Others: the sets of trained directions on the element centred at $(\frac{3}{64},\frac{3}{64})$ for different Galerkin iterations.}
\label{2dhighfrehelmcompare_dgf}
\end{figure}

Figure \ref{2dhighfrehelmcompare_dgf} (Left) shows the comparison of the errors when with training procedures and without the training procedure. Figure \ref{2dhighfrehelmcompare_dgf} shows the sets of trained directions on the element centred at $(\frac{3}{64},\frac{3}{64})$ for different Galerkin iterations.

The results shown in Figure \ref{2dhighfrehelmcompare_dgf} (Left) indicate that the accuracy of the approximation generated by the DGPWNN method with the training procedure is twice that of the approximation generated by the proposed method without the training procedure. Besides, we can see that, from Figure \ref{2dhighfrehelmcompare_dgf} (Others), the Adam-optimizer did substantially change the propagation directions.

\begin{figure}[H]
\begin{center}
\begin{tabular}{cc}
\epsfxsize=0.4\textwidth\epsffile{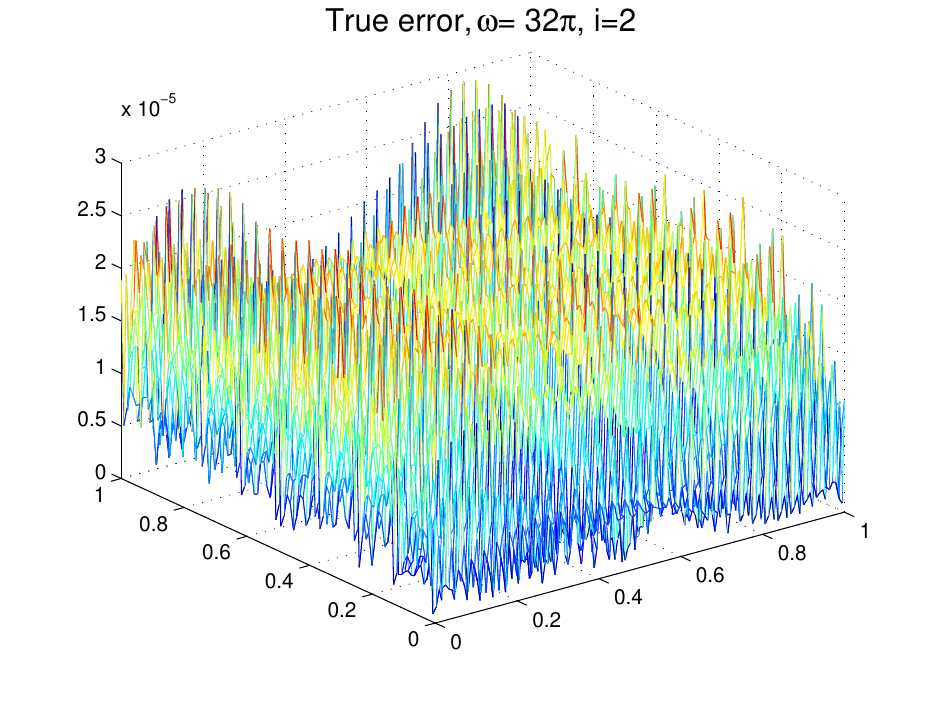}&
\epsfxsize=0.4\textwidth\epsffile{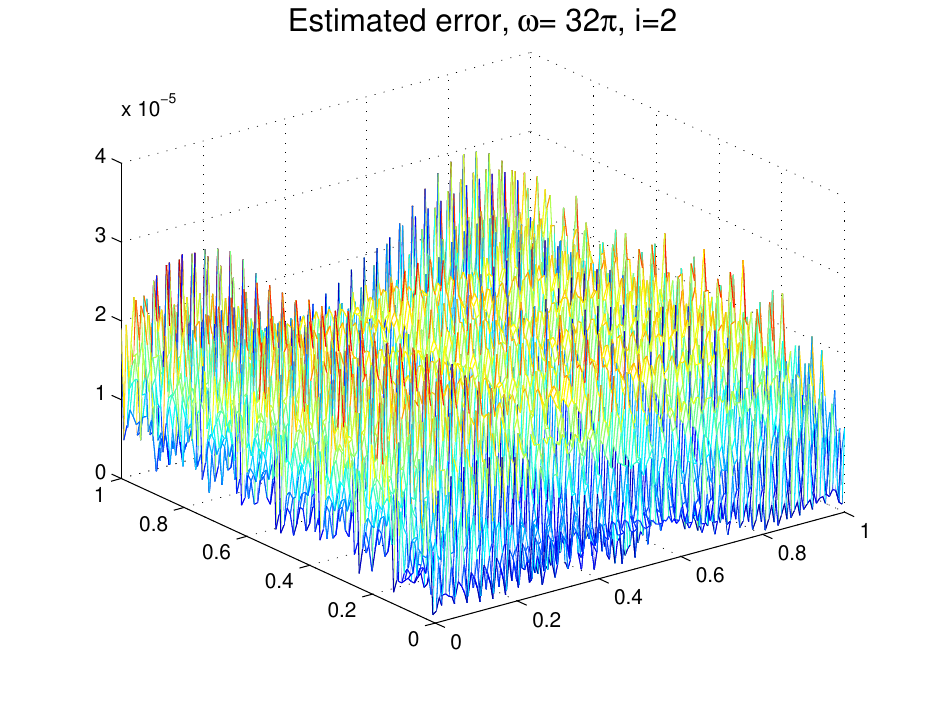}\\
\epsfxsize=0.4\textwidth\epsffile{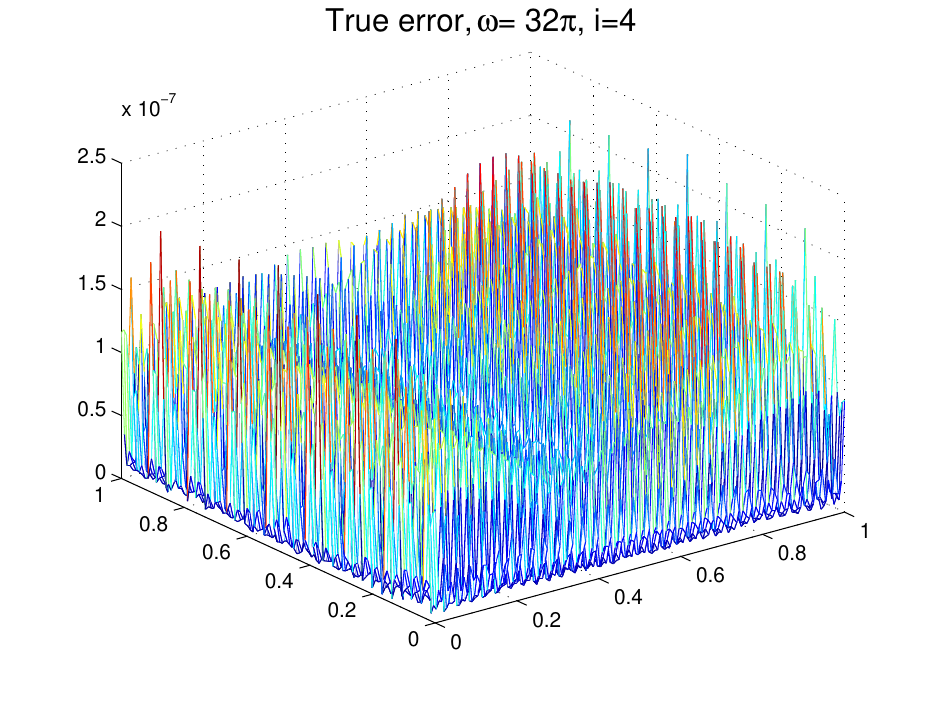}&
\epsfxsize=0.4\textwidth\epsffile{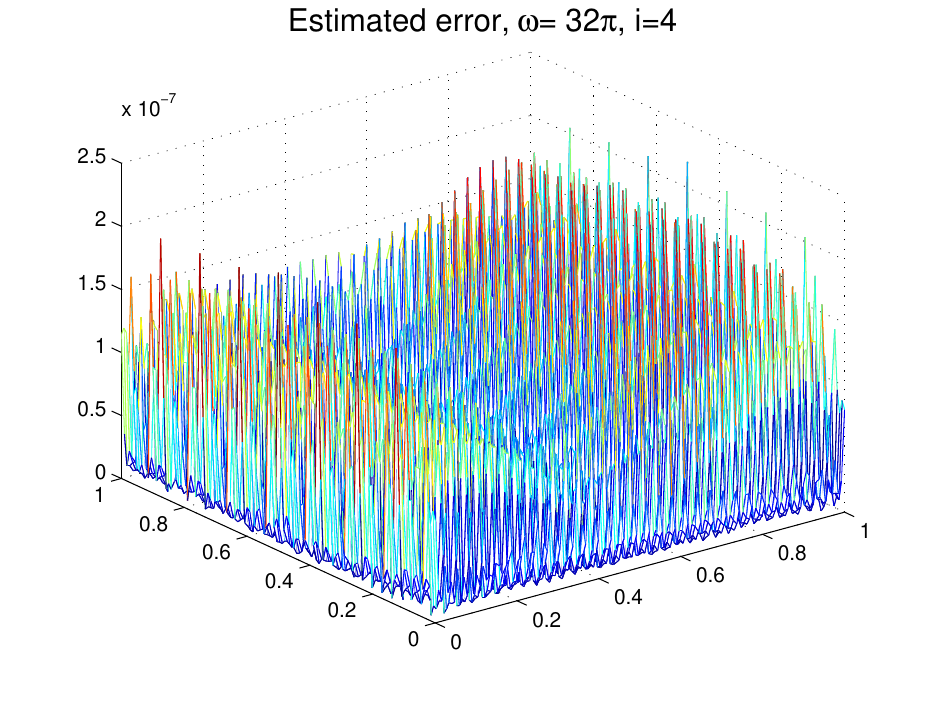}\\
\epsfxsize=0.4\textwidth\epsffile{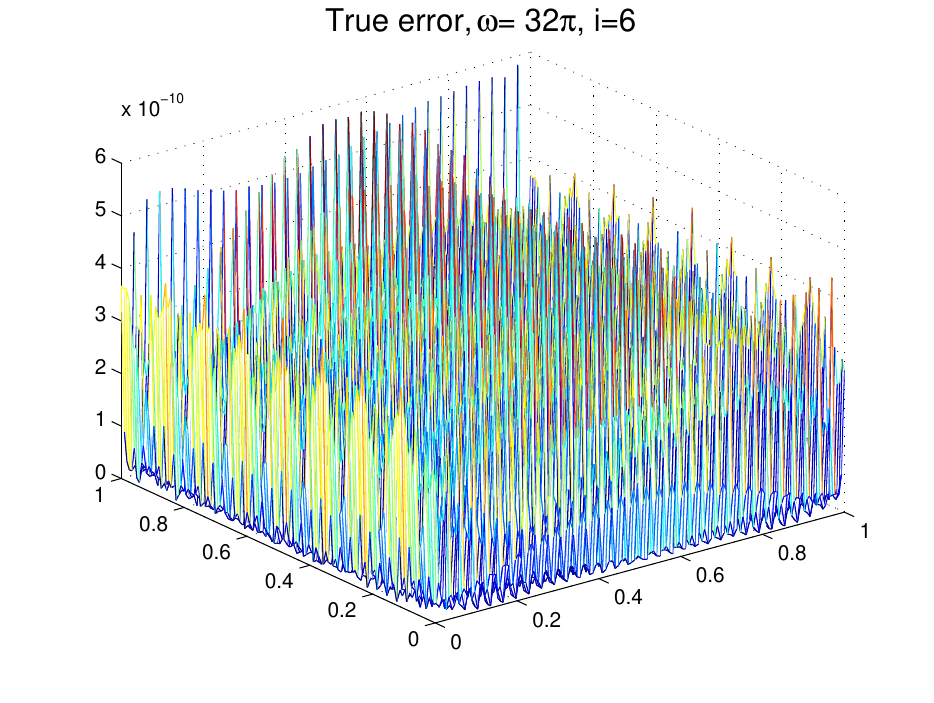}&
\epsfxsize=0.4\textwidth\epsffile{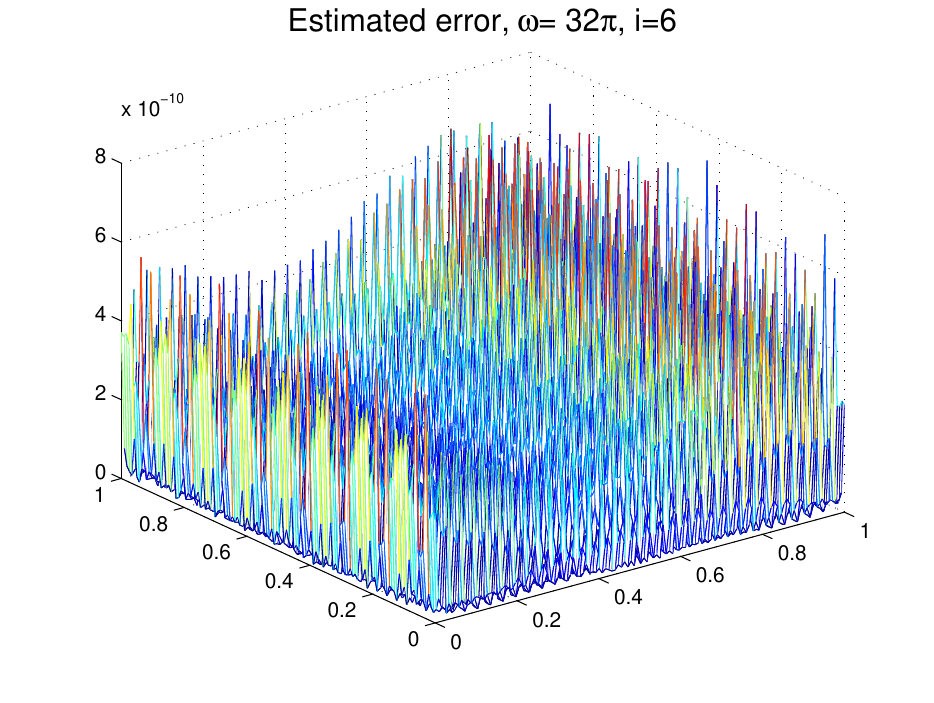}\\
\end{tabular}
\end{center}
  \caption{Displacement of a membrane (Helmholtz equation in two dimensions). Exact error $|u-u_{i-1}|$ (left column) and approximate error $\eta(u_{i-1}, \varphi_i^{NN}) ~|\varphi_i^{NN}|$ (right column) for i = 2, 4, 6. }
\label{2dhelarge}
\end{figure}

Figure \ref{2dhelarge} shows the exact error $|u-u_{i-1}|$ as well as the maximizer $\eta(u_{i-1}, \varphi_i^{NN}) ~|\varphi_i^{NN}|$ at several stages of the algorithm. Here we set $\omega=32\pi$.

It is apparent that initially, the low frequency components of the error are learned, with later iterations learning the high frequency error components. Besides, the pointwise error is consistent with the Galerkin iteration error in Figure \ref{2dmul_larwave_nn1}.

\subsubsection{Comparison with the PWLS method} \label{comdgnnpwls}
\begin{figure}[H]
\begin{center}
\begin{tabular}{cc}
\epsfxsize=0.4\textwidth\epsffile{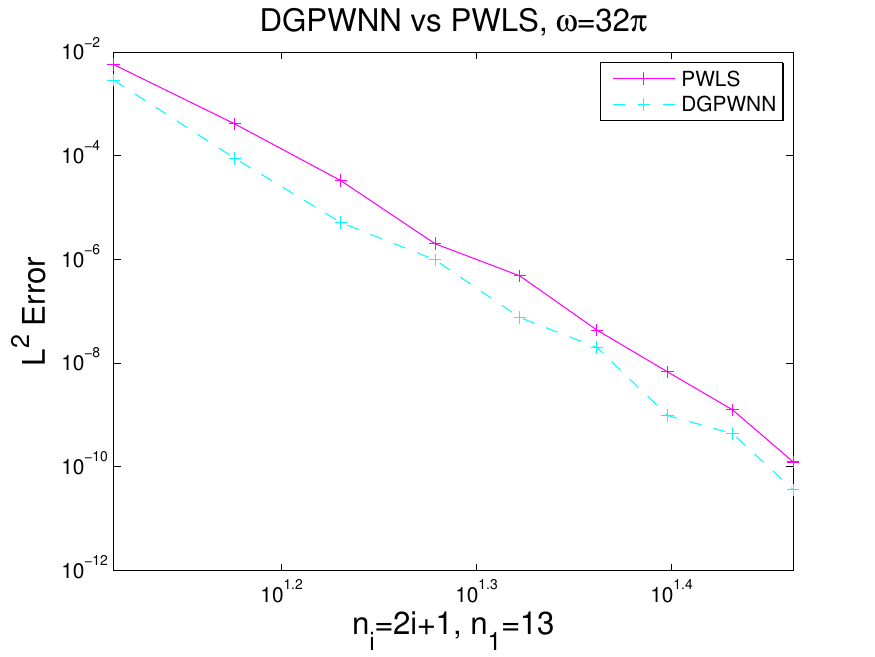}&
\epsfxsize=0.4\textwidth\epsffile{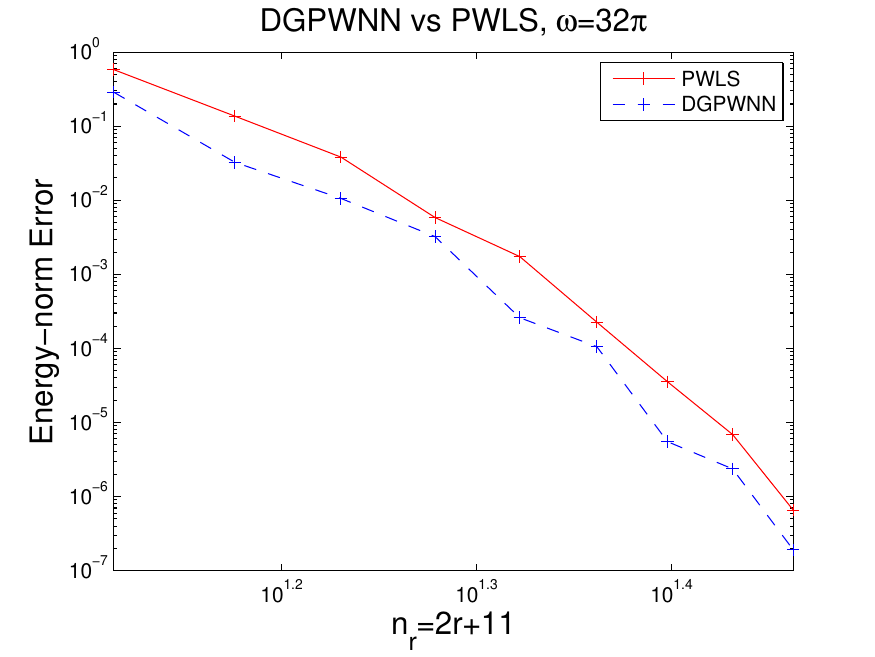}\\
\epsfxsize=0.4\textwidth\epsffile{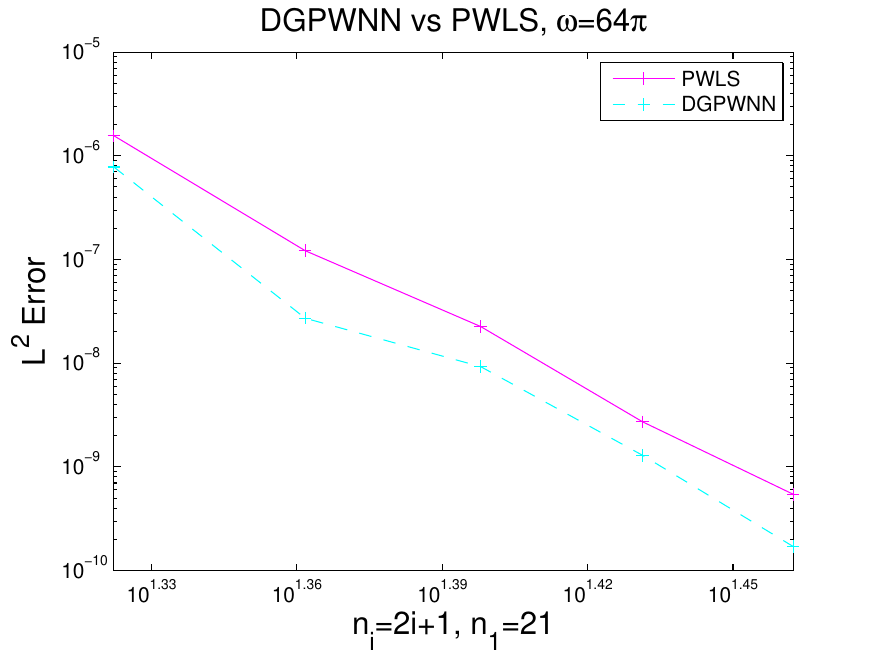}&
\epsfxsize=0.4\textwidth\epsffile{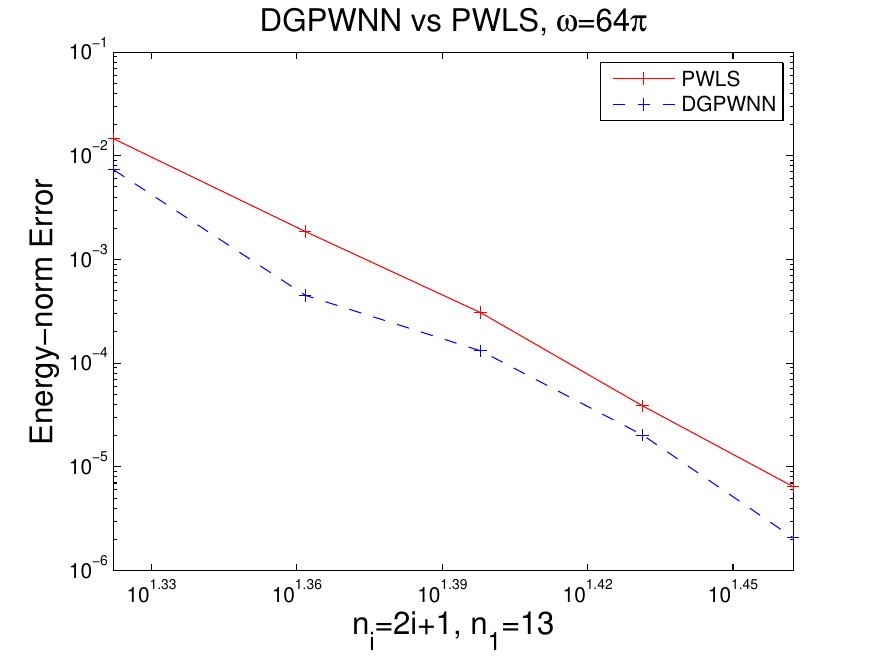}\\
\end{tabular}
\end{center}
 \caption{Comparison between the DGPWNN and the PWLS method. (Left) True errors in the $L^2$ norm with respect to $n_i =2i+1$. (Right) True errors in the  energy norm with respect to $n_i =2i+1$. }
\label{2dhelcompari}
\end{figure}

In this section, the proposed DGPWNN is compared with the PWLS method \cite{hy}. Figure \ref{2dhelcompari} shows the comparison of the errors in the two norms of the resulting approximations. In particular, for the DGPWNN, $n_i$ on the horizontal axis denotes the width of the discontinuous network for the $i-$th Galerkin iteration. While, for the PWLS, $n_i$ denotes the number of degrees of freedom per each element.

It is worth emphasizing that the errors of the approximations generated by the DGPWNN method are smaller than that generated by the PWLS method for every $n_i$. In particular, at the last iteration, specifically for the case of $\omega=32\pi, ~n_i=29$, the errors of the approximations in $L^2$ norm and the energy norms generated by the DGPWNN method reach to $3.61e-11$ and $1.91e-7$, respectively, while the errors of the approximations in $L^2$ norm and the energy norms generated by the PWLS method reach to $1.24e-10$ and $6.55e-7$, respectively.

Table \ref{computetime} shows the overall run time of Function AugmentBasis($u_j$) for each discontinuous Galerkin iteration $j$, and the computing time of linear system generated by the PWLS method with the number of plane wave basis functions $n_j$, respectively. The time unit is seconds.

\vskip 0.1in
\begin{center}
       \tabcaption{Comparison of computing time between DGPWNN and PWLS.}
\label{computetime}
       \vskip 0.1in
\begin{tabular}{|c|c|c|c|c|c|} \hline
   $n_j$  &  17 &  19 &  21 & 23 &  25     \\ \hline
\text{DGPWNN}  &   3.29e+3    & 4.48e+3  &  3.11e+3  & 4.50e+3 &  6.31e+3   \\ \hline
\text{PWLS}  &   1.44e+3 & 2.34e+3  & 3.10e+3  & 4.49e+3 &  6.49e+3   \\ \hline
   \end{tabular}
     \end{center}

\vskip 0.1in
 It can be seen from Table \ref{computetime} that, when $j$ increases, the computational cost associated with obtaining $c$ for the entire sequence of DGPWNN is comparable to the cost of the PWLS method with a single network of width $n_j$.
This is mainly because the main computational expense in Algorithm 1 is the computation of the entire sequence of activation coefficients $c^{(j)}$ defined by the linear system (\ref{galerlsq}) for discontinuous Galerkin iteration $j$ (refer to \cite[Section 2.4.1]{AD}), while the number of training epoches in Function AugmentBasis($u_j$) is decreasing sharply. Of course, the improvement of accuracy comes at the cost of a certain computational complexity.

\subsection{An example of three-dimensional Helmholtz equation}
The following test problem consists of a point source and associated
boundary conditions for homogeneous Helmholtz equations (see
\cite{HKM}):
\begin{equation}
\begin{split}
&u(r,r_0)= {1\over 4\pi} {e^{i\omega|r-r_0|}\over |r-r_0|} ~~ \text{in} ~~\Omega,  \\
&{\partial u\over \partial {\bf n}}+i\omega u=g \quad
\text{over}\quad\partial\Omega, \\
\end{split}\label{fig1}
\end{equation}
in a cubic computational domain $\Omega=[0,1]\times[0,1]\times[0,1]$. The location of the source is off-centred at
$r_0=(-1,-1,-1)$ and $r=(x,y,z)$ is an observation point. Note that the
analytic solution of the Helmholtz equation with such boundary
condition has a singularity at $r=r_0$.

Set $h\approx \mathcal{O}(\frac{\pi}{\omega})$, $\varepsilon=10^{-6}$, and $n_i=2n_1^2$ satisfying $n_1 \geq 4$. 

\begin{figure}[H]
\begin{center}
\begin{tabular}{cc}
\epsfxsize=0.4\textwidth\epsffile{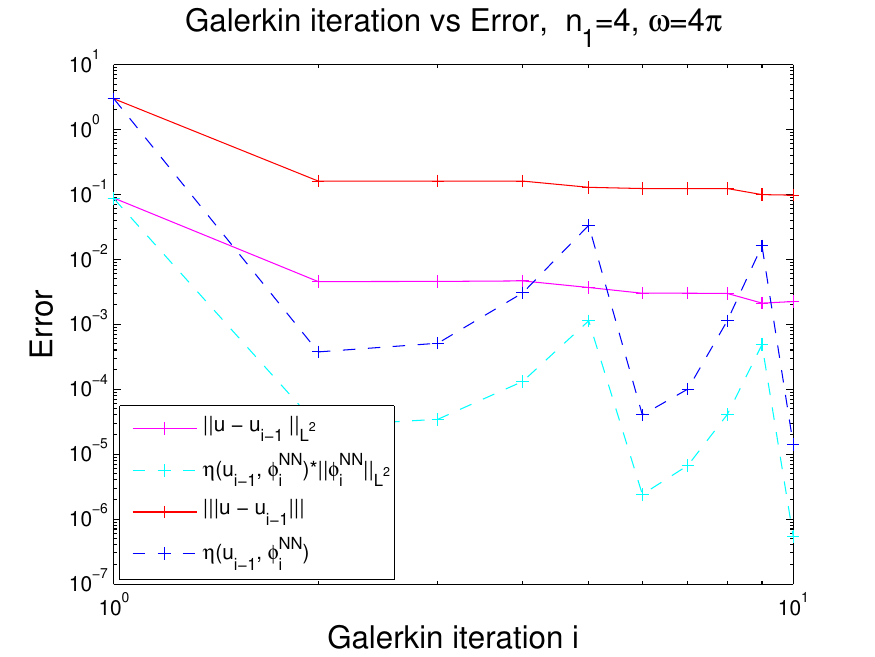}&
\epsfxsize=0.4\textwidth\epsffile{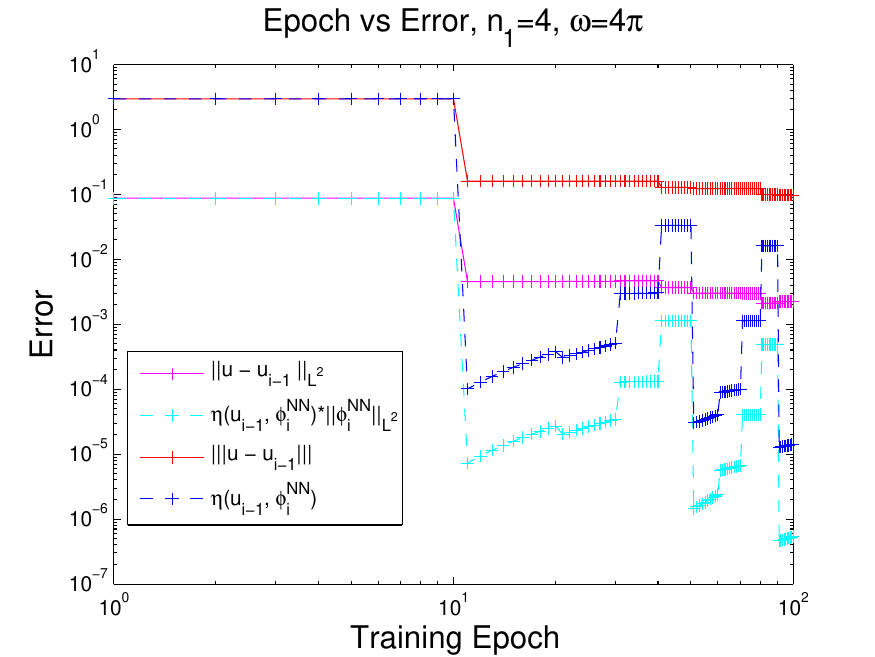}\\
\epsfxsize=0.4\textwidth\epsffile{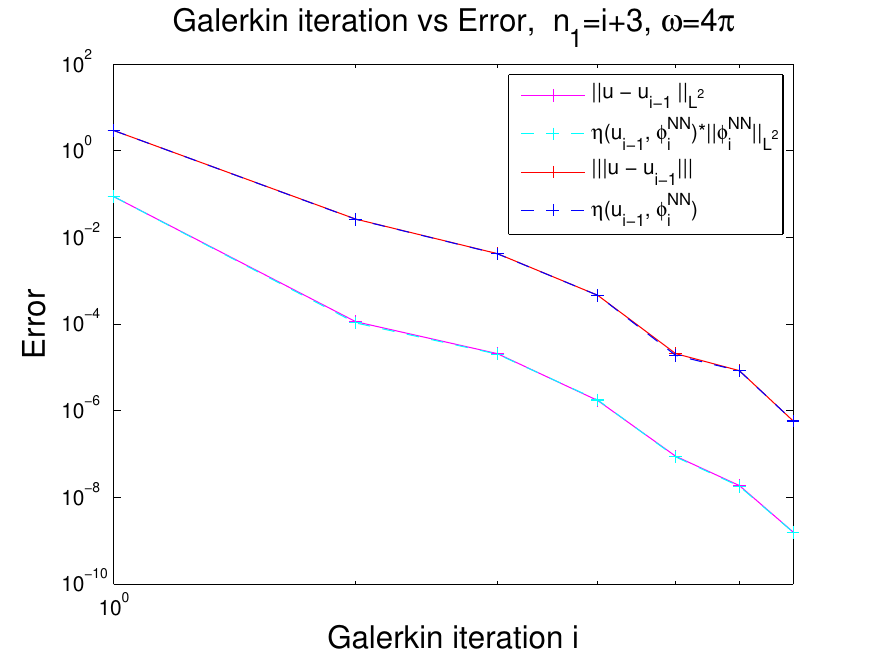}&
\epsfxsize=0.4\textwidth\epsffile{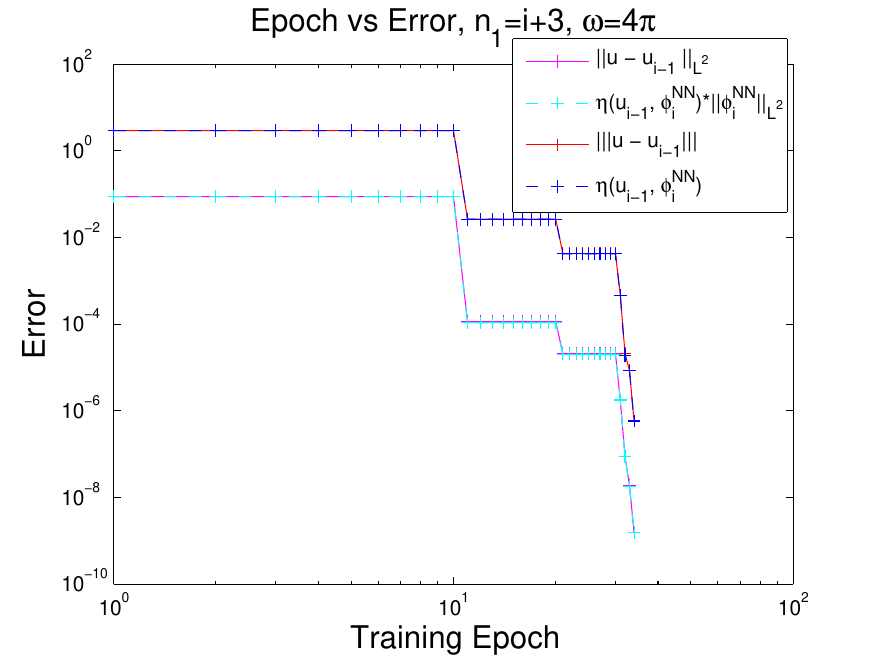}\\
\end{tabular}
\end{center}
 \caption{Displacement of a string (Isotropic Helmholtz equation in three dimension). (Left-Up) Estimated and true error in the $L^2$ and energy norms at each Galerkin iteration, $n_1 = 4$ for all i. (Right-Up) The piecewise constant segments (solid lines) denote the true errors for each Galerkin iteration $i$, while the dashed blue line (resp., cyan line) denotes the progress of the loss function (resp., $L^2$ error estimator $\eta(u_{i-1}, \varphi_i^{NN}) ~| | \varphi_i^{NN}| |_{L^2}$) in approximating the true errors within each Galerkin iteration with $n_1=4$. The x-axis thus denotes the cumulative training epoch over all Galerkin iterations. (Left-Bottom) and (Right-Bottom) The analogous results for $n_1 =i+3$. }
\label{3dhenn1}
\end{figure}

Figure \ref{3dhenn1} shows the true errors $||u-u_{i-1}||_{L^2}$ and $|||u-u_{i-1}|||$ at the end of each Galerkin iteration as well as the corresponding error estimates $\eta(u_{i-1}, \varphi_i^{NN}) ~| | \varphi_i^{NN}| |_{L^2}$ and $\eta(u_{i-1}, \varphi_i^{NN})$. We also provide the analogous results after each training epoch.

We observe first that the initial Galerkin iterations reduce the error substantially, with further decreases beginning to stagnate and increasingly poor matching of the estimator $\eta(u_{i-1}, \varphi_i^{NN})$ of the true error as $i$ grows larger whenever $n_i$ is constant, owing to the poor approximations $u_i$ to the analytic solution, which further, to the contrary, coincides with the validity of (\ref{equiestimator}) and (\ref{l2indi}); stagnation is not an issue when $n_i$ is increased with $i$. In particular, compared with the case that $n_i$ is constant, the case that $n_i$ is increased with $i$ provides more accurate approximations $u_i$ and estimators $\eta(u_{i-1}, \varphi_i^{NN}) ~| | \varphi_i^{NN}| |_{L^2}$ and $\eta(u_{i-1}, \varphi_i^{NN})$ of the true errors, with negligible discrepancies beginning to show in each Galerkin iteration. Besides, for the case of adaptive choice of the number $n_i$, only no more than ten outer discontinuous Galerkin iterations may guarantee the convergence of our plane wave neural network algorithm.

Next, Table \ref{3dhelmcond_multidomain} shows the condition number of the matrix $K^{(j)}$ associated with the linear system (\ref{dgappro}) at each discontinuous Galerkin iteration $j$ when $n_1 =j+3$. We observe that the condition number is approximately unity as expected.

\vskip 0.1in
\begin{center}
       \tabcaption{Condition number for the discontinuous Galerkin matrix at each iteration $j$.}
\label{3dhelmcond_multidomain}
       \vskip 0.1in
\begin{tabular}{|c|c|c|c|c|c|c|c|} \hline
   $j$  & 1 & 2 &  3 &  4 &  5 & 6 &  7    \\ \hline
$\text{cond}(K^{(j)})$  & 1.00 & 1.00 &  1.01 &  1.05 &  1.07 & 1.09&  1.15
 \\ \hline
   \end{tabular}
     \end{center}

\vskip 0.1in

\begin{figure}[H]
\begin{center}
\begin{tabular}{cc}
\epsfxsize=0.4\textwidth\epsffile{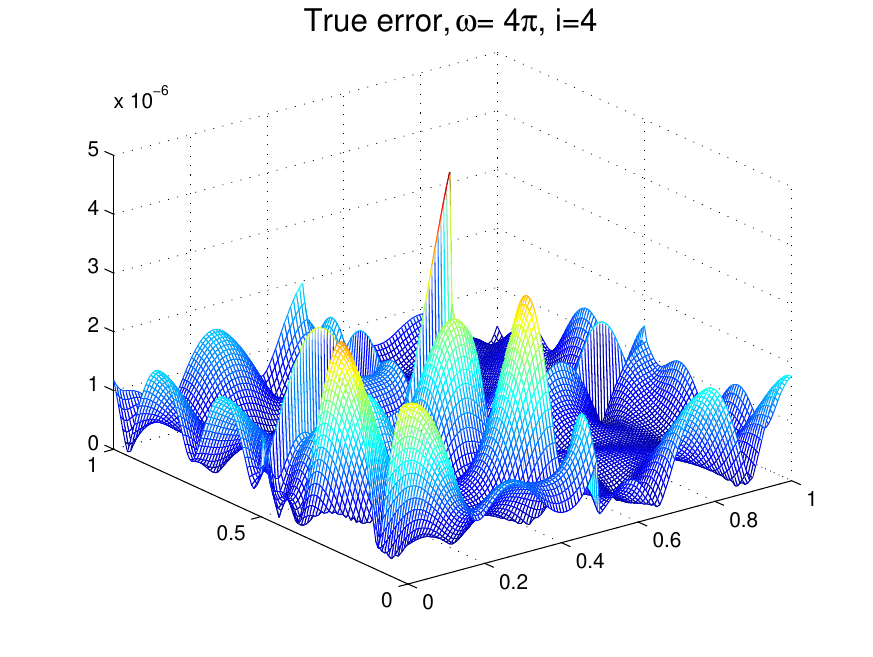}&
\epsfxsize=0.4\textwidth\epsffile{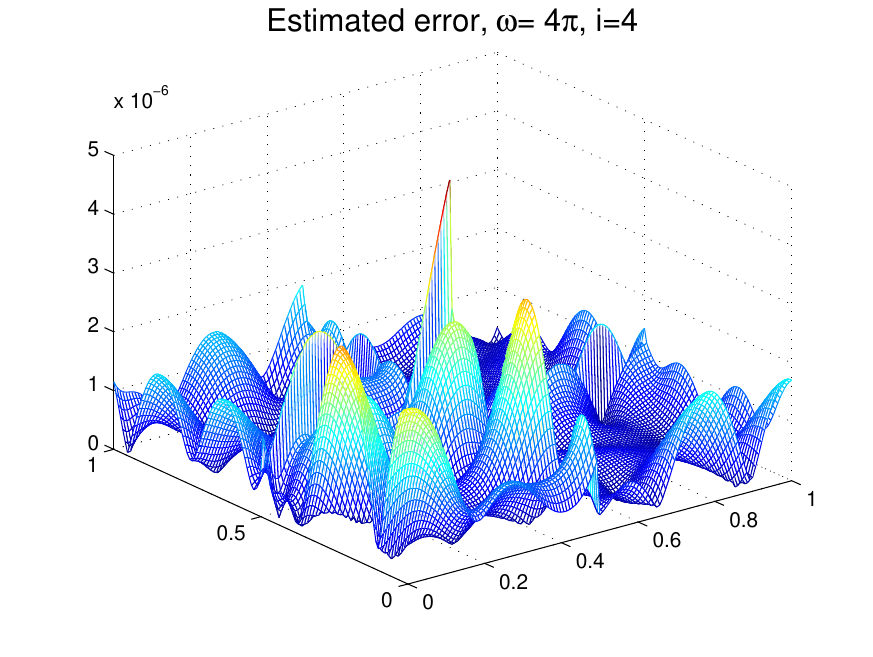}\\
\epsfxsize=0.4\textwidth\epsffile{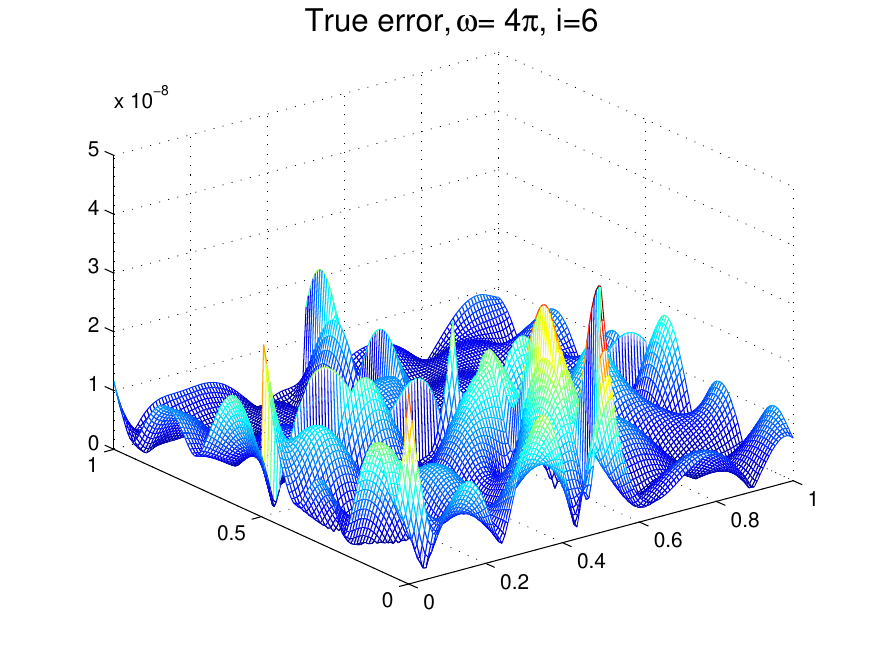}&
\epsfxsize=0.4\textwidth\epsffile{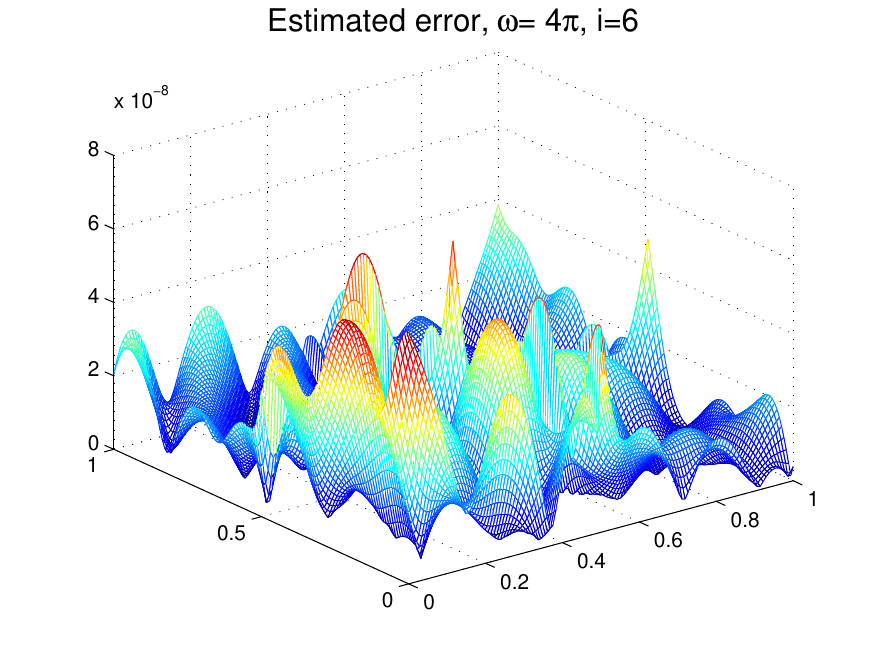}\\
\end{tabular}
\end{center}
 \caption{Displacement of a membrane (Helmholtz equation in three dimensions). Exact error $|u-u_{i-1}|$ (left column) and approximate error $\eta(u_{i-1}, \varphi_i^{NN}) ~|\varphi_i^{NN}|$ (right column) for i = 4, 6. }
\label{3dhenn2}
\end{figure}

Figure \ref{3dhenn2} shows the exact error $|u-u_{i-1}|$ as well as the maximizer $\eta(u_{i-1}, \varphi_i^{NN}) ~|\varphi_i^{NN}|$ at several stages of the algorithm for the case of cross-section $z=0.5$.

It is apparent that initially, the low frequency components of the error are learned, with later iterations learning the high frequency error components, and the pointwise error is consistent with the Galerkin iteration error in Figure \ref{3dhenn1}.

\subsection{An example of three-dimensional Maxwell's equations}

We compute the electric field due to an electric dipole source at
the point \( x_0=(0.6,0.6,0.6)\). The dipole point source can
be defined as the solution of a homogeneous Maxwell system
(\ref{maxeq}). The exact solution of the problems is
\begin{equation}
{E}_{\text{ex}}=-\text{i}\omega I \phi({ x},{ x}_0){ a}
+\frac{I}{{ i}\omega\varepsilon} \nabla(\nabla\phi\cdot{ a}),
\end{equation}
where
$$ \phi({x},{ x}_0) = \frac{\text{exp}(\text{i}\omega\sqrt{\varepsilon}|{ x}-{ x}_0|)}{4\pi|{x}-{ x}_0|}$$
and \(\Omega=[-0.5,0.5]^3\).

 Set $h = \frac{4\pi}{\omega}$ and $\varepsilon=10^{-6}$. Set the number of basis functions per element as $n_i=4n_1^2$ satisfying $n_1 \geq 4$, which indicates that the number of total degrees of freedom for linear systems (\ref{galerlsq}) is equal to $h^{-3}\times n_i$. Set the number of training epoches as 20 in Algorithm 2. The total computational time is about $2.67e+3$ seconds when $\omega= 8\pi$ and $n_i=4(i+3)^2$.

\begin{figure}[H]
\begin{center}
\begin{tabular}{cc}
\epsfxsize=0.4\textwidth\epsffile{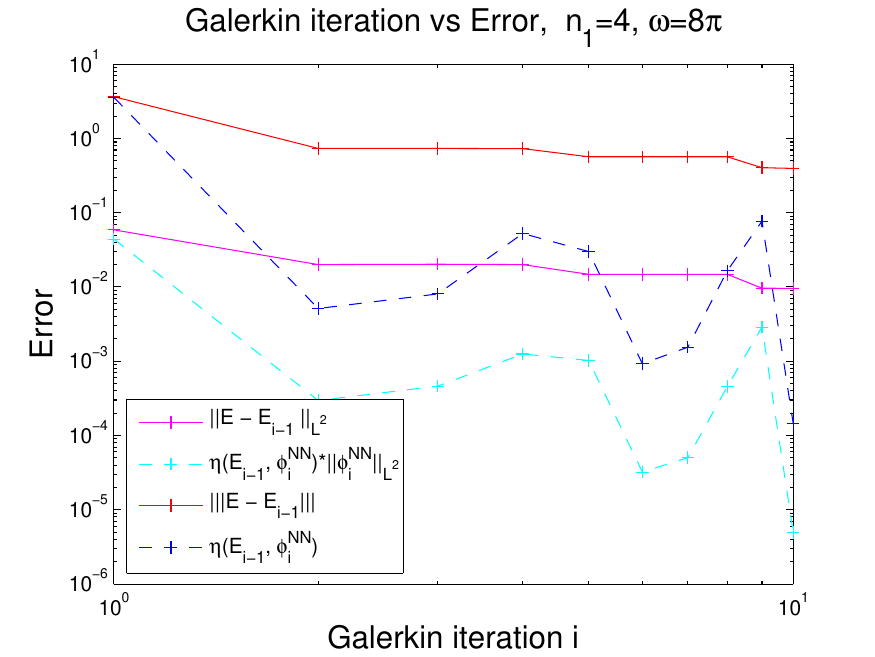}&
\epsfxsize=0.4\textwidth\epsffile{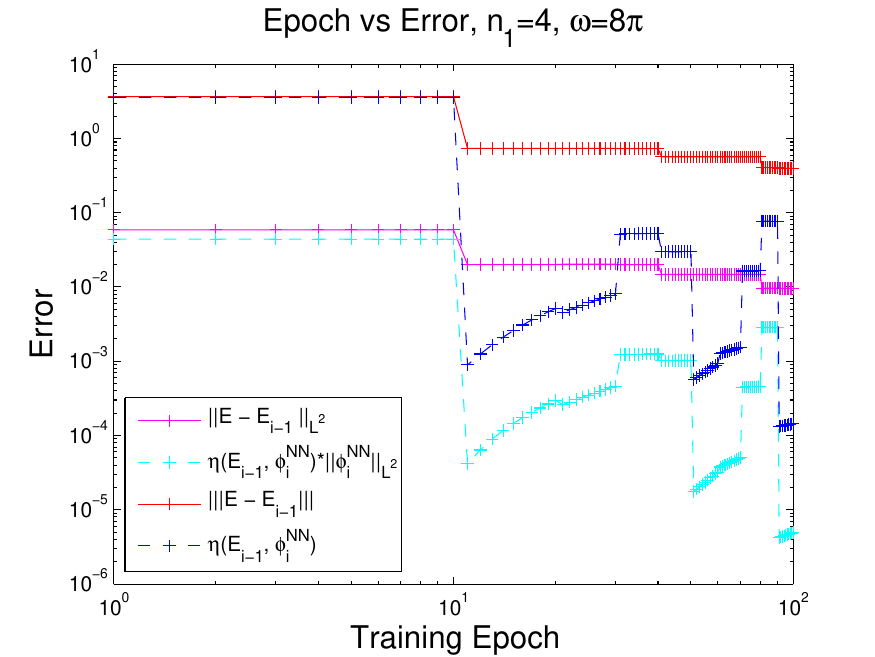}\\
\epsfxsize=0.4\textwidth\epsffile{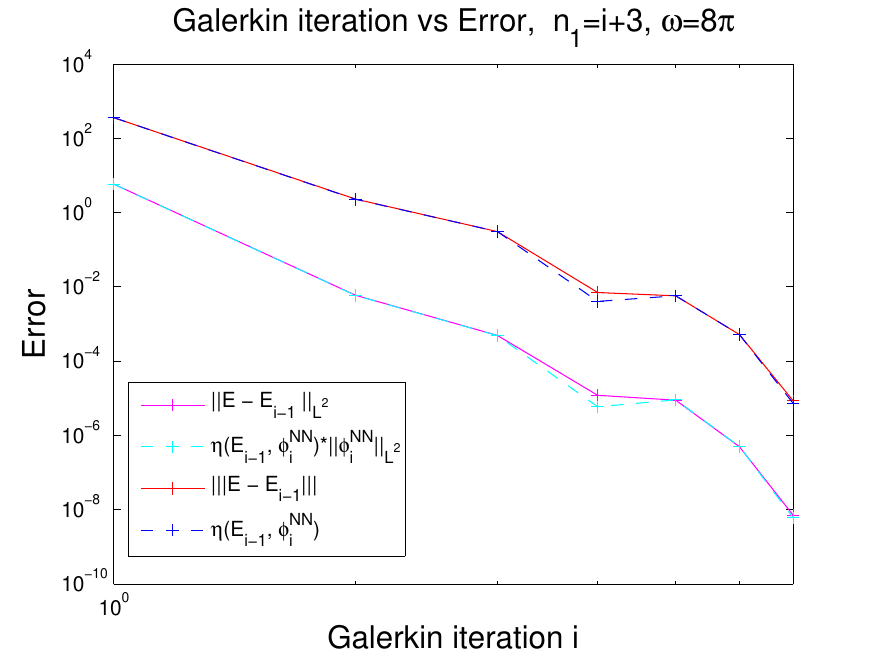}&
\epsfxsize=0.4\textwidth\epsffile{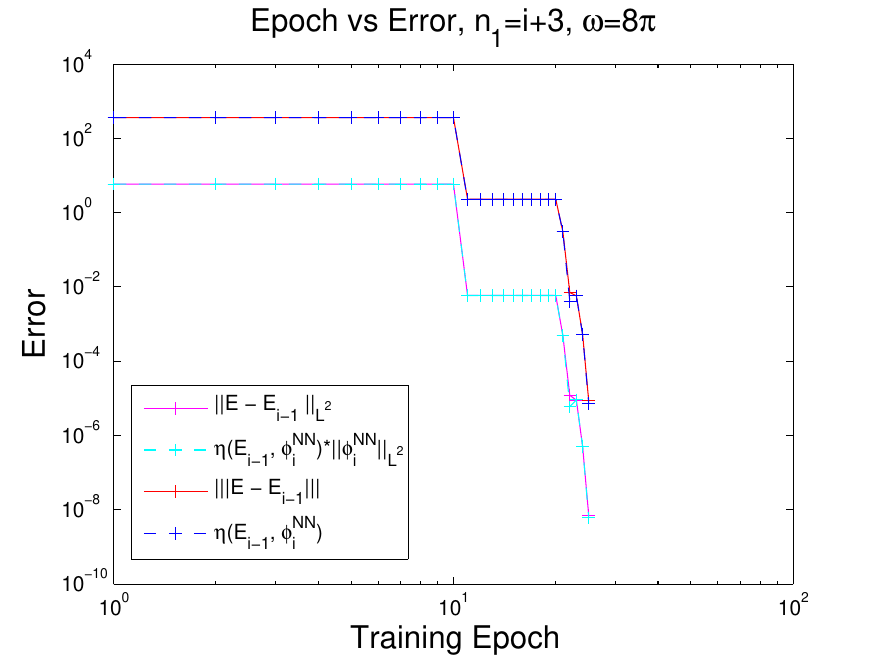}\\
\end{tabular}
\end{center}
 \caption{Displacement of a string (Maxwell's equations in three dimension). (Left-Up) Estimated and true error in the $L^2$ and energy norms at each Galerkin iteration, $n_1 = 4$ for all i. (Right-Up) The piecewise constant segments (solid lines) denote the true errors for each Galerkin iteration $i$, while the dashed blue line (resp., cyan line) denotes the progress of the loss function (resp., $L^2$ error estimator $\eta(E_{i-1}, \varphi_i^{NN}) ~| | \varphi_i^{NN}| |_{L^2}$) in approximating the true errors within each Galerkin iteration with $n_1=4$. The x-axis thus denotes the cumulative training epoch over all Galerkin iterations. (Left-Bottom) and (Right-Bottom) The analogous results for $n_1 =i+3$. }
\label{3dmaxwell}
\end{figure}

Figure \ref{3dmaxwell} shows the true errors $||E-E_{i-1}||_{L^2}$ and $|||E-E_{i-1}|||$ at the end of each Galerkin iteration as well as the corresponding error estimates $\eta(E_{i-1}, \varphi_i^{NN}) ~| | \varphi_i^{NN}| |_{L^2}$ and $\eta(E_{i-1}, \varphi_i^{NN})$. We also provide the analogous results after each training epoch.

We observe first that the initial Galerkin iterations reduce the error substantially, with further decreases beginning to stagnate and increasingly poor matching of the estimator $\eta(E_{i-1}, \varphi_i^{NN})$ of the true error as $i$ grows larger whenever $n_i$ is constant, owing to the poor approximations $u_i$ to the analytic solution,
  which further, to the contrary, coincides with the validity of (\ref{equiestimator}) and (\ref{l2indi}); stagnation is not an issue when $n_i$ is increased with $i$. In particular, compared with the case that $n_i$ is constant, the case that $n_i$ is increased with $i$ provides more accurate approximations $u_i$ and estimators $\eta(E_{i-1}, \varphi_i^{NN}) ~| | \varphi_i^{NN}| |_{L^2}$ and $\eta(E_{i-1}, \varphi_i^{NN})$ of the true errors, with negligible discrepancies beginning to show in each Galerkin iteration. Besides, for the case of adaptive choice of the number $n_i$, only no more than ten outer discontinuous Galerkin iterations may guarantee the convergence of our plane wave neural network algorithm.

Next, Table \ref{3dmaxcond_multidomain} shows the condition number of the matrix $K^{(j)}$ associated with the linear system (\ref{dgappro}) at each discontinuous Galerkin iteration $j$ when $n_1 =j+3$. We observe that the condition number is approximately unity as expected.

\vskip 0.1in
\begin{center}
       \tabcaption{Condition number for the discontinuous Galerkin matrix at each iteration $j$.}
\label{3dmaxcond_multidomain}
       \vskip 0.1in
\begin{tabular}{|c|c|c|c|c|c|c|c|} \hline
   $j$  & 1 & 2 &  3 &  4 &  5 & 6 &  7    \\ \hline
$\text{cond}(K^{(j)})$  & 1.00 & 1.00 &  1.00 &  1.02 &  1.05 & 1.11&  1.14
 \\ \hline
   \end{tabular}
     \end{center}

\vskip 0.1in

\begin{figure}[H]
\begin{center}
\begin{tabular}{cc}
\epsfxsize=0.4\textwidth\epsffile{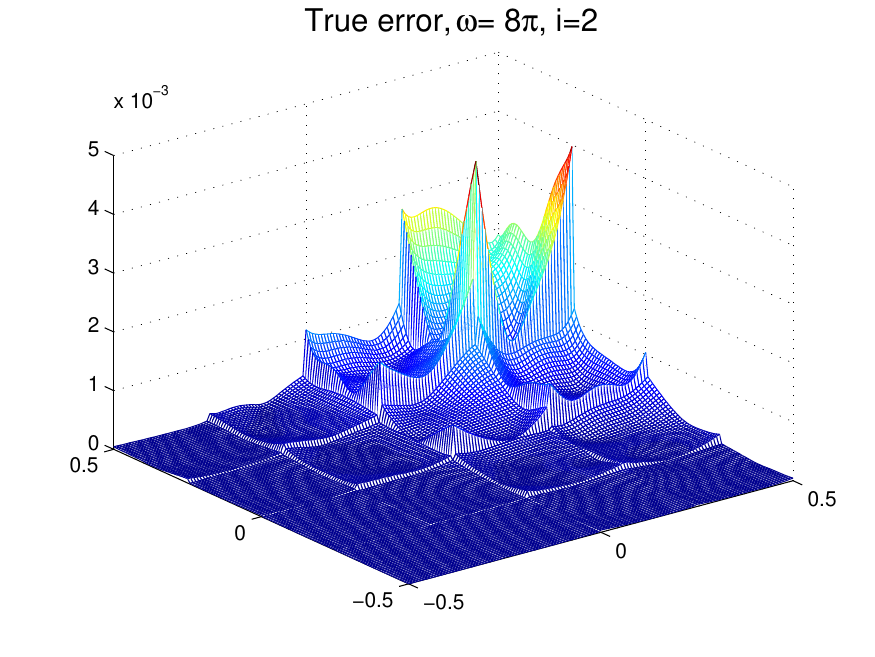}&
\epsfxsize=0.4\textwidth\epsffile{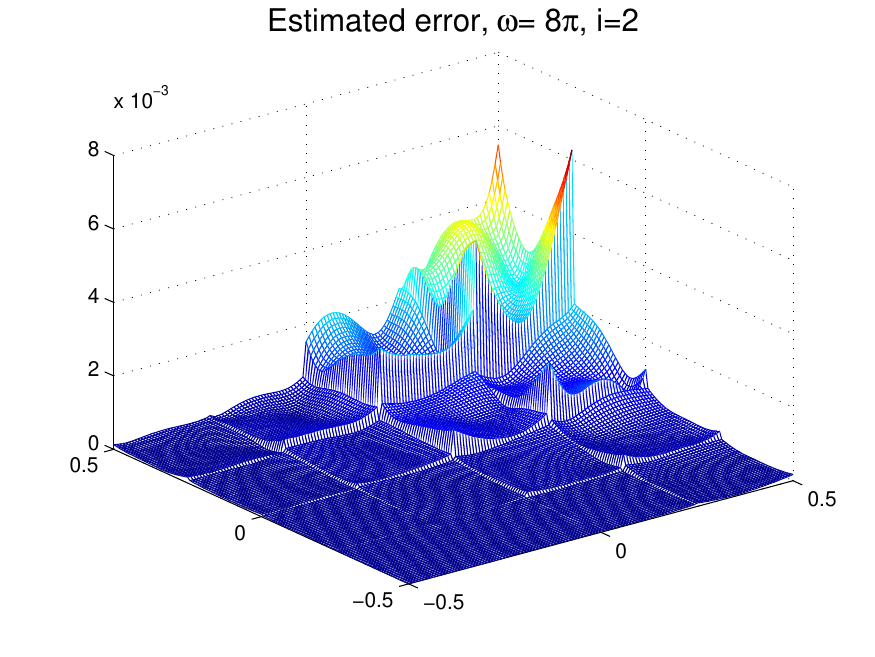}\\
\epsfxsize=0.4\textwidth\epsffile{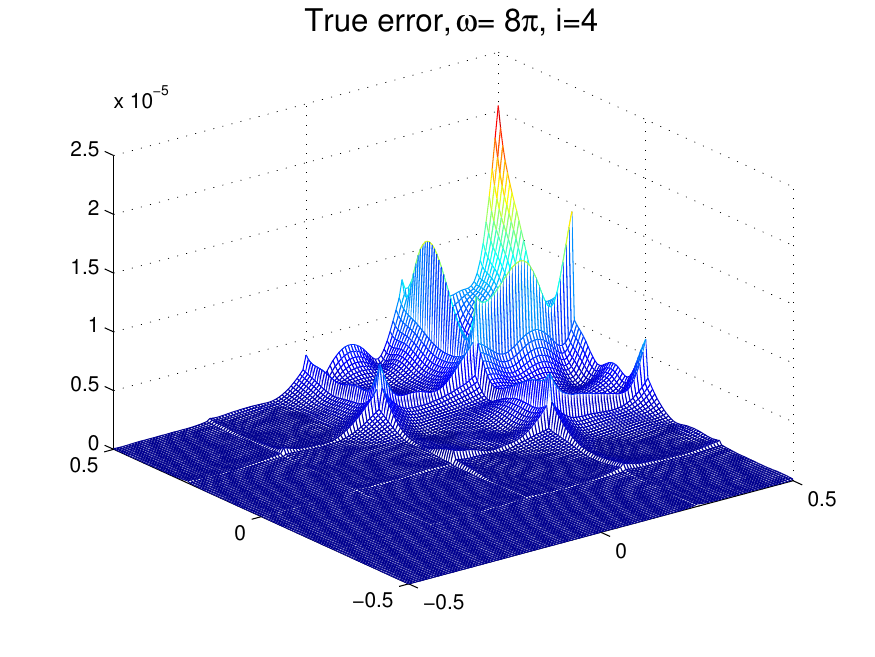}&
\epsfxsize=0.4\textwidth\epsffile{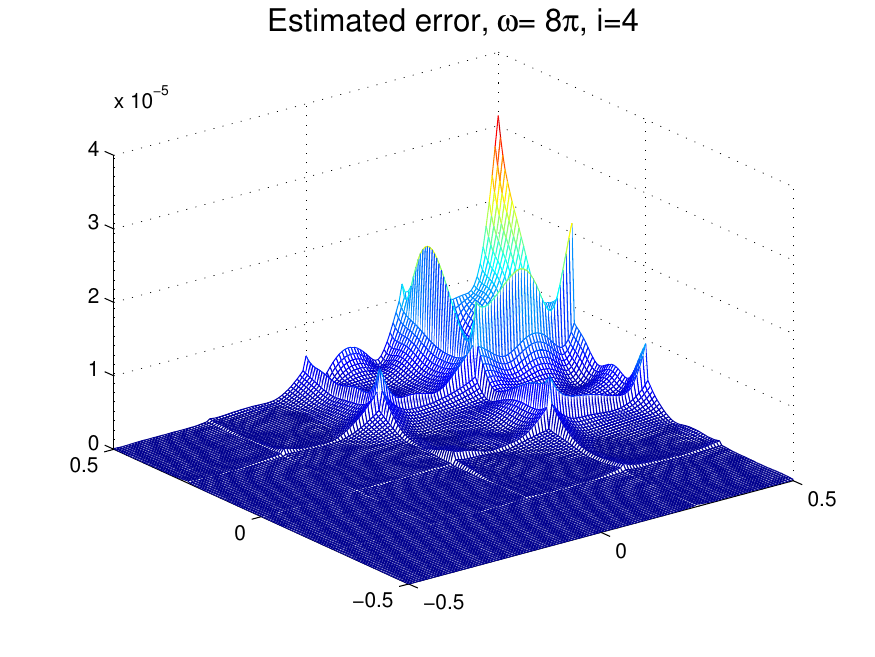}\\
\epsfxsize=0.4\textwidth\epsffile{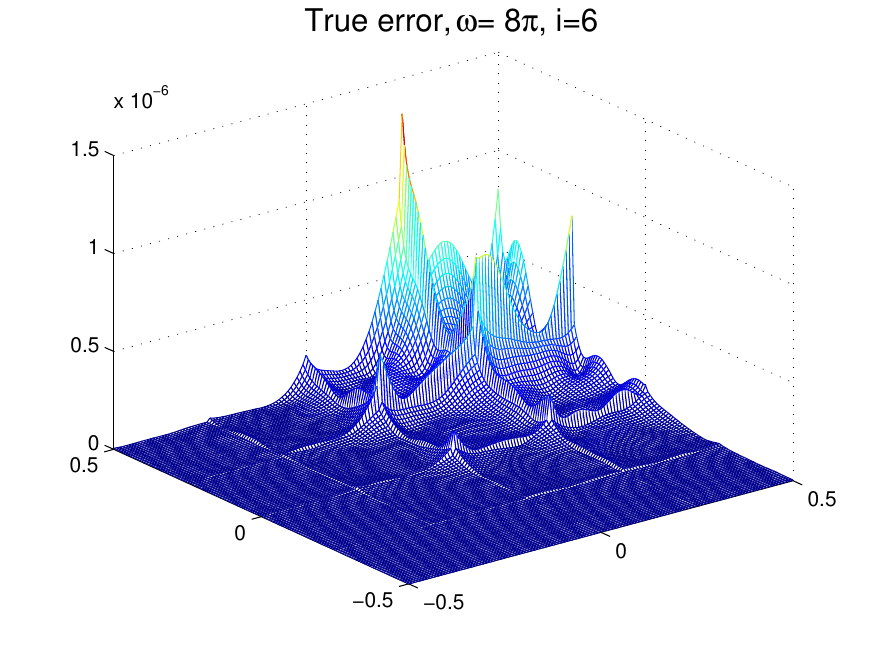}&
\epsfxsize=0.4\textwidth\epsffile{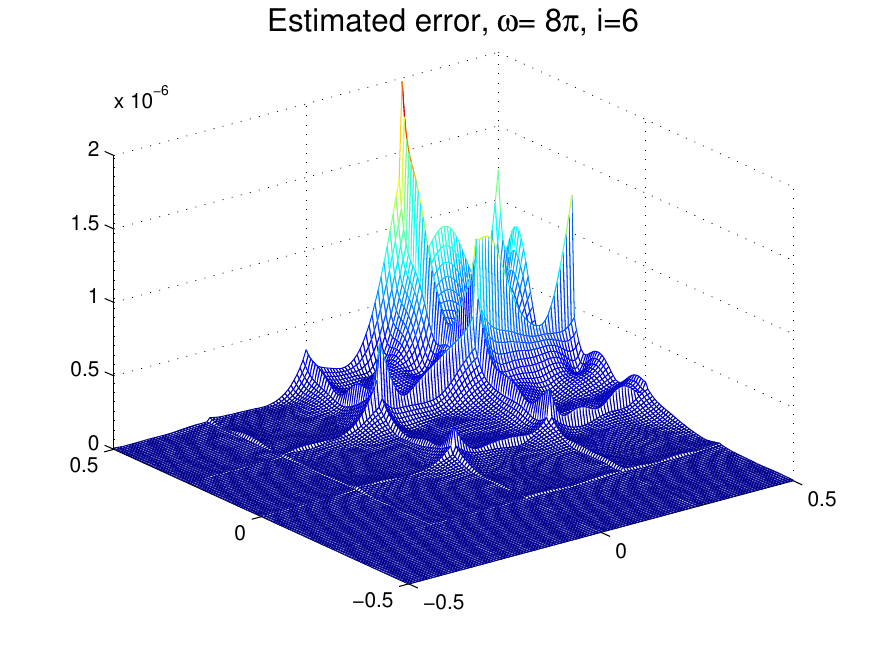}\\
\end{tabular}
\end{center}
 \caption{Displacement of a membrane (Maxwell's equations in three dimensions). Exact error $|E-E_{i-1}|$ (left column) and approximate error $\eta(E_{i-1}, \varphi_i^{NN}) ~|\varphi_i^{NN}|$ (right column) for i = 2, 4, 6. }
\label{3dmax2}
\end{figure}

Figure \ref{3dmax2} shows the exact error $|E-E_{i-1}|$ as well as the maximizer $\eta(E_{i-1}, \varphi_i^{NN}) ~|\varphi_i^{NN}|$ at several stages of the algorithm for the case of cross-section $z=0$.

It is apparent that the pointwise error is consistent with the Galerkin iteration error in Figure \ref{3dmaxwell}. Besides, we can found that, the closer to the singularity point ${x}_0=(0.6,0.6,0.6)$, the larger the pointwise error. To improve this situation, the local refinement on the meshwidth $h$ and the degrees of freedom $n_i$ per each element can be employed in the future.

\section{Summary} 
We have introduced the {\it discontinuous} Galerkin plane wave neural networks for the discretization of Helmholtz and Maxwell's equations. The sequence of discretized {\it discontinuous} subspaces can be adaptively generated by searching an approximate maximizer $\varphi_i^{NN}$ of the residual functional $r(u_{i-1}): V({\cal T}_h) \rightarrow \mathbb{R}$ expressed as the realization of a {\it discontinuous} plane wave neural network $V_{n_i}^{\sigma}$ with a single hidden layer. The unconditional convergence of the DGPWNN method are established only needing to select suitably distributed propagation directions, and also without the assumption on the boundedness of the neural network parameters. We have demonstrated that the solution of a sequence of small {\it discontinuous} neural networks with $hp-$refinement can be selected as the best linear combination of the resulting corrections; plane wave activation functions are employed to enhance the accuracy and accelerate the convergence. Simultaneously, the training key is to optimize a set of propagation directions by Adam algorithm. We emphasize that the performance of the algorithms is independent of wavenumbers, by choosing $h\approx \mathcal{O}(\frac{\pi}{\omega})$ and gradually increasing the width $n_i$ of the network at each discontinuous Galerkin iteration, and only no more than ten outer discontinuous Galerkin iterations may guarantee the convergence of our plane wave neural network algorithm.

In the future, we shall focus on the extension of the proposed {\it discontinuous} Galerkin plane wave neural networks to the models of nonhomogeneous wave equations, anisotropic wave equations and time-dependent wave equations.

\vskip0.2in

\noindent
{\bf Data Availability:} All data generated or analysed during this study are available from the corresponding author on reasonable request.

\noindent
{\bf Acknowledgements:} The work was supported in part by the Science Challenge Project (Grant No. TZ2025007),
National Natural Science Foundation of China (Grant No. 12371373) 
and Shandong Provincial Natural Science Foundation (Grant No. ZR2024MA059).



\section*{Declarations}

\noindent
{\bf Conflict of Interest:} The authors declared that they have no conflict of interest.


\vskip0.2in

\begin{thebibliography}{99}   



\bibitem{AD}
Ainsworth M. and Dong J., Galerkin neural networks: A framework for approximating variational equations with error control, SIAM J. Sci.
Comput., 43(2021), pp. A2474-A2501.



\bibitem{ALZOU}
 Ammari H., Li B. and Zou J., Mathematical analysis of electromagnetic scattering by dielectric nanoparticles with high refractive indices, Trans. Amer. Math. Soc., 376(2023), pp. 39-90.



 \bibitem{BZZOU}
 Bao G., Zhang H. and Zou J., Unique determination of periodic polyhedral structures by scattered electromagnetic fields, Trans. Amer. Math. Soc., 363(2011), pp. 4527-4551.



\bibitem{berg}
 Berg J. and Nystr\"{o}m K., A unified deep artificial neural network approach to partial differential equations in complex geometries, Neurocomputing, 317(2018), pp. 28-41.



\bibitem{CY}
 Cybenko G., Approximation by superpositions of a sigmoidal function, Math. Control Signals Systems, 2(1989), pp. 303-314.

\bibitem{eyu}
E W. and Yu B., The deep Ritz method: A deep learning-based numerical algorithm for
solving variational problems, Commun. Math. Stat., 6(2018), pp. 1-12.



\bibitem{FO}
 Folland G., Real analysis: Modern techniques and their applications, John Wiley $\&$ Sons, New York, 1999.

\bibitem{Git}
 Gittelson C., Hiptmair R. and Perugia I., Plane wave discontinuous Galerkin methods:
Analysis of the h-version, ESAIM Math. Model. Numer. Anal., 43(2009), pp. 297-331.


\bibitem{hexu}
 He J., Li L., Xu J., and Zheng C., ReLU deep neural networks and linear finite elements, J.
Comput. Math., 38(2020), pp. 502-527.


\bibitem{ref21}
 Hiptmair R., Moiola A. and Perugia I., Plane wave discontinuous
Galerkin methods for the 2D Helmholtz equation: analysis of the
$p$-version. SIAM J. Numer. Anal., 49(2011), pp. 264-284.

\bibitem{pwdg}
 Hiptmair R., Moiola A. and Perugia I., Error analysis of
Trefftz-discontinuous Galerkin methods for the time-harmonic Maxwell
equations, Math. Comp., 82(2013), pp. 247-268.

\bibitem{HMPsur}
 Hiptmair R., Moiola A. and Perugia I., A survey of Trefftz methods for the Helmholtz equation,
 Lect. Notes Comput. Sci. Eng., 114(2016), pp. 237-278.

\bibitem{HO}
 Hornik K., Approximation capabilities of multilayer feedforward networks, Neural Netw., 4(1991), pp. 251-257.


\bibitem{hy}
 Hu Q. and Yuan L., A weighted variational formulation based on plane wave basis for discretization of Helmholtz equations, Int. J. Numer. Anal. Model., 11(2014), pp. 587-607.

\bibitem{hy2}
 Hu Q. and Yuan L., A plane wave least-squares method for
time-harmonic Maxwell's equations in absorbing media, SIAM J. Sci.
Comput., 36(2014), pp. A1911-A1936.

\bibitem{hy3}
 Hu Q. and Yuan L., A plane wave method combined with local spectral
elements for nonhomogeneous Helmholtz equation and time-harmonic
Maxwell equations, Adv. Comput. Math., 44(2018), pp. 245-275.

\bibitem{HGA}
 Huttunen T., Gamallo P. and Astley R., Comparison of two wave element methods for the
Helmholtz problem, Commun. Numer. Meth. Engng., 25(2009), pp. 35-52.

\bibitem{HKM}
 Huttunen T., Kaipio J. and Monk P., The perfectly matched layer for
the ultra weak variational formulation of the 3D Helmholtz equation,
Int. J. Numer. Meth. Engng., 61(2004), pp. 1072-1092.



\bibitem{Kharazmi}
 Kharazmi E., Zhang Z. and Karniadakis G., Variational physics-informed neural networks for solving partial differential equations, arXiv:1912.00873, 2019.


\bibitem{LLPS}
  Leshno M.,   Lin V.,   Pinkus A. and  Schocken S., Multilayer feedforward networks with a
nonpolynomial activation function can approximate any function, Neural Netw., 6(1993), pp. 861-867.

\bibitem{Moiola}
 Moiola A., Trefftz-discontinuous Galerkin Methods for Time-Harmonic Wave Problems, PhD Dissertation, ETH Zurich, Switzerland, 2012.

\bibitem{MP}
 Moiola A. and Perugia I., A space-time Trefftz discontinuous Galerkin method for the acoustic wave equation in first-order formulation, Numer. Math., 138(2018), pp. 389-435.


\bibitem{mhp}
 Moiola A., Hiptmair R. and Perugia I., Plane wave approximation of
homogeneous Helmholtz solutions, Z. Angew. Math. Phys., 62(2011),
pp. 809-837.


\bibitem{ref12}
 Monk P. and Wang D., A least-squares method for the helmholtz
equation, Comput. Methods Appl. Mech. Engrg., 175(1999), pp.
121-136.

\bibitem{Par2022} Parolin E., Huybrechs D. and Moiola A., Stable approximation of Helmholtz solutions by evanescent plane waves, ESAIM: M2AN, 57(2023), pp. 3499-3536.


\bibitem{peng}
 Peng J., Wang J. and Shu S., Adaptive BDDC algorithms for the system arising from plane wave discretization
of Helmholtz equations, Int. J. Numer. Methods Eng., 116(2018), pp. 683-707.

\bibitem{peng2}
 Peng J., Shu S., Wang J. and Zhong L., Adaptive-Multilevel BDDC algorithm for three-dimensional plane wave Helmholtz systems, J. Comput. Appl. Math., 381(2021): 113011.



\bibitem{PRV}
 Petersen P., Raslan M. and Voigtlaender F., Topological properties of the set of functions
generated by neural networks of fixed size, Found. Comput. Math., 21(2021), pp. 375-444.

\bibitem{raissi}
Raissi M., Perdikaris P. and Karniadakis G., Physics-informed neural networks: A deep
learning framework for solving forward and inverse problems involving nonlinear partial differential equations, J. Comput. Phys., 378(2019), pp. 686-707.




\bibitem{Karniadakis}
 Shin Y., Zhang Z. and Karniadakis G., Error estimates of residual minimization using neural networks for linear PDEs, arXiv:2010.08019, 2020.

\bibitem{chung}
 Yeung T., Cheung K., Chung E., Fu S. and Qian J., Learning rays via deep neural network in a ray-based IPDG method for high-frequency Helmholtz equations in inhomogeneous media, J. Comput. Phys., 465(2022): 111380.

\bibitem{yuan2}
 Yuan L. and Hu Q., A PWDG method for the Maxwell system in anisotropic media with piecewise constant coefficient matrix, ESAIM: M2AN, 58(2024), pp. 1-22.

\bibitem{yuanyue}
 Yuan L., Wang X. and Yue X., A space-time Trefftz dG method for the second order time-dependent Maxwell system in anisotropic media, J. Comput. Math., 44(2026), pp. 394-426.



  \bibitem{yuanhu}
 Yuan L. and Hu Q., A discontinuous plane wave neural network method for Helmholtz equation and time-harmonic Maxwell's equations, Adv. Comput. Math., 51(2025): 18.



\bibitem{zang}
 Zang Y., Bao G., Ye X. and Zhou H., Weak adversarial networks for high-dimensional
partial differential equations, J. Comput. Phys., 411(2020): 109409.

\bibitem{luo}
 Zhang Y., Chen C., Shi N., Sun R. and Luo Z., Adam can converge without any modification on update rules, in Proceedings of the 36th International Conference on Neural Information Processing Systems (NIPS'22),
 pp. 28386-28399.

\bibitem{zhao}
 Zhao M., Zhu N. and Wang L., The electromagnetic scattering from multiple arbitrarily shaped cavities with inhomogeneous anisotropic media, J. Comput. Phys., 489(2023): 112274.


\end{thebibliography}
\end{document}